\input amstex
\documentstyle{amsppt}
\magnification=1200
\catcode`\@=11
\redefine\logo@{}
\catcode`\@=13
\pageheight{19cm}

\define \bn{\Bbb N}
\define \bz{\Bbb Z}
\define \bq{\Bbb Q}
\define \br{\Bbb R}
\define \bc{\Bbb C}

\define\Da{{\Cal D}}

\define\Ka{{\Cal K}}
\define \E{{\Cal E}}

\define\rk{\text{rk}~}


\define\Aut{\text{Aut}}

\define\Exc{\text{Exc}}
\define\0o{{\overline 0}}
\define\1o{{\overline 1}}

\define\wH{\widetilde{H}}
\define\wth{\widetilde{h}}

\TagsOnRight


\topmatter

\title
On correspondences of a K3 surface with itself. II
\endtitle

\author
Viacheslav V. Nikulin
\endauthor

\dedicatory To 60th Birthday of Igor Dolgachev
\enddedicatory

\address
Deptm. of Pure Mathem. The University of Liverpool, Liverpool
L69 3BX, UK;
\vskip1pt
Steklov Mathematical Institute,
ul. Gubkina 8, Moscow 117966, GSP-1, Russia
\endaddress
\email
vnikulin\@liv.ac.uk\ \
vvnikulin\@list.ru
\endemail

\abstract
Let $X$ be a K3 surface with a polarization $H$ of degree
$H^2=2rs$, $r,\,s\ge 1$, and the isotropic Mukai
vector $v=(r,H,s)$ is primitive. The moduli space of sheaves over $X$ with
the Mukai vector $v=(r,H,s)$ is again a K3 surface, $Y$.

We prove that $Y\cong X$ if the Picard lattice $N(X)$ has an element
$h_1$ with $(h_1)^2=f(v)$ and minor additional congruence conditions
modulo $N_i(v)$. All these conditions are exactly written, very efficient,
and they are necessary if $X$ is general with $\rk N(X)\le 2$.

Existence of such kind a criterion is very surprising, and it also
gives some geometric interpretation of elements in $N(X)$ with
a negative square.
Moreover, we describe all irreducible divisorial conditions on
moduli of $(X,H)$ which imply $Y\cong X$, and we prove that their number
is always infinite.

Thus, we treat in general problems considered in
\cite{MN1}, \cite{MN2} and \cite{N4}, where the additional condition
$H\cdot N(X)=\bz$ had been imposed.
\endabstract

\rightheadtext
{On correspondences of K3 with itself}
\leftheadtext{V.V. Nikulin}
\endtopmatter

\document

\head
0. Introduction
\endhead

Let $X$ be a K3 surface with a polarization $H$ of degree
$H^2=2rs$ where $r,\,s\,\in \bn$. Assume that the
isotropic Mukai vector $v=(r,H,s)$ is primitive.

Let $Y$ be the moduli space of sheaves (coherent and semi-stable
with respect to $H$) over $X$ with the isotropic
Mukai vector $v=(r,H,s)$.
The $Y$ (or, in special cases, its minimal resolution of singularities
which we denote by the same letter $Y$)
is again a K3 surface which is equipped with
a natural $nef$ element $h$ with $h^2=2ab$ where we denote
$c=\text{g.c.d}(r,s)$ and $a=r/c$, $b=s/c$ (see Sect. 2.1 below).
The surface $Y$ is isogenous to $X$ in the sense of Mukai.
The second Chern class of the
corresponding quasi-universal sheave gives then a 2-dimensional algebraic
cycle $Z\subset X\times Y$ and an algebraic correspondence between
$X$ and $Y$. See Mukai \cite{Mu1}---\cite{Mu5} and also Abe \cite{A} about 
these results.

Let $H$ be divisible by $d\in \bn$ where $\wH=H/d$ is primitive
in the Picard lattice $N(X)$ of $X$.
Primitivity of $v=(r,H,s)$ means that 
$\text{g.c.d}(r,d,s)=\text{g.c.d}(c,d)=1$. 
We have $d^2|ab$. Let $\gamma=\gamma (\wH)$ is defined by
$\wH\cdot N(X)=\gamma \bz$, i.e.
$H\cdot N(X)=\gamma d$. Clearly, $\gamma |(2rs/d^2)=\wH^2$.

We denote
$$
n(v)=\text{g.c.d}(r,s,d\gamma ).
$$
By Mukai, \cite{Mu2}, \cite{Mu3}, $T(X)\subset T(Y)$, and
$n(v)=[T(Y):T(X)]$ where $T(X)$ and $T(Y)$ are transcendental lattices
of $X$ and $Y$. We assume that
$$
n(v)=\text{g.c.d}(r,s,d\gamma)=\text{g.c.d}(c,d\gamma)=1.
\tag{0.1}
$$
Since $\text{g.c.d}(c,d)=1$, this is equivalent to $\text{g.c.d}(c,\gamma)=1$.
By Mukai \cite{Mu2}, \cite{Mu3}, the transcendental
periods
$(T(X), H^{2,0}(X))$ and $(T(Y),$ $H^{2,0}(Y))$ are
isomorphic in this case. We can expect that sometimes the surfaces $X$
and $Y$ are also isomorphic, and we then get a cycle $Z\subset X\times X$,
and a correspondence of $X$ with itself.  Thus, an interesting for
us question is

\proclaim{Question 1} When is $Y$ isomorphic to $X$?
\endproclaim

We want to answer this question in terms of Picard lattices
$N(X)$ and $N(Y)$ of
$X$ and $Y$. Then our question can be reformulated as follows:

\proclaim{Question 2} Assume that $N$ is a hyperbolic lattice,
$H_1\in N$ an element with square $2rs$. What are
conditions on $N$ and $H_1$ such that for any K3 surface
$X$ with Picard lattice $N(X)$ and s polarization $H\in N(X)$
the corresponding K3 surface $Y$ is isomorphic to $X$,
if the the pairs $(N(X), H)$ and $(N, H_1)$ are
isomorphic as abstract lattices with fixed elements?

In other words, what are conditions on $(N(X), H)$ as an abstract lattice
with an element $H$ which are sufficient for
$Y$ to be isomorphic to $X$, and they are necessary,
if $X$ is a general K3 surface with the Picard lattice $N(X)$?
\endproclaim

We answered this question in \cite{MN1}, \cite{MN2} and \cite{N4}
under the condition $d=\gamma=1$ (equivalently, $H\cdot
N(X)=\bz$): in \cite{MN1} for $r=s=2$; in \cite{MN2} for $r=s$; in
\cite{N4} for arbitrary $r$ and $s$.

The main surprising result of \cite{MN1}, \cite{MN2} and \cite{N4}
was that $Y\cong X$ if the Picard lattice $N(X)$ has an element
$h_1$ with some prescribed square $(h_1)^2$ and some minor
additional conditions. Moreover, these conditions are necessary to
have $Y\cong X$ for a general K3 surface $X$ with $\rho (X)=\rk N(X)=2$.
Thus, sometimes, elements of Picard lattice $N(X)$ deliver
important 2-dimensional algebraic cycles on $X\times X$. Moreover,
here $(h_1)^2$ can be negative, and this gives geometric meaning
for elements of the Picard lattice with negative square (it is
well-known only for $h_1^2=-2$; then $\pm h_1$ is effective).

Here we prove similar results in general.

We assume \thetag{0.1}. Then $d^2\vert ab$ and $\gamma \vert 2ab/d^2$.
Moreover, $\text{g.c.d}(a,b)=1$. We have $d=d_ad_b$
where $d_a=\text{g.c.d}(d,a)$ and $d_b=\text{g.c.d}(d,b)$. We define
$$
a_1={a\over d_a^2},\ b_1={b\over d_b^2}.
$$
We have $\gamma=\gamma_2\gamma_a\gamma_b$ where
$\gamma_a=\text{g.c.d}(a_1,\gamma)$,
$\gamma_b=\text{g.c.d}(b_1,\gamma)$, $\gamma_2=
\gamma/(\gamma_a\gamma_b)\,|\,2$.
We define
$$
a_2={a_1\over \gamma_a},\ b_2={b_1\over \gamma_b},\ e_2={2\over \gamma_2}.
$$

Due to Mukai \cite{Mu3} (see also Examples 2.3.4 and 2.3.5 below),
one has the following result for $\rho(X)=1$.
{\it If $\rho (X)=1$ and $X$ is a general K3 surface with its Picard lattice
(i. e. $\Aut (T(X), H^{2,0}(X))=\pm 1$),
then $Y\cong X$ if and only if $c=1$ and either $a_1=1$ or $b_1=1$.
In particular (e. g. see Lemma 2.1.1 below),
for a primitive Mukai vector $(r,H,s)$
one always has $Y\cong X$ if and only if $c=1$ and either $a_1=1$ or $b_1=1$.}

We prove (see Theorem 4.4) the following our main result for
$\rho(X)=2$. We denote as $\bz f(\wH)$ the orthogonal complement
to $\wH$ in the 2-dimensional lattice $N$ and use the invariants
$\gamma,\,\delta\in \bn$ and $\mu\in
\left(\bz/(2a_1b_1c^2/\gamma)\right)^\ast$ of the pair $\wH\in N$
(see Proposition 3.1.1). Here $\wH\cdot N=\gamma\bz$, $\det
N=-\gamma\delta$, $N=[\wH, f(\wH),
(\mu\wH+f(\wH))/(2a_1b_1c^2/\gamma)]$. One always has
$\delta\equiv \gamma \mu^2\mod {4a_1b_1c^2/\gamma}$. Moreover,
below $\wth_1=(p_1\wH+q_1f(\wH))/((2/\gamma_2)(a_1/\gamma_a)c)$
for $a$-series, and \newline
$\wth_1=(p_1\wH+q_1f(\wH))/((2/\gamma_2)(b_1/\gamma_b)c)$ for
$b$-series. Also we denote by $n^{(l)}$ the $l$-component of a
natural number $n$ for a prime $l$, i. e. $n^{(l)}=l^k\vert n$ and
\newline $\text{g.c.d}(n,n/n^{(l)})=1$.

\proclaim{Theorem 0.1} Let $X$ be a K3 surface and $H$ a
polarization of $X$ of degree $H^2=2rs$ where $r,\,s\in \bn$. Assume that
the Mukai vector $(r,H,s)$ is primitive. Let $Y$ be the moduli
space of sheaves on $X$ with the isotropic Mukai vector
$v=(r,H,s)$. Let $\wH=H/d$ be the corresponding primitive polarization.

We have $Y\cong X$ if there exists $\wth_1\in N(X)$ such that
$\wH$, $\wth_1$ belong to a 2-dimensional primitive sublattice
$N\subset N(X)$ such that $\wH\cdot N=\gamma \bz$, $\gamma >0$,
and
$$
\text{g.c.d}(c,d\gamma)=1,
$$
moreover, for one of $\epsilon=\pm 1$ the element $\wth_1$ belongs
to the $a$-series or to the $b$-series described below:

$\wth_1$ belongs to the $a$-series if
$$
\wth_1^2=\epsilon2b_1c\ and\ \wH\cdot \wth_1\equiv
0\mod{\gamma(b_1/\gamma_b)c},
$$
$$
\wH\cdot \wth_1\not\equiv 0\mod{\gamma(b_1/\gamma_b)cl_1},\
\wth_1/l_2\not\in N(X)
$$
for any prime $l_1$ such that $l_1^2\vert a_1$ and
$\text{g.c.d}(l_1,\gamma)=1$,
and any prime $l_2$ such that $l_2^2\vert b_1$ and
$\text{g.c.d}(l_2,\gamma)=1$,
and
$$p_1= {\wH\cdot
\wth_1\over \gamma (b_1/\gamma_b)c},\ \ q_1=-{f(\wH)\cdot \wth_1
\over \delta (b_1/\gamma_b)c}
$$
satisfy the singular condition (AS)
of $a$-series:
$$
\split &if\ odd\ prime\ l\vert \gamma\ and\ l^2\vert a_1,\
then\
q_1\not\equiv 0\mod l\ and\\
&either\ \delta\not\equiv 0\mod l\
or\ (\delta-\gamma \mu^2)\not\equiv 0\mod
(a_1^{(l)}/\gamma_a^{(l)})l;\\
&if\ odd\ prime\ l\vert \gamma\ and\ l\vert b_1,\ then\ q_1\equiv 0\mod
{\gamma_b^{(l)}};\\
&if\ odd\ prime\ l\vert \gamma\ and\ l^2\vert b_1,\  
then\ p_1\not\equiv 0\mod l;\\
&if\  2\vert \gamma,\ \gamma_2=1\ and\ 2\vert a_1,\ then\ p_1\equiv
1\mod 2;\\
&if\ 2\vert \gamma,\ \gamma_2=1\ and\ 4\vert a_1,\ then\
\delta-\gamma\mu^2\not\equiv 0\mod {(8a_1b_1c^2/\gamma)};\\
&if\ 2\vert \gamma,\ \gamma_2=1,\ and\ 2\vert b_1,\ then\ p_1-\mu
q_1\not\equiv 0\mod 4\ and\ q_1\equiv 0\mod {\gamma_b^{(2)}};\\
&if\ 2\vert \gamma,\ \gamma_2=2\ and\ 2\vert b_1\, ,\
then\ p_1\equiv 1\mod 2\ and\
q_1\equiv 0\mod \gamma^{(2)}/2.
\endsplit
$$

$\wth_1$ belongs to the $b$-series if
$$
\wth_1^2=\epsilon2a_1c\ and\ \wH\cdot \wth_1\equiv
0\mod{\gamma(a_1/\gamma_a)c},
$$
$$
\wH\cdot \wth_1\not\equiv 0\mod{\gamma(a_1/\gamma_a)cl_1},\
\wth_1/l_2\not\in N(X)
$$
for any prime $l_1$ such that
$l_1^2\vert b_1$ and $\text{g.c.d}(l_1,\gamma)=1$ and
any prime $l_2$ such that $l_2^2\vert a_1$ and
$\text{g.c.d}(l_2,\gamma)=1$,
and
$$p_1= {\wH\cdot
\wth_1\over \gamma (a_1/\gamma_a)c},\ \ q_1=-{f(\wH)\cdot
\wth_1\over \delta\,(a_1/\gamma_a)c}
$$
satisfy the singular condition (BS) of $b$-series:
$$
\split
&if\ odd\ prime\ l\vert \gamma\ and\ l\vert a_1,\ then\ q_1\equiv 0\mod
{\gamma_a^{(l)}};\\
&if\ odd\ prime\ l\vert \gamma\ and\ l^2\vert a_1,\  
then\ p_1\not\equiv 0\mod l;\\
&if\ odd\ prime\ l\vert \gamma\ and\ l^2\vert b_1,\
then\  q_1\not\equiv 0\mod l\ and\\
&either\ \delta\not\equiv 0\mod l\
or\ (\delta-\gamma \mu^2)\not\equiv 0\mod
(b_1^{(l)}/\gamma_b^{(l)})l;\\
&if\ 2\vert \gamma,\ \gamma_2=1,\ and\ 2\vert a_1,\ then\ p_1-\mu
q_1\not\equiv 0\mod 4\ and\ q_1\equiv 0\mod {\gamma_a^{(2)}};\\
&if\  2\vert \gamma,\ \gamma_2=1\ and\ 2\vert b_1,\ then\ p_1\equiv
1\mod 2;\\
&if\ 2\vert \gamma,\ \gamma_2=1\ and\ 4\vert b_1,\ then\
\delta-\gamma\mu^2\not\equiv 0\mod {(8a_1b_1c^2/\gamma)};\\
&if\ 2\vert\gamma,\ \gamma_2=2\ and\ 2\vert a_1\, ,\
then\ p_1\equiv 1\mod 2\ and\
q_1\equiv 0\mod \gamma^{(2)}/2.
\endsplit
$$

Moreover, one has formulae \thetag{4.23} and \thetag{4.24} in terms of $X$
for the canonical primitive $nef$ element $\wth$ of $Y$ defined by
$(-a,0,b)\mod \bz v$.

These conditions are necessary to have $Y\cong X$ if $\rho (X)\le
2$ and $X$ is a general K3 surface with its Picard lattice,
i. e. the automorphism group of the transcendental periods
$(T(X),H^{2,0}(X))$ is $\pm 1$.
\endproclaim

\smallpagebreak

As concrete examples, in Sect. 6 we specialize the theorem for
$\gamma=1$ and $\gamma=2$. The same can be done for any $\gamma$.

For the Mukai case when $c=1$ and either $a_1=1$ or $b_1=1$, one satisfies
conditions of Theorem 0.1 for $\wth_1=\wH$.

It seems many (if not all) known examples when $\rho(X)\ge 2$ and
$Y\cong X$ follow from this Theorem.
E.g. see \cite{C}, \cite{T1}---\cite{T3} and \cite{V}.

Like in \cite{MN1}, \cite{MN2} and \cite{N4} we also describe all
irreducible divisorial conditions on moduli of polarized K3
surfaces $(X,H)$ which imply $\wH\cdot N(X)=\gamma \bz$ and
$Y\cong X$. We show that they are labelled by pairs $(\pm
\mu,\delta)$ where $\pm \mu\in
\left(\bz/(2a_1b_1c^2/\gamma)\right)^\ast$, $\delta \in \bn$ and
$\delta\equiv \mu^2\gamma\mod 4a_1b_1c^2/\gamma$, moreover the
pair belongs to the $a$-series or to the $b$-series. It belongs to
the $a$-series if at
least for one $\epsilon=\pm 1$ the equation
$$
\gamma p_1^2-\delta q_1^2=\epsilon
2(2/\gamma_2)(a_1/\gamma_a)\gamma_b c
$$
has an integral solution $(p_1,q_1)$ where $(p_1,q_1)$ satisfy
conditions (A) of $a$-series \thetag{3.3.54}---\thetag{3.3.57}.
Similarly one can consider $b$-series changing $a$ and $b$ places.
See Sect. 4.

In Sect. 5, as an application, we prove that the number of
the irreducible divisorial conditions is infinite
if non-empty.
If $\gamma=1$, the same considerations as for
$d=\gamma=1$ in \cite{N4} show
that for any type of a primitive isotropic Mukai vector $(r,H,s)$
the number of divisorial conditions on moduli of K3 which imply
that $Y\cong X$ and $\gamma=1$ is always non-empty and infinite. In particular,
for any type of a primitive isotropic Mukai vector the number of
divisorial conditions on moduli of K3 which imply that $Y\cong X$
is always non-empty and infinite.

This paper generalizes to the general case results of \cite{N4}
(see also \cite{MN1} and \cite{MN2}) where a particular case
$d=\gamma=1$ had been considered.

 As in \cite{MN1}, \cite{MN2} and \cite{N4},
the fundamental tools to get the results above is the Global
Torelli Theorem for K3 surfaces proved by Piatetsky-Shapiro and
Shafarevich in \cite{PS}, and results of Mukai
\cite{Mu2}, \cite{Mu3}.  By results of \cite{Mu2}, \cite{Mu3},
we can calculate periods of $Y$ using periods of $X$; comparing the
periods, by the Global Torelli Theorem for K3 surfaces \cite{PS},
we can find out if $Y$ is isomorphic to $X$.

These paper treats in general problems considered in \cite{MN1},
\cite{MN2} and \cite{N4} where the additional condition
$\gamma=d=1$ had been imposed. It makes results of this paper much
more complicated. For instance, in these paper
we don't consider the question of non-emptiness of the divisorial
conditions on moduli for $\gamma >1$. It is more difficult
in the general setting of this paper. We hope to consider this
problem later.

\head
1. Preliminary notations and results about lattices and K3 surfaces
\endhead

\subhead
1.1. Some notations about lattices
\endsubhead
We use notations and terminology from \cite{N2} about lattices,
their discriminant groups and forms. A {\it lattice} $L$ is a
non-degenerate integral symmetric bilinear form. I. e. $L$ is a free
$\bz$-module equipped with a symmetric pairing $x\cdot y\in \bz$ for
$x,\,y\in L$, and this pairing should be non-degenerate. We denote
$x^2=x\cdot x$. The {\it signature} of $L$ is the signature of the
corresponding real form $L\otimes \br$.
The lattice $L$ is called {\it even}
if $x^2$ is even for any $x\in L$. Otherwise, $L$ is called {\it odd}.
The {\it determinant} of $L$ is defined to be $\det L=\det(e_i\cdot e_j)$
where $\{e_i\}$ is some basis of $L$. The lattice $L$ is {\it unimodular}
if $\det L=\pm 1$.
The {\it dual lattice} of $L$ is
$L^\ast=Hom(L,\,\bz)\subset L\otimes \bq$. The
{\it discriminant group} of $L$ is $A_L=L^\ast/L$. It has the order
$|\det L|$. The group $A_L$ is equipped with the
{\it discriminant bilinear form} $b_L:A_L\times A_L\to \bq/\bz$
and the {\it discriminant quadratic form} $q_L:A_L\to \bq/2\bz$
if $L$ is even. To get this forms, one should extend the form of $L$ to
the form on the dual lattice $L^\ast$ with values in $\bq$.

For $x\in L$, we shall consider
the invariant $\gamma (x)\ge 0$ where
$$
x\cdot L=\gamma (x)\bz .
\tag{1.1.1}
$$
Clearly, $\gamma (x)|x^2$ if $x\not=0$.

We denote by $L(k)$ the lattice obtained from a lattice $L$ by
multiplication of the form of $L$ by $k\in \bq$.
The orthogonal sum of lattices $L_1$ and $L_2$ is denoted by
$L_1\oplus L_2$.
For a symmetric integral matrix
$A$, we denote by $\langle A \rangle$ a lattice which is given by
the matrix $A$ in some bases. We denote
$$
U=\left(
\matrix
0&1\\
1&0
\endmatrix
\right).
\tag{1.1.2}
$$
Any even unimodular lattice of the signature $(1,1)$ is isomorphic to
$U$.

An embedding $L_1\subset L_2$ of lattices is called {\it primitive}
if $L_2/L_1$ has no torsion.
We denote by $O(L)$, $O(b_L)$ and $O(q_L)$ the automorphism groups of
the corresponding forms.
Any $\delta\in L$ with $\delta^2=-2$ defines
a reflection $s_\delta\in O(L)$ which is given by the formula
$$
x\to x+(x\cdot \delta)\delta,
$$
$x\in L$. All such reflections generate
the {\it 2-reflection group} $W^{(-2)}(L)\subset O(L)$.

\subhead
1.2. Some notations about K3 surfaces
\endsubhead
Here we remind some basic notions and results about K3 surfaces,
e. g. see \cite{PS}, \cite{S-D}, \cite{Sh}.
A K3 surface $X$ is a non-singular projective algebraic surface over
$\bc$ such that its canonical class $K_X$ is zero and the irregularity
$q_X=0$. We denote by $N(X)$ the {\it Picard lattice} of $X$ which is
a hyperbolic lattice with the intersection pairing
$x\cdot y$ for $x,\,y\in N(S)$. Since the canonical class $K_X=0$,
the space $H^{2,0}(X)$ of 2-dimensional holomorphic differential
forms on $X$ has dimension one over $\bc$, and
$$
N(X)=\{x\in H^2(X,\bz)\ |\ x\cdot H^{2,0}(X)=0\}
\tag{1.2.1}
$$
where $H^2(X,\bz)$ with the intersection pairing is a
22-dimensional even unimodular lattice of signature
$(3,19)$. The orthogonal lattice $T(X)$ to $N(X)$ in $H^2(X,\bz)$ is called
the {\it transcendental lattice of $X$.} We have
$H^{2,0}(X)\subset T(X)\otimes \bc$. The pair $(T(X), H^{2,0}(X))$ is
called the {\it transcendental periods of $X$}.
The {\it Picard number} of $X$ is
$\rho(X)=\rk N(X)$. A non-zero element $x\in N(X)\otimes \br$ is
called {\it nef} if $x\not=0$ and $x\cdot C\ge 0$ for any effective
curve $C\subset X$. It is known that an element $x\in N(X)$ is ample (i. e. it
defines a polarization)
if $x^2>0$, $x$ is $nef$, and the orthogonal complement
$x^\perp$ to $x$ in  $N(X)$ has no elements with square $-2$.
For any non-zero element $x\in N(X)$ with $x^2\ge 0$, there exists a reflection
$w\in W^{(-2)}(N(X))$ such that the element
$\pm w(x)$ is nef; it then is ample if $x^2>0$ and $x^\perp$ had no
elements with square $-2$ in $N(X)$. The $nef$ element $\pm w(x)$ is defined
canonically by $x$. It is called {\it the canonical $nef$ element of $x$}.

We denote by $V^+(X)$ the light cone of $X$, which is the half-cone of
$$
V(X)=\{x\in N(X)\otimes \br\ |\ x^2>0\ \}
\tag{1.2.2}
$$
containing a polarization of $X$. In particular, all $nef$ elements
$x$ of $X$ belong to $\overline{V^+(X)}$:
one has $x\cdot V^+(X)>0$ for them.

The reflection group $W^{(-2)}(N(X))$ acts in $V^+(X)$ discretely,
and its  fundamental
chamber is the closure $\overline{\Ka(X)}$ of the K\"ahler cone
$\Ka(X)$ of $X$. It is the same as the set of all $nef$ elements of $X$.
Its faces are orthogonal to the set $\Exc(X)$ of all exceptional curves $r$
on $X$ which are non-singular rational curves $r$ on $X$ with $r^2=-2$.
Thus, we have
$$
\overline{\Ka (X)}=\{0\not=x\in \overline{V^+(X)}\ |
\ x\cdot \Exc(X)\ge 0\,\}.
\tag{1.2.3}
$$
\par\bigskip

\head
2. Condition of $Y\cong X$ for a general K3 surface $X$
with a given Picard lattice
\endhead

\subhead
2.1. The correspondence
\endsubhead
Let $X$ be a smooth complex projective K3 surface with a
polarization
$H$ of degree $2rs$ where $r,\,s\in \bn$.

Assume that $H$ is divisible by $d\in \bn$ and $\wH=H/d$ is primitive
in $N(X)$.
Then $\wH^2=2rs/d^2$ and $d^2|rs$. We denote
$$
c=\text{g.c.d}(r,s),\ a=r/c,\ b=s/c.
\tag{2.1.1}
$$
We assume that the Mukai vector $(r,H,s)$ is primitive, i. e.
$$
\text{g.c.d}(r,s,d)=\text{g.c.d}(c,d)=1.
\tag{2.1.2}
$$

Let $Y$ be the moduli space of sheaves $\E$ (coherent and semi-stable
with respect to $H$) on $X$ with the primitive isotropic Mukai vector
$v=(r,H,s)$. Then $\rk \E=r$, $\chi (\E)=r+s$ and $c_1(\E)=H$.
The $Y$ (or, in special cases, its minimal resolution of
singularities which we denote by the same letter $Y$) is again a K3 surface.
See \cite{Mu1}---\cite{Mu5} and also \cite{A} about these results.

Let
$$
H^\ast(X,\bz)=H^0(X,\bz)\oplus H^2(X,\bz)\oplus H^4(X,\bz)
\tag{2.1.3}
$$
be the full cohomology lattice of $X$ equipped with the Mukai pairing
$$
(u,v)=-(u_0\cdot v_2+u_2\cdot v_0)+u_1\cdot v_1
\tag{2.1.4}
$$
for $u_0,v_0\in H^0(X,\bz)$, $u_1,v_1\in H^2(X,\bz)$,
$u_2,v_2\in H^4(X,\bz)$. We naturally identify
$H^0(X,\bz)$ and $H^4(X,\bz)$ with $\bz$. Then the Mukai pairing is
$$
(u,v)=-(u_0v_2+u_2v_0)+u_1\cdot v_1.
\tag{2.1.5}
$$
The element
$$
v=(r,H,s)=(r,H,\chi-r)\in H^\ast(X,\bz)
\tag{2.1.6}
$$
is isotropic, i.e. $v^2=0$, since $H^2=2rs$.
In this case (for a primitive $v$),
Mukai \cite{Mu2}---\cite{Mu5} (see also Abe \cite{A}) showed that $Y$ is a
K3 surface, and one has the natural identification
$$
H^2(Y,\,\bz) \cong (v^\perp/\bz v)
\tag{2.1.7}
$$
which also gives the isomorphism of the Hodge structures of $X$ and $Y$,
i. e. $H^{2,0}(Y)$ will be identified with the image of $H^{2,0}(X)$.
The $Y$ has the canonical $nef$ element $h$
defined by $(-a,0,b)\mod \bz v$ with $h^2=2ab$ (see Sect. 1.2).

In particular, \thetag{2.1.7} gives the embedding
$$
T(X)\subset T(Y)
\tag{2.1.8}
$$
of the transcendental lattices of the index
$$
[T(Y):T(X)]=n(v)=\min {|v\cdot x|}
\tag{2.1.9}
$$
where $x\in H^0(X,\bz)\oplus N(X)\oplus H^4(X,\bz)$ and $v\cdot x\not=0$
(see \cite{Mu2}, \cite{Mu3}). In this paper, we are interested in the
case when $Y\cong X$. By \thetag{2.1.9}, it may happen if $n(v)=1$ only.

We can introduce the invariant $\gamma=\gamma (\widetilde{H})\in \bn$
which is defined by
$$
\wH\cdot N(X)=\gamma \bz,
\tag{2.1.10}
$$
equivalently, $H\cdot N(X)=\gamma d\bz$. Clearly, $\gamma|2rs/d^2=\wH^2$,
and
$$
n(v)=\text{g.c.d}(r,s,\gamma d)=\text{g.c.d}(c,\gamma d).
\tag{2.1.11}
$$
Thus, $n(v)=1$, and it is possible to have $Y\cong X$ only if
$$
\text{g.c.d}(r,s,\gamma d)=\text{g.c.d}(c,\gamma)=\text{g.c.d}(c,d)=1.
\tag{2.1.12}
$$
This is exactly the case when, according to Mukai,  the transcendental periods
$$
(T(X), H^{2,0}(X))\cong (T(Y),H^{2,0}(Y))
\tag{2.1.13}
$$
are isomorphic.

From \thetag{2.1.7}, we obtain the following {\it specialization principle.}

We say that a K3 surface $X$ is {\it general (for its Picard lattice)} if
the automorphism group of the transcendental periods 
$\Aut(T(X), H^{2,0}(X))=\pm 1$. 
We have

\proclaim{Lemma 2.1.1} (The specialization principle.)
Assume that for a general K3 surface $X$ with $N=N(X)$ and a primitive
isotropic
Mukai vector $v=(r,H,s)$ where $H\in N$ is a polarization of $X$,
one has $Y\cong X$

Then the same is valid for any K3 surface
$X^\prime$ such that $H\in N\subset N(X^\prime)$ if $N\subset N(X^\prime)$
is a primitive sublattice in $N(X^\prime)$ and $H$ is a polarization of $X$.
\endproclaim

\demo{Proof} Since $Y\cong X$ and $X$ is general, there exist only two
isomorphisms of transcendental periods 
$(T(X), H^{2,0}(X))\cong (T(Y),H^{2,0}(Y))$ 
which are $\pm 1$ of the identification of the transcendental periods
$(T(X), H^{2,0}(X))=(T(Y), H^{2,0}(Y))$ which is defined by
Mukai identification \thetag{2.1.7}. Since $Y\cong X$, 
there exists an extension 
of this identification to the isomorphism $H^2(X,\bz)\cong H^2(Y, \bz)$ of
the cohomology lattices.

Now assume that $H\in N\subset N(X^\prime)$.
Let $Y^\prime$ be the corresponding moduli space of sheaves on $X^\prime$ with
the same Mukai vector $v$.

By local epimorphicity of the period map for K3, we can assume that
$N=N(X)$ is the Picard lattice of a K3 surface $X$ with the polarization
$H$, $X$ is general and the embedding $N(X)\subset N(X^\prime)$ extends to
the identification of the cohomology lattices $H^2(X,\bz)=H^2(X^\prime, \bz)$.
Then $T(X)\supset T(X^\prime)$ is a primitive sublattice.
By the Mukai identification \thetag{2.1.7}, the identification \thetag{2.1.7}
$(T(X^\prime), H^{2,0}(X^\prime))=(T(Y^\prime), H^{2,0}(Y^\prime))$ 
is extending to 
the identification of the transcendental periods
$(T(X),H^{2,0}(X))=(T(Y),H^{2,0}(Y))$ and to the identification
of the cohomology lattices $H^2(Y,\bz)=H^2(Y^\prime,\bz)$.
Since $X$ is general, $Y\cong X$, and the
identification above of their transcendental periods
extends to an isomorphism of the lattices
$H^2(X,\bz)\cong H^2(Y,\bz)$.
This gives the isomorphism
$H^2(X^\prime,\bz)=H^2(X,\bz)\cong H^2(Y,\bz)=H^2(Y^\prime,\bz)$
which extends the above isomorphism
$(T(X^\prime),H^{2,0}(X^\prime))\cong (T(Y^\prime),H^{2,0}(Y^\prime))$.

By global Torelli Theorem for K3 surfaces \cite{PS}, 
this defines an isomorphism 
$Y^\prime\cong X^\prime$. This finishes the proof.
\enddemo

\subhead
2.2. The characteristic map of a primitive element of a lattice
\endsubhead
Let $S$ be an even lattice and $P\in S$ its primitive element
with $P^2=2m\not=0$ and $\gamma (P)=\gamma|2m$ (in $S$), i. e.
$P\cdot S=\gamma \bz$.

We want to calculate the discriminant quadratic form of $S$.
Consider
$$
K(P)=P^\perp_{S}
\tag{2.2.1}
$$
the orthogonal complement to $P$ in $S$.
Put $P^\ast=P/2m$. Then any element $x\in S$ can be written as
$$
x=n\gamma P^\ast+k^\ast
\tag{2.2.2}
$$
where $n\in \bz$ and $k^\ast\in K(P)^\ast$, because
$$
\bz P\oplus K(P)\subset S\subset
S^\ast\subset \bz P^\ast\oplus K(P)^\ast
$$
and $P\cdot S=\gamma \bz$. Since $\gamma(P)=\gamma$,
the map $n\gamma P^\ast+[P] \to k^\ast +K(P)$ gives an isomorphism of
the groups
$\bz/{2m\over \gamma} \cong [\gamma P^\ast]/[P]\cong
[u^\ast(P)+K(P)]/K(P)$ where
$u^\ast(P)+K(P)$ has order $2m/\gamma$ in $A_{K(P)}=K(P)^\ast/K(P)$.
Clarify in \cite{N2}. It follows,
$$
S=[\bz P, K(P), \gamma P^\ast+u^\ast(P)].
\tag{2.2.3}
$$
The element $u^\ast(P)$ is defined canonically mod $K(P)$ by the
condition that $\gamma P^\ast + u^\ast \in S$. The element
$$
u^\ast(P)+K(P)\in K(P)^\ast/K(P)
\tag{2.2.4}
$$
is called {\it the canonical element of $P$}.
Since
$\gamma P^\ast+u^\ast(P)$ belongs to the even lattice $S$, it follows
$$
(\gamma P^\ast+u^\ast(P))^2=\frac{\gamma^2}{2m}
           +{u^\ast(P)}^2 \equiv 0 \mod 2.
\tag{2.2.5}
$$
For $n\in \bz$ and $k^\ast \in K(P)^\ast$, we have
$x=nP^\ast+k^\ast \in S^*$ if and only if
$$
(nP^\ast+k^\ast)\cdot (\gamma P^\ast+u^\ast(P))={n\gamma\over 2m}+
k^\ast\cdot u^\ast(P)\in \bz.
$$
It follows,
$$
\aligned
S^\ast &=
\{nP^\ast+k^\ast\ |\ n\in \bz,\ k^\ast \in K(P)^\ast,\
n\equiv -{2m\over \gamma} \ (k^\ast\cdot u^\ast(P))
\mod{ {2m\over \gamma}} \}\subset \\
& \subset \bz P^\ast+K(P)^\ast.
\endaligned
\tag{2.2.6}
$$
It gives the calculation of the discriminant group $A_S=S^\ast /S$ where
$S$ is given by \thetag{2.2.3} and $S^\ast$ is given by \thetag{2.2.6}.

We define {\it the canonical submodule}
$\widetilde{K}(P)^\ast\subset \bz\oplus K(P)^\ast$ by the condition
$$
\widetilde{K}(P)^\ast=\{n\in \bz,\, k^\ast \in K(P)^\ast\ |\
n\equiv -{2m\over \gamma} \ (k^\ast\cdot u^\ast(P))
\mod{ {2m\over \gamma}}\}.
\tag{2.2.7}
$$
Now we define the {\it characteristic map}
$$
\kappa (P):\widetilde{K}(P)^\ast\to A_S,\ \
\tag{2.2.8}
$$
by the condition
$$
\kappa (P)(n,k^\ast)=nP^\ast+k^\ast\mod S.
\tag{2.2.9}
$$
Obviously, the characteristic map is epimorphic. Its kernel is
$$
\widetilde{K}(P)^\ast_0=[(2m\bz,\,K),\bz(\gamma,\,u^\ast(P))]\cong S.
\tag{2.2.10}
$$
Thus, we correspond to a primitive $P\in S$ with $P^2=2m$ and
$\gamma (P)=\gamma$ the canonical triplet
$$
(K(P),u^\ast(P)+K(P),\kappa(P)).
\tag{2.2.11}
$$
This triplet is important because of the trivial  but very important for
us

\proclaim{Lemma 2.2.1} Let $P_1\in S$ and $P_2\in S$ are two
primitive elements of an even lattice S with $P_1^2\not=0$ and
$P_2^2\not=0$.

There exists an automorphism $f\in O(S)$ such that
$f(P_1)=P_2$ and $f$ gives $\pm 1$ on the discriminant group
$A_S=S^\ast/S$ if and only if $P_1^2=P_2^2$, $\gamma (P_1)=\gamma (P_2)$,
and there exists an isomorphism of lattices
$\phi: K(P_1)\to K(P_2)$ such that
$\phi^\ast (u^\ast(P_2)+K(P_2))=u^\ast(P_1)+K(P_1)$ and
$\widetilde{\phi}^\ast=(\text{id},\phi^\ast):\widetilde{K}(P_2)^\ast \to
\widetilde{K}(P_1)^\ast$ is $\pm$ commuting  with the characteristic maps
$\kappa (P_1)$ and $\kappa (P_2)$, i. e.
$$
\kappa (P_1)\widetilde{\phi^\ast}=\pm \kappa (P_2).
$$
\endproclaim

\demo{Proof} Trivial.
\enddemo

We also mention that
$$
\det (S)=\gamma^2\det K(P)/2m.
\tag{2.2.12}
$$
because $[S:\bz P\oplus K(P)]=2m/\gamma$

\subhead
2.3. Relation between periods of $X$ and $Y$
\endsubhead
Here we consider the case and notations of Sect. 2.1. Thus,
for the primitive isotropic Mukai vector $v=(r,H,s)$, $r,\,s\ge 1$,
$H^2=2rs$,
we assume that for $d\in \bn$, the element
$\wH=H/d$ is primitive in $N(X)$,
$\gamma (\wH)=\gamma$ (in $N(X)$).
We remind that $c=\text{g.c.d}(r,s)$, $a=r/c$, $b=s/c$,  and
$n(v)=\text{g.c.d}(r,s,d\gamma)=(c,d\gamma)=1$. It follows,
$d^2|ab$ and $\gamma d^2|2ab$. Thus, our data are defined by
$$
a,\,b,\,c,\,d,\,\gamma\,\in \bn,\ \text{such that\ }\
(a,b)=(d,c)=(d,\gamma)=1,\ d^2|ab,\ \gamma d^2|2ab,
\tag{2.3.1}
$$
and by a primitive polarization
$$
\wH\in N(X)\ \ \text{such that\ }\ \wH^2=2abc^2/d^2,\ \gamma(\wH)=\gamma.
\tag{2.3.2}
$$
Then, the Mukai vector $v=(r,H,s)=(ac,d\wH,bc)$.

Let us denote by $e_1$ the canonical generator of $H^0(X,\bz)$ and
by $e_2$ the canonical generator of $H^4(X,\bz)$. They generate
the sublattice $U$ in $H^\ast (X,\bz)$ with the Gram matrix $U$.
Consider Mukai vector $v=(re_1+se_2+H)$. We have
$$
N(Y)=v\,^\perp _{U\oplus N(X)}/\bz v.
\tag{2.3.3}
$$
Let us calculate $N(Y)$. Let $K(H)=K(\wH)=(H)^\perp_{N(X)}$. We denote
$H^\ast=H/(2rs)\in (\bz H)^\ast=\bz H^\ast$, and
$\wH^\ast=\wH/(2rs/d^2)=dH^\ast\in (\bz \wH)^\ast=\bz \wH^\ast$.
 Then we have an
embedding of lattices of finite index
$$
\bz \wH\oplus K(H)\subset N(X)\subset N(X)^\ast\subset
\bz \wH^\ast \oplus K(H)^\ast
\tag{2.3.4}
$$
We have the orthogonal decomposition up to finite index
$$
U\oplus \bz \wH\oplus K(H)\subset U\oplus N(X)\subset
U\oplus \bz \wH^\ast \oplus K(H)^\ast .
\tag{2.3.5}
$$
Let $f=x_1e_1+x_2e_2+yH^\ast+z^\ast\in v^\perp_{U\oplus N(X)}$,
$z^\ast \in K(H)^\ast$.
Then $-sx_1-rx_2+y=0$ since $f\in v^\perp$ and hence $(f,v)=0$.
Thus, $y=(sx_1+rx_2)$ and
$$
f=x_1e_1+x_2e_2+(sx_1+rx_2)H^\ast + z^\ast .
\tag{2.3.6}
$$
Here
$f\in U\oplus N(X)$ if and only if
$$x_1, x_2 \in \bz,\ \ sx_1+rx_2\equiv 0\mod d,\ \
{sx_1+rx_2\over d}\wH^\ast+z^\ast\in N(X).
\tag{2.3.7}
$$
Equivalently, by Sect. 2.2, we have
$$
sx_1+rx_2\equiv 0\mod d\gamma,\ \
z^\ast={sx_1+rx_2\over d\gamma}u^\ast(H)\mod K(H).
\tag{2.3.8}
$$
We denote
$$
h^\prime=(-a,b)\oplus 0\in U\oplus N(X).
$$
Clearly, $h^\prime \in v^\perp$ and $h=h^\prime \mod \bz v\in N(Y)$.
Thus, the orthogonal complement
contains
$$
[\bz v, \bz h^\prime, K(H)]
\tag{2.3.9}
$$
where $h^\prime =-ae_1+be_2$,
and \thetag{2.3.9} is a sublattice of finite index
in $(v^\perp)_{U\oplus N(X)}$. The generators $v$, $h^\prime$
and generators of $K(H)$ are free,
and we can rewrite $f$ above using these generators
with rational coefficients. We have
$$
e_1={v-ch^\prime-H\over 2r},\ e_2={v+ch^\prime -H\over 2s}.
\tag{2.3.10}
$$
It follows,
$$
f={sx_1+rx_2\over 2rs}v+{c(-sx_1+rx_2)\over 2rs}h^\prime +z^\ast
\tag{2.3.11}
$$
where $x_1$, $x_2$, $z^\ast$ satisfy \thetag{2.3.8}. Considering
$\mod \bz v$, we finally get
$$
N(Y)={c(-sx_1+rx_2)\over 2rs}h+{sx_1+rx_2\over d\gamma}u^\ast(\wH)+K(H),
\ \text{where}\  sx_1+rx_2\equiv 0\mod d\gamma.
\tag{2.3.12}
$$

Let us calculate the lattice $K(h)=h^\perp_{N(Y)}$. It is equal to
$(sx_1+rx_2)/(d\gamma)u^\ast(\wH)+K(H)$ where
$sx_1+rx_2\equiv 0\mod {d\gamma}$ and $-sx_1+rx_2=0$.
It follows, $x_1=ax$, $x_2=bx$, $x\in \bz$, and
$sx_1+rx_2=2abcx\equiv 0\mod d\gamma$, which is always true. Thus,
$K(h)=[2abc/(d\gamma)u^\ast(\wH)+K(H)]$. By Sect. 2.2,
$u^\ast(\wH)+K(H)$ has the order $2abc^2/(d^2\gamma)$.
We have $(2abc/(d\gamma),2abc^2/(d^2\gamma)=
2abc/(d^2\gamma)(d,c)$ where $2abc/(d^2\gamma)\in \bn$ and $(d,c)=1$.
It follows, $K(h)=[K(H),2abc/(d^2\gamma)u^\ast(\wH)]$, and the index
$[K(h):K(H)]=c$.

Let us show that $d|h$ in $N(Y)$ and $h/d$ is primitive in $N(Y)$.
By \thetag{2.3.12}, the primitive submodule in $N(Y)$, generated by $h$,
is $c(-sx_1+rx_2)/(2rs)h$ where $(sx_1+rx_2)/(d\gamma)u^\ast(\wH)\in K(h)$
and $sx_1+rx_2\equiv 0\mod d\gamma$. From calculation of $K(h)$ above,
we then get $(sx_1+rx_2)/(d\gamma)\equiv 0\mod{2abc/(d^2\gamma)}$. It
follows, $bx_1+ax_2\equiv 0\mod 2ab/d$. Let $d=d_ad_b$ where $d_a|a$ and
$d_b|b$. Then $(a/d_a)|x_1$ and $(b/d_b)|x_2$. It follows,
$x_1=(a/d_a)\widetilde{x}_1$ and $x_2=(b/d_b)\widetilde{x}_2$ where
$\widetilde{x}_1,\ \widetilde{x}_2 \in \bz$ and
$d_b\widetilde{x}_1+d_a\widetilde{x}_2\equiv 0\mod 2$. It follows that
the module $c(-sx_1+rx_2)/(2rs)h=
(-d_b\widetilde{x}_1+d_a\widetilde{x}_2)/(2d_ad_b)$ where
$d_b\widetilde{x}_1+d_a\widetilde{x}_2\equiv 0\mod 2$. It follows
that this module is $\bz(h/d)$. It proves the statement.
We denote
$$
\wth={h\over d}\ .
\tag{2.3.13}
$$
We had proved that $\wth$ is primitive in $N(Y)$ and $h=d\wth$.
We have $\wth^2=2ab/d^2$.

Let us show that $\gamma (\wth)=\gamma (H)=\gamma$. We remind that
$\gamma(\wth)\bz = \wth\cdot N(Y)$. By Sect. 2.2 and \thetag{2.3.12},
$\wth^2/\gamma(\wth)$ is equal to the index
$[[(sx_1+rx_2)/(d\gamma)u^\ast(\wH)+K(H)]:K(h)]$
where $sx_1+rx_2\equiv 0\mod d\gamma$. Such elements $x_1$, $x_2$
give $cd\gamma \bz$. Thus,
$[(sx_1+rx_2)/(d\gamma)u^\ast(\wH)+K(H)]=[cu^\ast(\wH)+K(H)]$. We had
proved that $K(h)=[2abc/(d^2\gamma)u^\ast(\wH)+K(H)]$. Thus,
the index is equal to $2ab/d^2\gamma$. Then we get
$\wth^2/\gamma(\wth)=2ab/(d^2\gamma (\wth))=2ab/(d^2\gamma)$.
It follows $\gamma (\wth)=\gamma$. We had also proved that
$u^\ast(\wth)=mcu^\ast(H)+K(h)$ where $m$ is defined
$\mod 2ab/(d^2\gamma)$. We shall calculate $m$ below.

Now we can rewrite \thetag{2.3.12} in the form \thetag{2.2.2}.
We denote $\wth^\ast = \wth /\wth^2=\wth/(2ab/d^2)$. We have
$$
N(Y)={-bx_1+ax_2\over d}\wth^\ast+{(bx_1+ax_2)\over d\gamma}cu^\ast(\wH)+K(H)
\tag{2.3.14}
$$
where $bx_1+ax_2\equiv 0\mod d\gamma$. We have proved that
the elements  $(-bx_1+ax_2)/d$ give $\bz \gamma$.

Let us write $d=d_ad_b$ where $d_a|a$ and $d_b|b$. Since $\gamma|2ab/d^2$
and $\text{g.c.d}(a,b)=1$, we can write $\gamma=\gamma_2\gamma_a\gamma_b$
where $\gamma_a=\text{g.c.d}(\gamma,a/d_a^2)$,
$\gamma_b=\text{g.c.d}(\gamma,b/d_b^2)$ and
$\gamma_2=\gamma/(\gamma_a\gamma_b)$. Clearly, $\gamma_2|2$.

We define
$m=m(a,b,d,\gamma)\mod 2ab/(d^2\gamma)$ by the conditions
$$
m\equiv -1\mod 2a/(d_a^2\gamma_a\gamma_2)\ \ \text{and}\ \
m\equiv 1\mod 2b/(d_b^2\gamma_b\gamma_2).
\tag{2.3.15}
$$
The congruences \thetag{2.3.15} define $m$ uniquely $\mod 2ab/(d^2\gamma)$.
Really, assume that $x\equiv 0\mod 2b/(d_b^2\gamma_b\gamma_2)$ and
$x\equiv 0 \mod 2a/(d_a^2\gamma_a\gamma_2)$. Then
$x/(2/\gamma_2)\equiv 0$ $\mod b/(d_b^2\gamma_b)$ and
$x/(2/\gamma_2)\equiv 0\mod a/(d_a^2\gamma_a)$. It follows that
$x/(2/\gamma_2)\equiv 0\mod ab/(d^2\gamma_a\gamma_b)$. Thus,
$x\equiv 0\mod 2ab/(d^2\gamma_a\gamma_b\gamma_2)$ where
$2ab/(d^2\gamma_a\gamma_b\gamma_2)=2ab/(d^2\gamma)$.
Thus $x\equiv 0\mod 2ab/(d^2\gamma)$. Clearly, there exists a
unique
$m(a,b)\mod 2ab$ defined by the condition
$$
m(a,b)\equiv -1\mod 2a\ \text{and}\ m(a,b)\equiv 1\mod 2b.
\tag{2.3.16}
$$
Then
$$
m(a,b,d,\gamma)\equiv m(a,b)\mod {{2ab\over d^2\gamma}},
\tag{2.3.17}
$$
and one can take \thetag{2.3.16} and \thetag{2.3.17} as definition
of $m(a,b,d,\gamma)$.

Let us prove that $u^\ast(\wth)+K(h)=m(a,b,d,\gamma)cu^\ast(\wH)+K(h)$.
We had proved that $u^\ast (\wth)+K(h)=mcu^\ast (\wH)+K(h)$
where $m$ is defined $\mod 2ab/(d^2\gamma)$.
To find $m \mod $ $2ab/(d^2\gamma)$, one should put
$(-bx_1+ax_2)/d=\gamma$ in \thetag{2.3.14} or
$$
-bx_1+ax_2=d\gamma.
\tag{2.3.18}
$$
From \thetag{2.3.18} we have $d_a\gamma_a|x_1$ and $d_b\gamma_b|x_2$.
From \thetag{2.3.14} and \thetag{2.3.18}, we get
$$
m\equiv {ax_2+bx_1\over d\gamma}
\equiv 1+{2bx_1\over d\gamma}\mod 2ab/(d^2\gamma),
\tag{2.3.19}
$$
and
$$
m\equiv -1+{2ax_2\over d\gamma} \mod 2ab/(d^2\gamma).
\tag{2.3.20}
$$
Since $d_a\gamma_a|x_1$, we get from \thetag{2.3.19} that
$m\equiv 1\mod 2b/(d_b^2\gamma_b\gamma_2)$, and from \thetag{2.3.20} that
$m\equiv -1\mod 2a/(d_a^2\gamma_a\gamma_2)$.
Thus, $m=m(a,b,d,\gamma )\mod 2ab/(d^2\gamma)$.
It proves the statement.

Since $h^2=2ab$ and $H^2=2abc^2$, we can formally put $h=H/c$.
By our construction, $K(H)\subset K(h)$ is a sublattice. Thus, we
can consider $N(X)$ and $N(Y)$ as extensions of finite index of
a common sublattice $\bz \wH+K(H)$.

Finally, we get the very important for us

\proclaim{Proposition 2.3.1} The Picard lattice of $X$ is
$$
N(X)=[\wH,K(H),\gamma \wH^\ast+u^\ast(\wH)]
\tag{2.3.21}
$$
where $\wH=H/d$ is primitive in $N(X)$ with $H^2=2abc^2$
and $\wH^\ast=d^2\wH/(2abc^2)$, the lattice $K(H)=H^\perp_{N(X)}$,
and  $u^\ast(\wH)+K(H)$ has the order $2abc^2/(d^2\gamma)$ in
\newline
$K(H)^\ast/K(H)$.

The Picard lattice of $Y$ is
$$
N(Y)=[\wth= h/d,K(h),\gamma \wth^\ast+u^\ast (\wth)],
\tag{2.3.22}
$$
where the element $h=(-a,0,b)\mod \bz v$, $h^2=2ab$ and $\wth=h/d$ is
primitive in $N(Y)$, the element $\wth^\ast=d^2\wth/(2ab)$, and
$u^\ast(\wth)+K(h)$ has the order $2ab/(d^2\gamma)$ in $K(h)^\ast/K(h)$.

They are related as follows:
$$
K(h)=h^\perp_{N(Y)}=[K(H),{2ab\over d^2\gamma}\,cu^\ast(\wH)],
\tag{2.3.23}
$$
and
$$
u^\ast(\wth)+K(h)=m(a,b)cu^\ast(\wH)+K(h)
\tag{2.3.24}
$$
where $m(a,b)\mod 2ab$ is defined by
$m(a,b)\equiv -1\mod 2a$ and $m(a,b)\equiv 1$ $\mod 2b$. To define
$u^\ast(\wth)+K(h)$ above, it is enough to consider
$m(a,b)\mod $ $2ab/(d^2\gamma)$.

We can formally put $h=H/c$, equivalently $\wth=\wH/c$. Then $N(X)$ and
$N(Y)$ become the  extensions of a common sublattice:
$$
N(X)\supset [\wH=c\,\wth,\,K(H)]\subset N(Y).
\tag{2.3.25}
$$
\endproclaim

Since $n(v)=1$, the transcendental lattices $T(X)=T(Y)$
are canonically identified in $v^\perp$. It follows that we have
the canonical identifications
$$
N(X)^\ast/N(X)=T(X)^\ast/T(X)=T(Y)^\ast/T(Y)=N(Y)^\ast/N(Y).
\tag{2.3.26}
$$
Here we use that discriminant groups of orthogonal complements in a
unimodular lattice are canonically isomorphic.
For example, here the identification \newline
$N(X)^\ast/N(X)=T(X)^\ast/T(X)$ is given by $n^\ast+N(X)\to t^\ast +T(X)$,
if $n^\ast+t^\ast \in H^2(X,\bz)$.

Let us calculate the identification
$$
N(X)^\ast/N(X)=N(Y)^\ast/N(Y).
\tag{2.3.27}
$$
Obviously (from the description above), it is given by the canonical maps
$$
N(X)^\ast \leftarrow (U\oplus N(X))^\ast\supset
{(U\oplus N(X))^\ast}^\perp_v\rightarrow (v^\perp)^\ast_0\leftarrow
(v^\perp/\bz v)^\ast.
\tag{2.3.28}
$$
Here ${(U\oplus N(X))^\ast}^\perp_v=\{x\in (U\oplus N(X))^\ast|x\cdot v=0\}$,
$(v^\perp)=\{x\in U\oplus N(X)|x\cdot v=0\}$ and
$(v^\perp)^\ast_0=\{x\in (v^\perp)^\ast |x\cdot v=0\}$.

Let $f^\ast=n\wH^\ast+k^\ast \in N(X)^\ast$ where $n\in\bz$ and
$k^\ast \in K(H)^\ast$. The element $\widetilde{f}^\ast=x_1e_1+x_2e_2+
n\wH^\ast+k^\ast$ is its lift to $(U\oplus N(X))^\ast$
where $x_1,\,x_2\in \bz$. We have $\widetilde{f}^\ast\in
{(U\oplus N(X))^\ast}^\perp_v$ if $-sx_1-rx_2+dn=0$. It follows that
$n=c(bx_1+ax_2)/d$ where $bx_1+ax_2\equiv 0\mod d$, and
$\widetilde{f}^\ast=x_1e_1+x_2e_2+
\left(c(bx_1+ax_2)/d\right)\wH^\ast+k^\ast$. It follows that
$$
f^\ast={c(bx_1+ax_2)\over d}\wH^\ast+k^\ast,\ \
bx_1+ax_2\equiv 0\mod d .
\tag{2.3.29}
$$
Like in \thetag{2.3.11}, \thetag{2.3.12}
and \thetag{2.3.14}, we finally get that the corresponding to
$f^\ast+N(X)$ element in $N(Y)^\ast /N(Y)$ is
$\epsilon(f^\ast)+N(Y)$ where
$$
\epsilon(f^\ast)={-bx_1+ax_2\over d}\wth^\ast+k^\ast\in N(Y)^\ast.
\tag{2.3.30}
$$
Thus, the identification \thetag{2.3.27} is given by
$f^\ast=c\left((bx_1+ax_2)/d\right)\wH^\ast+k^\ast+N(X)\rightarrow
\epsilon(f^\ast)=(-bx_1+ax_2)\wth^\ast+k^\ast+N(Y)$ where
$bx_1+ax_1\equiv 0\mod d$. We have $x_1=d_a\widetilde{x}_1$ and
$x_2=d_b\widetilde{x}_2$ where $\widetilde{x}_1,\,\widetilde{x}_2\in \bz$.
Then
$$
f^\ast=c\left((b/d_b)\widetilde{x}_1+
(a/d_a)\widetilde{x}_2\right)\wH^\ast+
k^\ast \rightarrow \epsilon(f^\ast)=
\left(-(b/d_b)\widetilde{x}_1+
(a/d_a)\widetilde{x}_2\right)\wth^\ast+k^\ast.
$$

Let us denote
$n_1=(b/d_b)\widetilde{x}_1+(a/d_a)\widetilde{x}_2$ and
$n_2=-(b/d_b)\widetilde{x}_1+(a/d_a)\widetilde{x}_2$. We have
$n_1\equiv m(a,b)n_2\mod {2ab/d^2}$. Really,
$n_1-m(a,b)n_2\equiv (2a/d_a)\widetilde{x}_2\equiv
0\mod 2a/d_a^2$ because
$m(a,b)\equiv -1\mod 2a/d_a^2$. We have
$n_1-m(a,b)n_2\equiv (2b/d_b)\widetilde{x}_1\equiv 0\mod
2b/d_b^2$ because
$m(a,b)\equiv 1\mod 2b/d_b^2$.
It follows the statement. Here we consider $m(a,b)\mod 2ab/d^2$.

The expression $(b/d_b)\widetilde{x}_1+(a/d_a)\widetilde{x}_2$
gives all integers $n\in \bz$. Thus, we have proved that
if $f^\ast=cn\wH^\ast+k^\ast\in N(X)^\ast$ where
$n\in \bz$, $k^\ast\in K(H)^\ast$, then
$\epsilon (f^\ast)=m(a,b)\wth^\ast+k^\ast\in N(Y)$.
Thus, we have proved

\proclaim{Proposition 2.3.2} In notations of Proposition 2.3.1,
the canonical identification
$\epsilon:N(X)^\ast/N(X)\cong N(Y)^\ast/N(Y)$
of the discriminant groups given by periods:
$$
\epsilon: N(X)^\ast/N(X)=T(X)^\ast/T(X)=T(Y)^\ast/T(Y)=N(Y)^\ast/N(Y)
\tag{2.3.31}
$$
is given by
$$
\epsilon:cn\wH^\ast+k^\ast+N(X)\mapsto m(a,b)n\wth^\ast+k^\ast
\tag{2.3.32}
$$
where $n\in\bz$, $k^\ast \in K(H)^\ast$. Here one should consider
$m(a,b)\mod~2ab/d^2$.

Equivalently, the characteristic maps
$\kappa(\wH):\widetilde{K}(H)\to N(X)^\ast/N(X)$ and
$\kappa (\wth):\widetilde{K}(h)\to N(Y)^\ast/N(Y)$ are related
as follows:
$$
\epsilon(\kappa(\wH)((cn,k^\ast)+\widetilde{K}(H)_0)=
\kappa(\wth)(m(a,b)n,k^\ast)+\widetilde{K}(h)_0
\tag{2.3.33}
$$
(here one should consider $m(a,b)\mod 2ab/d^2$).
Here we set $\widetilde{K}(H)=\widetilde{K}(\wH)$ and
$\widetilde{K}(h)=\widetilde{K}(\wth)$ (see \thetag{2.2.7}).
\endproclaim

This finishes the calculation of the periods of $N(Y)$ in terms of
the periods of $N(X)$.

Applying Lemma 2.2.1 and Propositions 2.3.1 and 2.3.2, by Global
Torelli Theorem for K3 surfaces \cite{PS}, we get

\proclaim{Theorem 2.3.3} Assume that $X$ is a K3 surface
with a polarization $H$ with $H^2=2rs$, $r,\,s\ge 1$, and
a primitive Mukai vector $v=(r,H,s)$ with the invariants
$$
(a,b,c,d,\gamma)
\tag{2.3.34}
$$
introduced above, i.e. $\wH=H/d$ is primitive.
Let $Y$ be the moduli space of coherent sheaves on
$X$ with the Mukai vector $v=(r,H,s)$.
We denote by $(K(\wH), u^\ast(\wH),$ $\kappa(\wH))$
the invariants \thetag{2.2.11} of $\wH\in N(X)$.

The transcendental periods $(T(X),H^{2,0}(X))$ and
$(T(Y), H^{2,0}(Y))$ are isomorphic if and only if
$$
n(v)=\text{g.c.d}(c,\,d\gamma)=1
\tag{2.3.35}
$$
(this is Mukai's result). We denote by $(K(H)=K(\wH), u^\ast(\wH),
\kappa(\wH))$ the invariants \thetag{2.2.11} of $\wH\in N(X)$.

Assume that \thetag{2.3.35} is valid. Then $Y\cong X$,
if the following conditions (a), (b) and (c) are valid:

(a) there exists a primitive $\widetilde{h}\in N(X)$ with
$\widetilde{h}^2=2ab/d^2$ and $\gamma(\widetilde{h})=\gamma$.

(b) There exists an embedding
$\phi:K(H)\subset K(\widetilde{h})$ of lattices such that
$\phi^\ast (K(\wth))=[K(H), 2abc/(d^2\gamma)u^\ast(\wH)]$ and
$\phi^\ast(u^\ast(\wth))+\phi^\ast(K(\wth))=m(a,b)cu^\ast (\wH)+
\phi^\ast(K(\wth))$. Here
$m(a,b)\mod~2ab/(d^2\gamma)$ is considered.

(c) There exists a choice of $\pm$ such that
$\kappa (\wth)(m(a,b)n,z^\ast)=\pm \kappa(H)(cn,\phi^\ast(z^\ast))$ if
$(cn,\phi^\ast(z^\ast))\in \widetilde{K}(H)^\ast=\widetilde{K}(\wH)^\ast$.
Here $m(a,b)\mod 2ab/d^2$ is considered.

The conditions (a), (b) and (c) are necessary for a K3 surface $X$ with
$\rho (X)\le 19$ which is general for its Picard lattice $N(X)$
in the following sense: the automorphism group of the transcendental periods
$(T(X),H^{2,0}(X))$ is $\pm 1$. If $\rho (X)=20$, then always
$Y\cong X$ if \thetag{2.3.35} holds.
\endproclaim

\example{Example 2.3.4} Let us assume that $\rho(X)=1$. Thus,
$N(X)=\bz \wH$. Assume conditions (a), (b) and (c) of Theorem 2.3.3
satisfy. Then $\wth=\pm \wH$. It follows $c=1$. The lattices $K(H)$ and
$K(\wth)$ are zero, and then (b) is valid. The discriminant group
$N(X)=\bz \wH^\ast/\bz \wH\equiv \bz/(2ab/d^2)\bz$. The condition (c)
is valid if and only if $m(a,b)\equiv \pm 1\mod 2ab/d^2$. This is
true if and only if $(a/d_a^2)=1$ or $(b/d_b^2)=1$.

Thus, $X\cong Y$ if $c=1$ and either $a_1=a/d_a^2=1$ or $b_1=b/d_b^2=1$.
These conditions are necessary to have $Y\cong X$,
if $X$ is a general K3 surface with $\rho(X)=1$. We recover the result
of Mukai from \cite{Mu3}.
\endexample

\example{Example 2.3.5} Following Example 2.3.5, let us consider
the case when we can satisfy conditions of Theorem 2.3.3 taking
$\wth=\pm H$ and $\phi=\pm id$. Again we get $c=1$.
(b) satisfies, if and only if
$m(a,b)\equiv \pm 1\mod 2ab/(d^2\gamma)$. This is
equivalent to either $2a/(d_a^2\gamma_a\gamma_2)\le 2$ or
$2b/(d_b^2\gamma_b\gamma_2)\le 2$. (c) satisfies, if and
only if $m(a,b)\equiv \pm 1\mod 2ab/d^2$.
Thus either $a/d_a^2=1$ or $b/d_b^2=1$. Thus, always $Y\cong X$,
if $c=1$ and either $a_1=a/d_a^2=1$ or $b_1=b/d_b^2=1$. This is just a
specialization (see Lemma 2.1.1)
of the $\rho=1$ case above. This result is also due to
Mukai \cite{Mu3}
\endexample

In Sect. 3 we consider $\rho=2$. We shall analyse when we can
satisfy conditions of Theorem 2.3.3 in this case. By specialization
(see Lemma 2.1.1) of these cases, we shall get results about K3 surfaces
with any Picard number $\rho\ge 2$.

\head
3. Conditions of $Y\cong X$ for a general K3
surface $X$ with $\rho=2$
\endhead

\subhead
3.1. Main results for $\rho(X)=2$
\endsubhead
Here we apply results of Sect. 2 to $X$ and $Y$ with
Picard number 2. Thus, we assume that
$\rho (X)=\rk N(X)=2$.

We start with some preliminary considerations on a primitive element
$P\in S$ of an even hyperbolic lattice $S$ of $\rk S=2$. We assume that
$P^2=2n$, $n\in \bn$, and $\gamma (P)=\gamma |2n$.

Let
$$
K(P)=P^\perp_{S}=\bz f(P)
\tag{3.1.1}
$$
and $f(P)^2=-t$ where $t>0$ is even. Then $\pm f(P)\in S$ is defined
uniquely by $P$. Below we set $f=f(P)$.

By elementary considerations, we have
$$
S=[\bz P,\, \bz f,\, {\gamma (\mu P+f)\over 2n}]
\tag{3.1.2}
$$
where
$$
\text{g.c.d}(\mu,{2n\over \gamma})=1.
\tag{3.1.3}
$$
The element
$$
\pm \mu \mod {{2n\over \gamma}} \in   (\bz/{2n\over \gamma})^\ast
\tag{3.1.4}
$$
is {\it the invariant of the pair $P\in S$}
up to isomorphisms of lattices with a primitive vector $P$ of
$P^2=2n$ and $\gamma(P)=\gamma$.
If $f$ changes to $-f$, then $\mu\mod 2n/\gamma$ changes to
$-\mu\mod 2n/\gamma$.

We have $(\gamma (\mu P+f)/(2n))^2=
\gamma^2(\mu^2-t/2n)/2n\equiv 0\mod 2$. It follows
$2n\mu^2-t\equiv 0 \mod 8n^2/\gamma^2$. It follows that for
some $\delta \in \bn$ we have
$$
f^2=-{2n\delta \over \gamma},\ \
n\delta\equiv 0\mod \gamma ,\ \ \text{and}\ \
\delta\equiv \mu^2\gamma \mod {{4n\over \gamma}}.
\tag{3.1.5}
$$
We have
$$
\det S=-\gamma\delta.
\tag{3.1.6}
$$

Any element $z\in S$ can be written as
$z=\gamma (xP+yf)/2n$ where $x \equiv \mu y \mod 2n/\gamma$. We have
$$
z^2={\gamma x^2-\delta y^2\over (2n/\gamma)}
\tag{3.1.7}
$$

It is convenient to put
$$
\widetilde{n}={2n\over \gamma}\ .
\tag{3.1.8}
$$
Thus, the considered above case of a primitive $P\in S$ where $S$
is an even hyperbolic lattice of $\rk S=2$ is described by
the invariants
$$
\widetilde{n},\ \gamma,\ \delta,\ \pm \mu\in (\bz/\widetilde{n})^\ast,
\tag{3.1.9}
$$
where $\widetilde{n},\,\gamma,\,\delta \in \bn$. The invariants \thetag{3.1.9}
must satisfy
$$
\widetilde{n}\gamma\equiv \widetilde{n}\delta\equiv 0\mod 2,\
\delta\equiv \mu^2\gamma\mod 2\widetilde{n}.
\tag{3.1.10}
$$
Then $P^2=\widetilde{n}\gamma$, $f^2=-\widetilde{n}\delta$, $P\perp f$, and
$$
S=\{{xP+yf\over \widetilde{n}}\ |\ x,\,y\in \bz,\ x\equiv \mu y\mod \widetilde{n}\}.
\tag{3.1.11}
$$
We have
$$
z^2={\gamma x^2-\delta y^2\over \widetilde{n}}\,.
\tag{3.1.12}
$$
Moreover,
$$
\det S=-\delta\gamma.
\tag{3.1.13}
$$
We denote
$$
P^\ast={P\over \widetilde{n}\gamma},\ f^\ast={f\over \widetilde{n}\delta}.
\tag{3.1.14}
$$
Then
$$
S^\ast=\{vP^\ast+wf^\ast\ |\ \mu v-w\equiv 0\mod \widetilde{n}\}
\tag{3.1.15}
$$
and
$$
S=\{vP^\ast+wf^\ast\ |\ v\equiv 0\mod \gamma,\ w\equiv 0\mod \delta,\
{v\over \gamma}\equiv {\mu w\over \delta }\mod \widetilde{n}\ \}.
\tag{3.1.16}
$$
Here $v,\,w\in \bz$. From \thetag{3.1.15}, $w=\mu v+\widetilde{n}t$,
and
$$
S^\ast=v(P^\ast+\mu f^\ast)+t\widetilde{n}f^\ast,\ \ v,\,t\in \bz,
\tag{3.1.17}
$$
and
$v(P^\ast+\mu f^\ast)+t\widetilde{n}f^\ast\in S$ if and only if
$$
v\equiv 0\mod \gamma,\ \mu v+\widetilde{n}t\equiv 0\mod \delta,\ \
\delta v\equiv
\gamma \mu(\mu v+\widetilde{n} t)\mod~{\delta \gamma \widetilde{n}}.
\tag{3.1.18}
$$
We have
$$
u^\ast(P)={\mu^{-1}f\over \widetilde{n}}+\bz f=\mu^{-1}\delta f^\ast+\bz f.
\tag{3.1.19}
$$

We remind notations we have used in Sect. 2:
$$
a_1={a\over d_a^2},\ b_1={b\over d_b^2}
\tag{3.1.20}
$$
where $d_a=\text{g.c.d}(d,a)$ and $d_b=\text{g.c.d}(d,b)$. We put
$$
a_2={a_1\over \gamma_a},\ b_2={b_1\over \gamma_b},\ e_2={2\over \gamma_2}.
\tag{3.1.21}
$$
where $\gamma_a=\text{g.c.d}(a_1,\gamma)$,
$\gamma_b=\text{g.c.d}(b_1,\gamma)$, $\gamma_2=\gamma/(\gamma_a\gamma_b)$.
Then $\gamma_2|2$.

Applying calculations above to $S=N(X)$ and primitive $P=\wH$ with
$n=2a_1b_1c^2$, $\widetilde{n}=2a_1b_1c^2/\gamma=e_2a_2b_2c^2$, and
$\gamma (\wH)=\gamma$, we get

\proclaim{Proposition 3.1.1} Let $X$ be a K3 surface with Picard
number $\rho = 2$ equipped with a primitive
polarization (or vector) $\wH\in N(X)$ of degree
$\wH^2=2a_1b_1c^2$ and $\gamma (\wH)=\gamma|2a_1b_1$.

Let $K(\wH)=(\wH)^\perp=\bz f(\wH)$. We have
$f(\wH)^2=-2a_1b_1c^2\delta/\gamma$ where $\det N(X)=-\delta\gamma$.

For some $\mu\in (\bz/(2a_1b_1c^2/\gamma))^\ast$ (the $\pm \mu$ is
the invariant of the pair $\wH\in N(X)$) where
$\delta\equiv\mu^2\gamma\mod 4a_1b_1c^2/\gamma$, one has
$$
N(X)=[\wH,f(\wH),{(\mu \wH+f(\wH))\over 2a_1b_1c^2/\gamma}],
\tag{3.1.22}
$$
$$
N(X)=\{z={x\wH+yf(\wH)\over 2a_1b_1c^2/\gamma } \ |\ x,y\in \bz\ \text{and}\
x \equiv \mu y \mod {{2a_1b_1c^2\over \gamma}} \}.
\tag{3.1.23}
$$
We have
$$
z^2={\gamma x^2-\delta y^2\over 2a_1b_1c^2/\gamma}.
\tag{3.1.24}
$$
For any primitive element $P \in N(X)$ with
$P^2=2a_1b_1c^2$, $\gamma (P)=\gamma$ and the same invariant $\pm \mu$,
there exists an automorphism
$\phi\in O(N(X))$ such that $\phi(\wH)=P$.
\endproclaim

Applying calculations above to $S=N(Y)$ and primitive $P=\wth\in N(Y)$
with $n=2a_1b_1$, $\widetilde{n}=2a_1b_1/\gamma=e_2a_2b_2$ and
$\gamma(\wth)=\gamma$, we get

\proclaim{Proposition 3.1.2} Let $Y$ be a K3 surface with Picard
number $\rho = 2$ equipped with a primitive
polarization (or vector) $\wth\in N(Y)$ of degree
$\wth^2=2a_1b_1$ and $\gamma (\wth)=\gamma|2a_1b_1$.

Let $K(\wth)=(\wth)^\perp=\bz f(\wth)$. We have
$f(\wth)^2=-2a_1b_1\delta/\gamma$ where $\det N(Y)=-\delta\gamma$.

For some $\nu\in (\bz/(2a_1b_1/\gamma))^\ast$ (the $\pm \nu$ is
the invariant of the pair $\wth\in N(Y)$) where
$\delta\equiv \nu^2\gamma \mod 4a_1b_1/\gamma$, one has
$$
N(Y)=[\wth,f(\wth),{(\nu \wth+f(\wth))\over 2a_1b_1/\gamma}],
\tag{3.1.25}
$$
$$
N(Y)=\{z={x\wth+yf(\wth)\over 2a_1b_1/\gamma } \ |\ x,y\in \bz\ \text{and}\
x \equiv \nu y \mod {{2a_1b_1\over \gamma}} \}.
\tag{3.1.26}
$$
We have
$$
z^2={\gamma x^2-\delta y^2\over 2a_1b_1/\gamma}.
\tag{3.1.27}
$$
For any primitive element $P \in N(Y)$ with
$P^2=2a_1b_1$, $\gamma (P)=\gamma$ and the same invariant $\pm \nu$,
there exists an automorphism
$\phi\in O(N(Y))$ such that $\phi(\wth)=P$.
\endproclaim

The crucial statement is

\proclaim{Theorem 3.1.3}
Let $X$ be a K3 surface, $\rho (X)=2$ and $H$ a polarization of
$X$ of degree $H^2=2rs$, $r,\,s\ge 1$, and Mukai vector $(r,H,s)$ is
primitive. Let $Y$ be the moduli space
of sheaves on $X$ with the isotropic Mukai vector $v=(r,H,s)$
and the canonical $nef$ element $h=(-a,0,b)\mod \bz v$. We
assume that
$$
\text{g.c.d}(c,d\gamma)=1.
$$

With notations of Propositions {3.1.1},
all elements
$$
\widetilde{h}={x\wH+yf(\wH)\over 2a_1b_1c^2/\gamma}\in N(X)
$$
with square
$\widetilde{h}^2=2a_1b_1$ satisfying Theorem 2.2.3 are in one to one
correspondence with integral solutions $(x,y)$ of the equation
$$
\gamma x^2-\delta y^2=4a_1^2b_1^2 c^2/\gamma
\tag{3.1.28}
$$
which satisfy conditions (i) --- (v) below:

(i)
$$
x\equiv \mu y\mod 2a_1b_1c^2/\gamma\, ,
\tag{3.1.29}
$$
$$
\mu \gamma x\equiv \delta y\mod 2a_1b_1c^2;
\tag{3.1.30}
$$

(ii) $(x,y)$ belongs to one of $a$-series (the sign $+$) or
$b$-series (the sign $-$) of solutions defined below:
$$
\pm m(a,b)\mu x+(\delta y/\gamma)\equiv 0\mod 2a_1b_1/\gamma\,,
\tag{3.1.31}
$$
$$
x\pm m(a,b)\mu y\equiv 0\mod 2a_1b_1/\gamma\, ,
\tag{3.1.32}
$$
$$
\pm m(a,b)\mu x+(\delta y/\gamma )\equiv \mu (x\pm m(a,b)\mu y)\mod
(2a_1b_1c^2/\gamma)(2a_1b_1/\gamma)\,;
\tag{3.1.33}
$$

(iii) there exists a choice of $\beta=\pm 1$ such that
$$
\cases
&\left({m(a,b)\gamma x\pm \mu \gamma y\over 2a_1b_1}-
\beta c\right)\equiv 0\mod \gamma\\
&\left({\delta m(a,b)y\pm \mu \gamma x\over 2a_1b_1}-
\beta\mu c\right)\equiv 0\mod \delta\\
&\delta\left({m(a,b)\gamma x\pm \mu \gamma y\over 2a_1b_1}-
\beta c\right)\equiv
\mu \gamma \left({\delta m(a,b)y\pm \mu \gamma x\over 2a_1b_1}-
\beta\mu c\right)\mod 2a_1b_1c^2\delta
\endcases
\tag{3.1.34}
$$
and
$$
\cases
&c\delta y\equiv 0\mod \gamma\\
&cx-(\pm \beta) {2a_1b_1c^2\over \gamma}\equiv 0\mod \delta\\
&\delta y\equiv
\mu \left(\gamma x -(\pm \beta) 2a_1b_1c\right)\mod 2a_1b_1c\delta
\endcases
\tag{3.1.35}
$$
where $+$ is taken for $a$-series, and $-$ is taken for $b$-series.

(iv) the pair $(x,y)$ is $\mu$-primitive:
$$
\text{g.c.d}\left(x,\,y,\,{x-\mu y\over 2a_1b_1c^2/\gamma}\right)=1;
\tag{3.1.36}
$$

(v) $\gamma (\wth)=\gamma$, equivalently
$$
\text{g.c.d}
\left(\gamma x,\,\delta y,\,
{\mu \gamma x-\delta y\over 2a_1b_1c^2/\gamma }\right)=
\gamma.
\tag{3.1.37}
$$

In particular (by Theorem 2.3.3),
for a general $X$ with $\rho(X)=2$
we have $Y\cong X$ if and only if the equation
$\gamma x^2-\delta y^2=4a_1^2b_1^2 c^2/\gamma$
has an integral solution $(x,y)$ satisfying conditions (i)---(v)
above. Moreover, a $nef$ primitive
element $P=(x\wH+f(\wH))/(2a_1b_1c^2/\gamma)$ with
$P^2=2a_1b_1$ and $\gamma (P)=\gamma$
defines a pair $(X,P)$ which is isomorphic to the
$(Y,\wth)$ if and only if $(x,y)$ satisfies the conditions (ii) and (iii)
(it satisfies conditions (i), (iv) and (v) since it corresponds
to a primitive element of $N(X)$ with $\gamma (P)=\gamma$).
\endproclaim

\demo{Proof} We denote
$$
H^\ast={\wH\over 2a_1b_1c^2},\
f(\wH)^\ast={\gamma f(\wH)\over 2a_1b_1c^2\delta}
\tag{3.1.38}
$$
where $K(H)=\bz f(\wH)=H^\perp$ in $N(X)$. By \thetag{3.1.19}, we have
$$
u^\ast(\wH)={\mu^{-1}\gamma f(\wH)\over 2a_1b_1c^2}+\bz f(\wH)=
\mu^{-1}\delta f(\wH)^\ast+\bz f(\wH).
\tag{3.1.39}
$$
Let
$$
\wth={x\wH+yf(\wH)\over 2a_1b_1c^2/\gamma } \in N(X)
\tag{3.1.40}
$$
satisfies conditions of Theorem 2.2.1. Then $x\, , y\in \bz$ and
$x\equiv \mu y\mod 2a_1b_1c^2/\gamma$. We get \thetag{3.1.29}
in (i). Moreover, $\wth$ is primitive which is equivalent to
(iv). We also have $\gamma (\wth)=\gamma$. It is equivalent to
$$
\text{g.c.d}\left(\wH\cdot \wth,\, f(\wH)\cdot \wth,
(\mu \wH+f(\wH))/(2a_1b_1c^2/\gamma)\cdot \wth \right)=\gamma.
$$
It follows \thetag{3.1.30} in (i), and (v).
We have $\widetilde{h}^2=2a_1b_1$. This is equivalent to
$\gamma x^2-\delta y^2=4a_1^2b_1^2c^2/\gamma$.

Consider $K(\wth)=\wth^\perp$ in $N(X)$. Let us denote
$$
f(\wth)={(\delta y/\gamma)\wH+xf(\wH)\over 2a_1b_1c^2/\gamma}.
\tag{3.1.41}
$$
The element $f(\wth)\in N(X)$ because of (i).
We have $f(\wth)\perp \wth$ and $f(\wth)^2=-2a_1b_1\delta/\gamma$.
Since $\gamma (\wth)=\gamma$, it follows that
$K(\wth)=\bz f(\wth)$ and $f(\wth)$ is primitive.
We have
$$
N(X)=\left[\wth,f(\wth),{\nu \wth+f(\wth)\over 2a_1b_1/\gamma}\right]
\tag{3.1.42}
$$
where $\nu\in \left(\bz/(2a_1b_1/\gamma)\right)^\ast$ (according to
Proposition 3.1.2). We denote
$$
\wth^\ast={\wth\over 2a_1b_1},\ \ f(\wth)^\ast={\gamma f(\wth)\over
2a_1b_1\delta}.
\tag{3.1.43}
$$
By \thetag{3.1.19}, we have
$$
u^\ast(\wth)={\nu^{-1}f(\wth)\over 2a_1b_1/\gamma}+\bz f(\wth)=
\nu^{-1}\delta f(\wth)^\ast+\bz f(\wth)\,.
\tag{3.1.44}
$$

There exists a unique (up to $\pm 1$)
embedding
$$
\phi:K(\wH)=\bz f(\wH)\to K(\wth)=\bz f(\wth),\ \ \phi(f(\wH))=\pm cf(\wth).
\tag{3.1.45}
$$
of one-dimensional lattices. Its dual is defined by
$\phi^\ast(f(\wth)^\ast )=\pm c f(\wH)^\ast$.

We have
$$
\phi^\ast(K(\wth))=\bz \phi^\ast(f(\wth))=\bz f(\wH)/c=
$$
$$
[K(\wH), \,(2a_1b_1c/\gamma ) u^\ast(\wH)]=
[K(\wH),\,(2abc/(d^2\gamma))u^\ast(\wH)]
$$
because of \thetag{3.1.39}. This gives the first part of (b)
in Theorem 2.3.3.

We have
$$
\phi^\ast(u^\ast (\wth))+\phi^\ast (K(\wth))=
\nu^{-1}\delta \phi^\ast(f(\wth)^\ast)+\bz \phi^\ast(f(\wth))=
\pm \nu^{-1}\delta c f(\wH)^\ast + \bz f(\wH)/c=
$$
$$
\pm \nu^{-1}\delta c f(\wH)^\ast +
\bz (2a_1b_1c\delta/\gamma )f(\wH)^\ast.
$$
On the other hand, by \thetag{3.1.39}
$$
m(a,b)cu^\ast(\wH)+\bz(2a_1b_1c\delta/\gamma)f(\wH)^\ast=
m(a,b)\mu^{-1}\delta c f(\wH)^\ast+\bz (2a_1b_1c\delta/\gamma)f(\wH)^\ast.
$$
Thus, by \thetag{3.1.44}, second part of (b) in Theorem 2.3.3 is
$\pm \nu^{-1}\equiv m(a,b)\mu^{-1}\mod 2a_1b_1/\gamma$.
Equivalently,
$$
\nu\equiv \pm m(a,b)\mu \mod {{2a_1b_1\over \gamma}}.
\tag{3.1.46}
$$
Thus, for $\nu$ given by \thetag{3.1.46} one has (this is the definition of
$\nu$)
$$
{\nu \wth + f(\wth)\over 2a_1b_1/\gamma} =
\left({\nu x+(\delta y/\gamma)\over 2a_1b_1/\gamma }\wH+
{x+\nu y\over 2a_1b_1/\gamma}f(\wH)\right)/(2a_1b_1c^2/\gamma)\in N(X).
\tag{3.1.47}
$$
This is equivalent for $\nu\equiv m(a,b)\mu\mod 2a_1b_1/\gamma $ to
$$
\cases
\pm m(a,b)\mu x+(\delta y/\gamma)\equiv 0\mod 2a_1b_1/\gamma \\
x\pm m(a,b)\mu y\equiv 0\mod 2a_1b_1/\gamma\\
\pm m(a,b)\mu x+(\delta y/\gamma )\equiv \mu (x \pm m(a,b)\mu y)\mod
(2a_1b_1c^2/\gamma)(2a_1b_1/\gamma)
\endcases .
\tag{3.1.48}
$$
This gives (ii).

Let us consider the condition (c) of Theorem 2.3.3. For a choice of
$\beta=\pm 1$, one has
$$
m(a,b)n\wth^\ast+z^\ast\equiv \beta (cn\wH^\ast+\phi^\ast(z^\ast))\mod N(X),
\ \ n\in \bz,\ z\in K(\wth)^\ast
\tag{3.1.49}
$$
if
$$
cn\wH^\ast+\phi^\ast(z^\ast) \in N(X)^\ast.
\tag{3.1.50}
$$
Let $z^\ast=kf(\wth)^\ast$, $k\in \bz$,
then $\phi^\ast (kf(\wth)^\ast)=\pm k c f(\wH)^\ast$
and $cn\wH^\ast+\phi^\ast(z^\ast)=cn\wH^\ast \pm kc f(\wH)^\ast$.
By \thetag{3.1.15}, $cn\wH^\ast \pm kc f(\wH)^\ast\in N(X)^\ast$ if
and only if $\mu c n \mp kc\equiv 0\mod 2a_1b_1c^2/\gamma$. This is
equivalent to
$k\equiv \pm \mu n\mod 2a_1b_1c/\gamma$. Like in \thetag{3.1.17}, we get
$k=\pm  \mu n+(2a_1b_1c/\gamma)t$ where $n,\,t\in \bz$. We get
$$
cn\wH^\ast \pm kc f(\wH)^\ast=cn(\wH^\ast+\mu f(\wH)^\ast)\pm
t(2a_1b_1c^2/\gamma )f(\wH)^\ast.
\tag{3.1.51}
$$
Thus, it is enough to check
$$
m(a,b)n\wth^\ast+kf(\wth)^\ast
\equiv \beta (cn\wH^\ast\pm kc f(\wH)^\ast)\mod N(X),
\tag{3.1.52}
$$
where $(n,k)=(1,\pm \mu)$ or $(n,k)=(0,2a_1b_1c/\gamma)$. Thus, one
should check for one of $\beta=\pm 1$ that
$$
m(a,b)\wth^\ast\pm \mu f(\wth)^\ast
\equiv \beta (c\wH^\ast+c\mu f(\wH)^\ast)\mod N(X)
\tag{3.1.53}
$$
and
$$
{2a_1b_1cf(\wth)^\ast\over \gamma}
\equiv
\pm \beta {2a_1b_1c^2f(\wH)^\ast\over \gamma}\mod N(X).
\tag{3.1.54}
$$

We have
$$
m(a,b)n\wth^\ast+kf(\wth)^\ast=
m(a,b)n{\wth\over 2a_1b_1}+k{\gamma f(\wth)\over 2a_1b_1\delta}=
$$
$$
{1\over 2a_1b_1\delta}\left(\delta m(a,b)n\wth+\gamma kf(\wth)\right)=
$$
$$
{1\over 2a_1b_1\delta}
\left(\delta m(a,b)n{x\wH+yf(\wH)\over 2a_1b_1c^2/\gamma}+
\gamma k {(\delta y /\gamma)\wH+x f(\wH)\over 2a_1b_1c^2/\gamma}\right)=
$$
$$
{1\over 2a_1b_1\delta}
\left(\delta m(a,b)n(\gamma x\wH^\ast+\delta y f(\wH)^\ast)
+ \gamma k (\delta y \wH^\ast +\delta x f(\wH)^\ast \right)=
$$
$$
{m(a,b)n\gamma x+\gamma k y\over 2a_1b_1}\wH^\ast+
{\delta m(a,b)ny+\gamma k x\over 2a_1b_1}f(\wH)^\ast.
\tag{3.1.55}
$$
Thus, \thetag{3.1.53} is
$$
\left({m(a,b)\gamma x\pm \mu \gamma y\over 2a_1b_1}-
\beta c\right)\wH^\ast+
\left({\delta m(a,b)y\pm \mu \gamma x\over 2a_1b_1}-
\beta\mu c\right)f(\wH)^\ast\in N(X)
\tag{3.1.56}
$$
and \thetag{3.1.54} is
$$
c y\wH^\ast+
\left(c x -(\pm \beta) {2a_1b_1c^2\over \gamma}\right)f(\wH)^\ast \in N(X).
\tag{3.1.57}
$$
Here one should take $+$ if $(x,y)$ belongs to $a$-series, and
one should take $-$ if $(x,y)$ belongs to $b$-series.

By \thetag{3.1.16}, we can reformulate \thetag{3.1.56} as
\thetag{3.1.34} in (iii), and we can reformulate
\thetag{3.1.57} as \thetag{3.1.35} in (iii).

This finishes the proof.
\enddemo

Now we analyse conditions of Theorem 3.1.3. The most important are
congruences $\mod \delta$ since $\delta$ is not bounded by a constant
depending on $(r,s)$.

Let us consider the congruence
$cx-(\pm \beta)(2a_1b_1c^2/ \gamma )\equiv 0\mod \delta$
in \thetag{3.1.35}. We know that
$\delta\equiv \gamma \mu^2\mod (2a_1b_1/\gamma)c^2$ where
$\mu \mod(2a_1b_1/\gamma)c^2$ is invertible,
$\gamma \vert 2a_1b_1$ and $\text{g.c.d}(\gamma,c)=1$. It
follows that $\text{g.c.d}(c,\delta)=1$. Thus, the considered congruence
is equivalent to
$$
x\equiv \pm {2a_1b_1c\over \gamma}\mod \delta
\tag{3.1.58}
$$
since $\beta=\pm 1$. Later we shall see that all congruences
$\mod \delta$ which follow from conditions of Theorem 3.1.3 are
consequences of \thetag{3.1.58}. Thus, \thetag{3.1.58} is the
most important condition of Theorem 3.1.3.

Let us consider all integral $(x,y)$ which satisfy \thetag{3.1.28} and
\thetag{3.1.58}, i. e.
$$
\gamma x^2-\delta y^2={4a_1^2b_1^2 c^2\over \gamma}\ \ \text{and}\
x\equiv \pm {2a_1b_1c\over \gamma}\mod \delta.
\tag{3.1.59}
$$
We apply the main trick used in \cite{MN1}, \cite{MN2} and \cite{N4}.
Considering $\pm (x,y)$, we can assume that
$x\equiv 2a_1b_1c/\gamma \mod \delta$,
i. e.  $x=2a_1b_1c/\gamma - k\delta$ where $k\in \bz$.
We have
$\gamma^2x^2=4a_1^2b_1^2c^2-4a_1b_1ck\delta\gamma + \gamma^2 k^2\delta^2=
\gamma \delta y^2+4a_1^2b_1^2c^2$. It follows
$$
\delta={y^2+4a_1b_1c k\over \gamma k^2}.
\tag{3.1.60}
$$
Consider a prime $l$ such that $l\not\vert 2a_1b_1c$. Assume that
$l^{2t+1}\vert k$, but $l^{2t+2}\not\vert k$. We have $k\vert y^2$.
Then $l^{2t+1}\vert y^2$. Then $l^{2t+2}\vert y^2$. Since
$k^2\vert y^2+4a_1b_1ck$, it follows that
$l^{2t+2}| y^2+4a_1b_1ck$
We then get $l^{2t+2}|k$. We get a contradiction. It follows that
$k=-\alpha q^2$ where $q\in \bz$ (we can additionally assume that $q\ge 0$),
$\alpha\vert 2a_1b_1c$ and $\alpha$
is square-free. Remark that $\alpha$ can be negative.

From \thetag{3.1.60}, we get
$$
\delta={y^2-4a_1b_1c\alpha q^2\over \gamma \alpha^2q^4}.
\tag{3.1.61}
$$
It follows $\alpha q \vert y$ and $y=\alpha pq$ where
$p\in \bz$. From \thetag{3.1.61},
$$
\delta={p^2-4a_1b_1c/\alpha\over \gamma q^2}.
\tag{3.1.62}
$$
Equivalently,
$$
p^2-\gamma\delta q^2={4a_1b_1c\over \alpha}.
\tag{3.1.63}
$$
If integral $p$, $q$ satisfy \thetag{3.1.63}, we have
$$
(x,\,y)=
\pm \left({2a_1b_1c\over \gamma}+\alpha \delta q^2,\,\alpha pq\right).
\tag{3.1.64}
$$
satisfy \thetag{3.1.59}.
We call solutions \thetag{3.1.64} of \thetag{3.1.59} as associated
solutions. Thus, we get the very important for us

\proclaim{Theorem 3.1.4}
All integral solutions $(x,y)$ of
$$
\gamma x^2-\delta y^2={4a_1^2b_1^2 c^2\over \gamma}\ \ \text{and}\
x\equiv \pm {2a_1b_1c\over \gamma} \mod \delta
\tag{3.1.65}
$$
are associated solutions
$$
(x,\,y)=
\pm \left({2a_1b_1c\over \gamma}+\alpha \delta q^2,\,\alpha pq\right)
\tag{3.1.66}
$$
to integral solutions $\alpha$, $(p,q)$ of
$$
\alpha\vert 2a_1b_1c\ \text{where $\alpha$ is
square-free, and}\ \
p^2- \gamma \delta q^2={4a_1b_1c\over \alpha}\,.
\tag{3.1.67}
$$
Any solution $(x,y)$ of \thetag{3.1.65} can be written in the form
\thetag{3.1.66} where $\alpha$, $(p,q)$ is solution of \thetag{3.1.67}. Any
solution of \thetag{3.1.67} gives a solution \thetag{3.1.66} of
\thetag{3.1.65}.

Solutions of \thetag{3.1.65} and \thetag{3.1.67} are in one-to-one
correspondence if we additionally assume that $q\ge 0$.
\endproclaim

Now we can write $(x,y)$ of Theorem 3.1.3 in the form
\thetag{3.1.66} as associated solutions to \thetag{3.1.67}.
Putting that $(x,y)$ to relations (i)---(v) of Theorem 3.1.3, we
get some relations on $\alpha$ and $(p,q)$. They are a finite
number of congruences $\mod N_i$ where $N_i$ depend only on
$(r,s)$ (or $(a,b,c)$). All $N_i$ are bounded by functions
depending only on $(r,s)$. These congruences have a lot of
relations between them and with \thetag{3.1.67}. All together they
give many very strong restrictions on $\alpha$ and $(p,q)$.
We analyse them below.

We fix
$$
\mu\in \left(\bz/(2a_1b_1c^2/\gamma)\right)^\ast
\tag{3.1.68}
$$
and consider $\delta \in \bn$ such that
$$
\delta\equiv\mu^2\gamma\mod {{4a_1b_1c^2\over \gamma}}
\tag{3.1.69}
$$

The relation \thetag{3.1.29} in (i) of Theorem 3.1.3 is equivalent to
$$
2a_1b_1c+\alpha\gamma \delta q^2\equiv \gamma \mu \alpha p q\mod
2a_1b_1c^2.
\tag{3.1.70}
$$
By \thetag{3.1.67}, we have
$$
-4a_1b_1c+\alpha p^2-\alpha \gamma \delta q^2=0.
\tag{3.1.71}
$$
Taking sum, we get
$$
2a_1b_1c\equiv \alpha p(p-\mu \gamma q)\mod 2a_1b_1c^2.
\tag{3.1.72}
$$
The relation \thetag{3.1.72} is equivalent to \thetag{3.1.29}.
Taking 2\thetag{3.1.70}+\thetag{3.1.71}, we get
$$
\alpha p^2+\alpha \gamma \delta q^2\equiv
2\gamma\mu\alpha p q\mod 4a_1b_1c^2
\tag{3.1.73}
$$
which is also equivalent to \thetag{3.1.29}.
From \thetag{3.1.69}, we get
$$
\alpha p^2+\alpha \mu^2\gamma^2 q^2\equiv 2\gamma \mu \alpha pq
\mod {{4a_1b_1c^2}} \tag{3.1.74}
$$
and
$$
\alpha (p-\mu \gamma q)^2\equiv 0\mod {{4a_1b_1c^2}}. \tag{3.1.75}
$$
This is equivalent to \thetag{3.1.73}. It follows $\alpha (p-\mu
\gamma q)^2\equiv 0\mod 4c^2$. Since $\alpha$ is square-free, it
follows
$$
2c\vert (p-\mu \gamma q). \tag{3.1.76}
$$
From \thetag{3.1.72}, we get
$$
2a_1b_1\equiv \alpha p{p-\mu \gamma q\over c}\mod {2a_1b_1c}\,
\tag{3.1.77}
$$
where $(p-\mu \gamma q)/c$ is an integer.
It follows, $\alpha \vert 2a_1b_1$.
Thus, we get

\proclaim{Proposition 3.1.5} The condition \thetag{3.1.29} of
Theorem 3.1.3 is equivalent to
$$
\alpha (p-\mu \gamma q)^2\equiv 0 \mod 4a_1b_1c^2. \tag{3.1.78}
$$
We have
$$
\alpha \vert 2a_1b_1. \tag{3.1.79}
$$
The condition \thetag{3.1.78} is
also equivalent to
$$
2a_1b_1c\equiv \alpha p(p-\mu \gamma q)\mod 2a_1b_1c^2.
\tag{3.1.80}
$$
\endproclaim

Considering the condition \thetag{3.1.30} of Theorem 3.1.3, we
similarly get

\proclaim{Proposition 3.1.6} The condition \thetag{3.1.30} of
Theorem 3.1.3 is equivalent to
$$
\mu \alpha p^2+\mu \alpha \gamma \delta q^2\equiv 2\alpha \delta pq
\mod 4a_1b_1c^2.
\tag{3.1.81}
$$
Congruences \thetag{3.1.78} and \thetag{3.1.81} are equivalent
$\mod (4a_1b_1c^2/\gamma)$, and \thetag{3.1.81} is equivalent to
\thetag{3.1.78} together with
$$
2\alpha pq(\delta-\gamma \mu^2)\equiv 0\mod 4a_1b_1c^2
\tag{3.1.82}
$$
where $\delta-\gamma\mu^2 \equiv 0\mod 4a_1b_1c^2/\gamma$.
\endproclaim

\demo{Proof} To get \thetag{3.1.82}, consider
$\mu$ \thetag{3.1.78}\ -\ \thetag{3.1.81}.
\enddemo

Now consider \thetag{3.1.31} of Theorem 3.1.3. Let us consider
the $a$-series (the sign $+$).
Since $m(a,b)\equiv -1\mod 2a_1/(\gamma_2\gamma_a)$, we get
$-\mu x+(\delta y/\gamma)\equiv 0\mod 2a_1/(\gamma_2\gamma_a)$.
It is satisfied because of \thetag{3.1.30}. Since
$m(a,b)\equiv 1\mod 2b_1/(\gamma_2\gamma_b)$, we get
$\mu x+(\delta y/\gamma)\equiv 0\mod 2b_1/(\gamma_2\gamma_b)$.
By \thetag{3.1.30}, we get
$\mu x-(\delta y/\gamma)\equiv 0\mod 2b_1/(\gamma_2\gamma_b)$.
Thus, we get
$2\mu x\equiv 0\mod 2b_1/(\gamma_2\gamma_b)$. If $\gamma_2=1$, this
is equivalent to $\mu x\equiv 0\mod b_1/\gamma_b$ and
$x\equiv 0\mod b_1/\gamma_b$. If $\gamma_2=2$, then $b_1/\gamma_b$ is
odd, and we get $2\mu x\equiv 0\mod (b_1/\gamma_b)$ which is equivalent to
$x\equiv 0\mod b_1/\gamma_b$. Thus, at any case we get
$(b_1/\gamma_b)|x$, equivalently, $x\equiv 0\mod (b_1/\gamma_b)$.
By \thetag{3.1.66}, this is equivalent to
$\alpha \delta q^2\equiv 0\mod (b_1/\gamma_b)$. We have
$\delta\equiv \mu^2\gamma \mod (b_1/\gamma_b)$ and
$\mu\mod (b_1/\gamma_b)$ is invertible. Thus, we get
$\alpha \gamma q^2\equiv 0\mod (b_1/\gamma_b)$ which is equivalent to
$\alpha(\gamma_b q)^2\equiv 0\mod b_1$.

Let us consider \thetag{3.1.32} of Theorem 3.1.3. We consider
the $a$-series. We get
$x-\mu y \equiv 0\mod 2a_1/(\gamma_2\gamma_a)$. It satisfies
because of \thetag{3.1.29}. We get
$x+\mu y\equiv 0\mod 2b_1/(\gamma_2\gamma_b)$.
Using \thetag{3.1.29}, we similarly get that this is equivalent to
$(b_1/\gamma_b)|x$.

Thus, we finally get

\proclaim{Proposition 3.1.7} The conditions \thetag{3.1.31} and
\thetag{3.1.32} in (ii) of Theorem 3.1.3 are equivalent to
$$
x\equiv 0\mod {{b_1\over \gamma_b}}
\tag{3.1.83}
$$
or to
$$
\alpha(\gamma_b q)^2\equiv 0\mod b_1
\tag{3.1.84}
$$
for $a$-series (the sign $+$),
and it is equivalent to
$$
x\equiv 0\mod {{a_1\over \gamma_a}}
\tag{3.1.85}
$$
or to
$$
\alpha (\gamma_a q)^2\equiv 0\mod a_1
\tag{3.1.86}
$$
for $b$-series (the sign $-$).
\endproclaim

Consider \thetag{3.1.33} in Theorem 3.1.3. We consider the
$a$-series. We get \newline
$\mod 2a_1/(\gamma_2\gamma_a)$ that
$$
2\mu \gamma x\equiv (2\delta+(\mu^2\gamma -\delta))y\mod
(2a_1b_1c^2)\left(2a_1/(\gamma_2\gamma_a)\right).
$$
We get $\mod 2b_1/(\gamma_2\gamma_b)$ that
$$
(\delta-\mu^2\gamma)y\equiv 0\mod
(2a_1b_1c^2)(2b_1/\gamma_2\gamma_b).
$$
Using \thetag{3.1.66} (and \thetag{3.1.71} to make the
relations homogeneous), we finally get

\proclaim{Proposition 3.1.8} The condition
\thetag{3.1.33} of Theorem 3.1.3 is equivalent to
$$
2\mu \gamma x\equiv (\delta+\mu^2\gamma)y\mod
(2a_1b_1c^2)\left(2a_1/(\gamma_2\gamma_a)\right)
\tag{3.1.87}
$$
and
$$
(\delta-\mu^2\gamma)y\equiv 0\mod
(2a_1b_1c^2)\left(2b_1/(\gamma_2\gamma_b)\right)
\tag{3.1.88}
$$
or
$$
\mu \alpha p^2-\alpha (\delta+\mu^2\gamma)pq+\mu\alpha\gamma\delta q^2
\equiv 0
\mod (2a_1b_1c^2)\left(2a_1/(\gamma_2\gamma_a)\right)
\tag{3.1.89}
$$
and
$$
\alpha(\delta-\mu^2\gamma)pq\equiv 0\mod
(2a_1b_1c^2)\left(2b_1/(\gamma_2\gamma_b)\right)
\tag{3.1.90}
$$
for the $a$-series (the sign $+$), and it is equivalent to
$$
2\mu \gamma x\equiv (\delta+\mu^2\gamma)y\mod
(2a_1b_1c^2)\left(2b_1/(\gamma_2\gamma_b)\right)
\tag{3.1.91}
$$
and
$$
(\delta-\mu^2\gamma)y\equiv 0\mod
(2a_1b_1c^2)\left(2a_1/(\gamma_2\gamma_a)\right)
\tag{3.1.92}
$$
or
$$
\mu \alpha p^2-\alpha (\delta+\mu^2\gamma)pq+\mu\alpha\gamma\delta q^2
\equiv 0 \mod
(2a_1b_1c^2)\left(2b_1/(\gamma_2\gamma_b)\right)
\tag{3.1.93}
$$
and
$$
\alpha(\delta-\mu^2\gamma)pq\equiv 0\mod
(2a_1b_1c^2)\left(2a_1/(\gamma_2\gamma_a)\right)
\tag{3.1.94}
$$
for the $b$-series (the sign $-$).
\endproclaim

Consider the condition \thetag{3.1.35} of Theorem 3.1.3.
Consider the $a$-series (the sign $+$). Second condition in
\thetag{3.1.35} gives $x\equiv \beta (2a_1b_1c/\gamma)\mod \delta$.
By \thetag{3.1.66}, we get
$(x,y)=\beta(2a_1b_1c/\gamma+\alpha\delta q^2,\,\alpha pq)$. Then
third condition in \thetag{3.1.35} is equivalent to
$$
\beta \delta \alpha pq\equiv
\mu(\beta 2a_1b_1c+\beta\alpha \gamma\delta q^2-\beta 2a_1b_1c)\mod
2a_1b_1c\delta.
$$
This is equivalent to
$\alpha q (p-\mu \gamma q)\equiv 0\mod 2a_1b_1c$. Since
$\gamma\vert 2a_1b_1c$, it follows the first condition in \thetag{3.1.35}
which is $c\delta \beta \alpha pq\equiv 0\mod \gamma$.

For $b$-series we get the same. Thus, we have

\proclaim{Proposition 3.1.9} The condition \thetag{3.1.35} of
Theorem 3.1.3 is equivalent to
$$
\alpha q (p-\mu \gamma q)\equiv 0\mod 2a_1b_1c.
\tag{3.1.95}
$$
\endproclaim

Consider \thetag{3.1.34} of Theorem 3.1.3. We consider the
$a$-series. Then
$$
(x,y)=\beta({2a_1b_1c\over \gamma}+\alpha\delta q^2,\,\alpha pq)
\tag{3.1.96}
$$
The first relation of \thetag{3.1.34} gives
$$
m(a,b)\gamma x + \mu \gamma y\equiv 2a_1b_1\beta c
\mod 2a_1b_1\gamma.
$$
Using \thetag{3.1.96}, we get
$$
m(a,b)\alpha\gamma\delta q^2+\mu\alpha \gamma pq\equiv 2a_1b_1c(1-m(a,b))
\mod 2a_1b_1\gamma.
$$
This is equivalent to
$$
-\alpha \delta q^2+\mu \alpha pq\equiv {4a_1b_1c\over \gamma}\mod 2a_1
$$
and
$$
\alpha\delta q^2+\mu\alpha pq\equiv 0 \mod 2b_1.
$$
Using \thetag{3.1.71}, we can rewrite the first relation in the
homogeneous form
$$
\alpha p(p-\mu\gamma q)\equiv 0\mod 2a_1\gamma.
$$

The second relation of \thetag{3.1.34} gives
$$
\delta m(a,b)y+\mu\gamma x - 2a_1b_1\beta\mu c\equiv 0\mod 2a_1b_1\delta.
$$
By \thetag{3.1.96} it is equivalent to
$$
m(a,b)\alpha pq+\alpha\mu\gamma q^2\equiv 0\mod 2a_1b_1.
$$
This is equivalent to
$$
-\alpha pq+\alpha\mu\gamma q^2\equiv 0\mod 2a_1
$$
and
$$
\alpha pq+\alpha\mu\gamma  q^2\equiv 0\mod 2b_1.
$$

The third relation in \thetag{3.1.34} is
$$
\split
&\delta\left(m(a,b)\gamma x+\mu \gamma y-2\beta a_1b_1c\right)\equiv\\
&\mu \gamma \left(\delta m(a,b)y+\mu\gamma x-2\mu \beta a_1b_1c \right)
\mod (2a_1b_1c^2\delta)(2a_1b_1).
\endsplit
$$
Using \thetag{3.1.96}, one can calculate that it is equivalent to
$$
\split
&\mu \alpha\gamma pq\left(m(a,b)-1\right)+
\alpha\gamma q^2\left(\gamma \mu^2-m(a,b)\delta\right)\\
&\equiv
2a_1b_1c\left(m(a,b)-1\right)\mod (2a_1b_1c^2)(2a_1b_1).
\endsplit
$$
This is equivalent to
$$
4a_1b_1c\equiv \alpha\gamma q\left(2\mu p-(\gamma\mu^2+\delta) q\right)
\mod(2a_1b_1c^2)(2a_1)
$$
and
$$
\alpha\gamma q^2(\delta-\gamma\mu^2)\equiv 0\mod (2a_1b_1c^2)(2b_1).
$$
Using \thetag{3.1.71}, we can rewrite the first relation in the
homogeneous form
$$
\alpha(p-\mu \gamma q)^2\equiv 0\mod (2a_1b_1c^2)(2a_1).
$$

Thus, we get

\proclaim{Proposition 3.1.10} The condition \thetag{3.1.34} of
Theorem 3.1.3 is equivalent to the system of congruences:

\noindent
For the $a$-series (the sign $+$)
$$
\alpha p(p-\mu\gamma q)\equiv 0\mod 2a_1\gamma,
\tag{3.1.97}
$$
$$
\alpha\delta q^2+\alpha\mu pq\equiv 0 \mod 2b_1,
\tag{3.1.98}
$$
$$
-\alpha pq+\alpha\mu\gamma q^2\equiv 0\mod 2a_1,
\tag{3.1.99}
$$
$$
\alpha pq+\alpha\mu\gamma q^2\equiv 0\mod 2b_1,
\tag{3.1.100}
$$
$$
\alpha(p-\mu \gamma q)^2\equiv 0\mod (2a_1b_1c^2)(2a_1),
\tag{3.1.101}
$$
$$
\alpha\gamma q^2(\delta-\gamma\mu^2)\equiv 0\mod (2a_1b_1c^2)(2b_1).
\tag{3.1.102}
$$

\noindent
For the $b$-series (the sign $-$):
$$
\alpha p(p-\mu\gamma q)\equiv 0\mod 2b_1\gamma,
\tag{3.1.103}
$$
$$
\alpha\delta q^2+\alpha \mu pq\equiv 0 \mod 2a_1,
\tag{3.1.104}
$$
$$
-\alpha pq+\alpha\mu\gamma q^2\equiv 0\mod 2b_1,
\tag{3.1.105}
$$
$$
\alpha pq+\alpha\mu\gamma q^2\equiv 0\mod 2a_1,
\tag{3.1.106}
$$
$$
\alpha(p-\mu \gamma q)^2\equiv 0\mod (2a_1b_1c^2)(2b_1),
\tag{3.1.107}
$$
$$
\alpha\gamma q^2(\delta-\gamma\mu^2)\equiv 0\mod (2a_1b_1c^2)(2a_1).
\tag{3.1.108}
$$
\endproclaim

Consider the condition (iv) of Theorem 3.1.3.
It means that the corresponding element $\wth$ is primitive. Since
$\wth^2=2a_1b_1$ and the lattice $N(X)$ is even, it is not valid
only if $\wth/l\in N(X)$ for some prime $l$ such that $l^2\vert a_1b_1$.
Thus, \thetag{3.1.36} is not valid if and only if
$$
x\equiv y\equiv {x-\mu y\over 2a_1b_1c^2/\gamma}\equiv 0\mod l
$$
for some prime $l$ such that $l^2\vert a_1b_1$.
Using \thetag{3.1.66} and \thetag{3.1.71},
we get

\proclaim{Proposition 3.1.11} The condition (iv) of Theorem 3.1.3 is
equivalent to the non-existence of a prime $l$ such
that $l^2\vert a_1b_1$ and
$$
x\equiv y\equiv {x-\mu y\over 2a_1b_1c^2/\gamma}\equiv 0\mod l.
\tag{3.1.109}
$$
Equivalently, the system of congruences
$$
\cases
\alpha p^2+\alpha \gamma\delta q^2\equiv 0\mod 2\gamma l\\
\alpha pq\equiv 0\mod l\\
\alpha p^2+\alpha \gamma\delta q^2\equiv 2\alpha \gamma\mu pq
\mod 4a_1b_1c^2 l
\endcases
\tag{3.1.110}
$$
is not satisfied for any prime $l$ such that $l^2\vert a_1b_1$.
\endproclaim

Consider the condition (v) of Theorem 3.1.3. This is equivalent to
$\gamma (\wth)=\gamma$ where $\wth\cdot N(X)=\gamma(\wth) \bz$.
All other conditions of Theorem 3.1.3 give that
$\gamma \vert \wth \cdot N(X)$, and $\gamma\vert \gamma(\wth)$. Since
$\wth^2=2a_1b_1$, it follows that $(\gamma(\wth)/\gamma)\,
\vert\, 2a_1b_1/\gamma$.
Equivalently, (v) is not satisfied if and only if
for some prime $l\vert 2a_1b_1/\gamma$ one has
$$
x\equiv {\delta y\over \gamma} \equiv
{\mu \gamma x-\delta y\over 2a_1b_1c^2}\equiv 0\mod \l.
$$
Using \thetag{3.1.66} and \thetag{3.1.71},
we get

\proclaim{Proposition 3.1.12} The condition (v) of Theorem 3.1.3 is
equivalent to the non-existence of a prime $l\,\vert\, (2a_1b_1/\gamma)$
such that
$$
x\equiv {\delta y\over \gamma} \equiv
{\mu \gamma x-\delta y\over 2a_1b_1c^2}\equiv 0\mod \l.
\tag{3.1.111}
$$
Equivalently, the system of congruences
$$
\cases
\alpha p^2+\alpha \gamma\delta q^2\equiv 0\mod 2\gamma l\\
\delta \alpha pq\equiv 0\mod \gamma l\\
\mu \alpha p^2+\mu \alpha \gamma\delta q^2\equiv 2\alpha\delta pq
\mod 4a_1b_1c^2 l
\endcases
\tag{3.1.112}
$$
is not satisfied for any prime $l\,\vert\, (2a_1b_1/\gamma)$.
\endproclaim

Now we collect analysed conditions of Theorem 3.1.3 all together.
We divide them in {\it general conditions: $(G^\prime)$}
which are valid for both
$a$ and $b$-series, {\it conditions $(A^\prime)$} which are valid
for the $a$-series,
and {\it conditions $(B^\prime)$}
which are valid for the $b$-series.

\smallpagebreak

\noindent
{\bf $(G^\prime)$: General conditions}
$$
\alpha (p-\mu \gamma q)^2\equiv 0 \mod 4a_1b_1c^2,
\tag{3.1.113}
$$
$$
2 \alpha pq (\delta-\gamma \mu^2) \equiv 0 \mod 4a_1b_1c^2,
\tag{3.1.114}
$$
$$
\alpha q (p-\mu \gamma q)\equiv 0\mod 2a_1b_1c,
\tag{3.1.115}
$$
$$
\cases
\alpha p^2+\alpha \gamma\delta q^2\equiv 0\mod 2\gamma l\\
\alpha pq\equiv 0\mod l\\
\alpha p^2+\alpha \gamma\delta q^2\equiv 2\alpha \gamma\mu pq
\mod 4a_1b_1c^2 l
\endcases
\tag{3.1.116}
$$
is not satisfied for any prime $l$ such that $l^2\vert a_1b_1$,
$$
\cases
\alpha p^2+\alpha \gamma\delta q^2\equiv 0\mod 2\gamma l\\
\delta \alpha pq\equiv 0\mod \gamma l\\
\mu \alpha p^2+\mu \alpha \gamma\delta q^2\equiv 2\alpha\delta pq
\mod 4a_1b_1c^2 l
\endcases
\tag{3.1.117}
$$
is not satisfied for any prime $l\,\vert\, (2a_1b_1/\gamma)$.

\smallpagebreak

\noindent
{\bf $(A^\prime)$: Conditions of the $a$-series}

$$
\alpha (\gamma_b\,q)^2\equiv 0\mod b_1, \tag{3.1.118}
$$
$$
\mu \alpha p^2-\alpha (\delta+\mu^2\gamma)pq+\mu\alpha\gamma\delta q^2
\equiv 0
\mod (2a_1b_1c^2)\left(2a_1/(\gamma_2\gamma_a)\right)
\tag{3.1.119}
$$
$$
\alpha(\delta-\mu^2\gamma)pq\equiv 0\mod
(2a_1b_1c^2)\left(2b_1/(\gamma_2\gamma_b)\right),
\tag{3.1.120}
$$
$$
\alpha p(p-\mu\gamma q)\equiv 0\mod 2a_1\gamma,
\tag{3.1.121}
$$
$$
\alpha\delta q^2+\alpha\mu pq\equiv 0 \mod 2b_1,
\tag{3.1.122}
$$
$$
-\alpha pq+\alpha\mu\gamma q^2\equiv 0\mod 2a_1,
\tag{3.1.123}
$$
$$
\alpha pq+\alpha\mu\gamma q^2\equiv 0\mod 2b_1,
\tag{3.1.124}
$$
$$
\alpha(p-\mu \gamma q)^2\equiv 0\mod (2a_1b_1c^2)(2a_1),
\tag{3.1.125}
$$
$$
\alpha\gamma q^2(\delta-\gamma\mu^2)\equiv 0\mod (2a_1b_1c^2)(2b_1).
\tag{3.1.126}
$$

\smallpagebreak

\noindent
{\bf $(B^\prime )$: Conditions of the $b$-series}

$$
\alpha (\gamma_a\,q)^2\equiv 0\mod a_1, \tag{3.1.127}
$$
$$
\mu \alpha p^2-\alpha (\delta+\mu^2\gamma)pq+\mu\alpha\gamma\delta q^2
\equiv 0 \mod
(2a_1b_1c^2)\left(2b_1/(\gamma_2\gamma_b)\right)
\tag{3.1.128}
$$
$$
\alpha(\delta-\mu^2\gamma)pq\equiv 0\mod
(2a_1b_1c^2)\left(2a_1/(\gamma_2\gamma_a)\right),
\tag{3.1.129}
$$
$$
\alpha p(p-\mu\gamma q)\equiv 0\mod 2b_1\gamma,
\tag{3.1.130}
$$
$$
\alpha\delta q^2+\alpha \mu pq\equiv 0 \mod 2a_1,
\tag{3.1.131}
$$
$$
-\alpha pq+\alpha\mu\gamma q^2\equiv 0\mod 2b_1,
\tag{3.1.132}
$$
$$
\alpha pq+\alpha\mu\gamma q^2\equiv 0\mod 2a_1,
\tag{3.1.133}
$$
$$
\alpha(p-\mu \gamma q)^2\equiv 0\mod (2a_1b_1c^2)(2b_1),
\tag{3.1.134}
$$
$$
\alpha\gamma q^2(\delta-\gamma\mu^2)\equiv 0\mod (2a_1b_1c^2)(2a_1).
\tag{3.1.135}
$$

\subhead
3.2. Simplification of the conditions $(G^\prime)$,
$(A^\prime)$ and $(B^\prime)$
\endsubhead
We have the following fundamental result which completely
determines $\alpha$ up to multiplication by $\pm 1$.

\proclaim{Lemma 3.2.1} For the $a$-series, we have that the
square-free $\alpha \vert b_1$ and
$$
b_1/|\alpha|=\widetilde{b}_1^2,\ \widetilde{b}_1>0,
$$
is a square.

Respectively, for the $b$-series, we have that the square-free
$\alpha \vert a_1$ and
$$
a_1/|\alpha|=\widetilde{a}_1^2,\ \widetilde{a}_1>0,
$$
is a square.
\endproclaim

\demo{Proof} First let us prove that
$$
\alpha\vert a_1b_1. \tag{3.2.1}
$$
Otherwise, $2a_1b_1/\alpha$ is odd. Let us consider the
$a$-series. By \thetag{3.1.125}, we have $\alpha(p-\mu\gamma
q)^2\equiv 0\mod {4c^2}$. Since $\alpha$ is square-free, it
follows that $p-\mu \gamma q\equiv 0\mod {2c}$. It follows
$p+\mu\gamma q\equiv 0\mod 2$. Then $p^2-\mu^2\gamma^2\equiv 0\mod
{4c}$. We have $\mu^2\gamma\equiv \delta\mod 4a_1b_1c^2/\gamma$.
Then $\mu^2\gamma^2\equiv \gamma\delta\mod {4c^2}$. Thus, we
obtain $p^2-\gamma\delta q^2\equiv 0\mod 4c$. On the other hand,
$p^2-\gamma\delta q^2=4a_1b_1c/\alpha\equiv 2c\mod 4c$ if
$2a_1b_1/\alpha$ is odd. We get a contradiction. It proves
\thetag{3.2.1}.

Now let us consider the $a$-series, and let us prove that
$$
\alpha \vert b_1. \tag{3.2.2}
$$
Otherwise, for a prime $l$ one has $l\vert \alpha$, $l\vert a_1$,
but $l$ does not divide $b_1$.

By \thetag{3.1.125}, we have $\alpha(p-\mu\gamma q)^2\equiv 0\mod
{4a_1^2c^2}$. Since $\alpha$ is square-free, it follows that
$p-\mu \gamma q\equiv 0\mod {2a_1c}$. It follows $p+\mu\gamma
q\equiv 0\mod 2$. Then $p^2-\mu^2\gamma^2\equiv 0\mod {4a_1c}$. We
have $\mu^2\gamma\equiv \delta\mod 4a_1b_1c^2/\gamma$. Then
$\mu^2\gamma^2\equiv \gamma\delta\mod {4a_1c^2}$. Thus, we obtain
$p^2-\gamma\delta q^2\equiv 0\mod {4a_1c}$. Thus, we obtain
$p^2-\gamma\delta q^2=4a_1b_1c/\alpha\equiv 0\mod {4a_1c}$.
Equivalently, $a_1b_1/\alpha\equiv 0\mod a_1$. This gives a
contradiction if for a prime $l$ one has $l\vert \alpha$, $\l\vert
a_1$, but $l$ does not divide $b_1$. This proves \thetag{3.2.2}.

Now let us prove that
$$
b_1/|\alpha| \tag{3.2.3}
$$
is a square. Otherwise, for a prime $l$ we have $l^{2t-1}\vert
b_1/\alpha$, but $l^{2t}$ does not divide $b_1/\alpha$ where $t\ge
1$. By \thetag{3.1.118}, $(\gamma q)^2\equiv 0\mod {b_1/\alpha}$.
It follows that
$$
l^t\vert \gamma q. \tag{3.2.4}
$$
By \thetag{3.1.113}, we have
$$
p-\gamma \mu q\equiv 0\mod {l^t}. \tag{3.2.5}
$$
By \thetag{3.2.4} and \thetag{3.1.5}, we obtain $l^t|p$.

From $p^2-\gamma\delta q^2=4a_1(b_1/\alpha)c$ we then get a
contradiction if $l$ does not divide $2c$.

We have $$ \gamma \delta\equiv\gamma^2\mu^2\mod {4a_1b_1c^2}.
\tag{3.2.6}
$$

From $p^2-\gamma\delta q^2=4a_1(b_1/\alpha)c$ and \thetag{3.2.6}
we then get
$$
p^2-\gamma^2\mu^2q^2\equiv (p-\gamma\mu q)(p+\gamma\mu q)\equiv
4a_1(b_1/\alpha)c\mod {4a_1b_1c^2}.
$$
Then, using \thetag{3.2.5}, we get
$$
{p-\gamma \mu q\over 2a_1cl^t}(p+\gamma\mu q)\equiv
2(b_1/\alpha)/l^t\mod {2(b_1/l^t)c} \tag{3.2.7}
$$
where $(p-\gamma \mu q)/(2a_1cl^t)$ is an integer.

If $l\vert c$ and $l$ is odd, \thetag{3.2.5} and \thetag{3.2.7}
give a contradiction because $l^t\vert p+\gamma\mu q$, $l^t\vert
2(b_1/l^t)c$, but $l^t$ does not divide $2(b_1/\alpha)/l^t$.

Now assume that $l=2\vert c$. If $2^{t+1}\vert \gamma\mu q$, then
$2^{t+1}\vert p$ by \thetag{3.2.5}, and we get a contradiction in
the same way. Assume that $2^{t+1}$ does not divide $\gamma\mu q$.
Then $\gamma\mu q/2^t$ and $p/2^t$ are both odd, and $2^{t+1}\vert
p+\gamma\mu q$. It also leads to a contradiction in the same way.

Now assume that $l=2$ and $c$ is odd. By \thetag{3.1.126}, we
obtain
$$
\gamma\delta q^2\equiv (\gamma \mu q)^2\mod 2^{4t}. \tag{3.2.8}
$$
Assume that $2^{t+1}\vert \mu\gamma q$. By \thetag{3.2.8}, then
$\gamma\delta q^2\equiv 0\mod 2^{2t+2}$. By \thetag{3.2.5},
$2^{t+1}\vert p$ and $2^{2t+2}\vert p^2$. Then $2^{2t+1}\vert
p^2-\gamma\delta q^2=4a_1(b_1/\alpha)c$ which gives a
contradiction because $4a_1(b_1/\alpha)c$ is divisible by
$2^{2t+1}$ only.

Now assume that $2^{t+1}$ does not divide $\mu\gamma q$. By
\thetag{3.2.8}, we get $\gamma\delta q^2\equiv 2^{2t}\mod
2^{2t+2}$. By \thetag{3.2.5} then $2^t\vert p$, but $2^{t+1}$ does
not divide $p$. It follows that $p^2\equiv 2^{2t}\mod 2^{2t+2}$.
Then $2^{2t+2}\vert p^2-\gamma\delta q^2=4a_1(b_1/\alpha)c$ which
again leads to a contradiction.

This finishes the proof of the theorem.
\enddemo

\proclaim{Lemma 3.2.2} Assume that
$$
\text{g.c.d}(\delta-\gamma \mu^2,\
4a_1b_1c^2)=(4a_1b_1c^2/\gamma)\gamma_0
$$
where $\gamma_0\vert \gamma$.

For any $u\vert a_1b_1c^2/(\gamma_a\gamma_b)$ and
$\text{g.c.d}(u,\gamma/\gamma_0)=1$ we can choose
$$\mu\in
(\bz/((2a_1b_1c^2/\gamma)\gamma_0u))^\ast
$$
such that
$$
\delta\equiv \gamma\mu^2\mod {(4a_1b_1c^2/\gamma)\gamma_0u}.
$$
\endproclaim

\demo{Proof} Assume $\mu_0\mod 2a_1b_1c^2/\gamma \in
(\bz/(2a_1b_1c^2/\gamma))^\ast$ and
$$
\delta\equiv \gamma \mu_0^2\mod {4a_1b_1c^2/\gamma}.
$$
Taking $\mu=\mu_0+(2a_1b_1c^2/\gamma)k$, we get
$$
\delta-\gamma\mu^2=\delta-\gamma \mu_0^2-4\mu_0a_1b_1c^2k-
(4a_1^2b_1^2c^4/\gamma)k^2.
$$
Then
$$
(\delta-\gamma\mu^2)/(4a_1b_1c^2/\gamma)\equiv
(\delta-\gamma\mu_0^2)/(4a_1b_1c^2/\gamma)+\gamma\mu_0 k\mod
\gamma_0u.
$$
Since $(\gamma/\gamma_0)\mu_0$ are invertible $\mod u$, we can
choose $k$ such that
$$
(\delta-\gamma\mu^2)/(4a_1b_1c^2/\gamma)\equiv 0\mod \gamma_0u .
$$
It follows the statement.
\enddemo

Further we consider the $a$-series (similar results will be valid
for $b$-series).

The congruence \thetag{3.1.125} implies \thetag{3.1.113}.
The \thetag{3.1.125} is equivalent to
$$
p-\mu \gamma q\equiv 0\mod {2a_1\widetilde{b}_1c}\,. \tag{3.2.9}
$$
Thus, \thetag{3.1.113} and \thetag{3.1.125} are equivalent to
\thetag{3.2.9}.

The congruence \thetag{3.1.118} is equivalent to $\gamma_b q\equiv
0\mod {\widetilde{b}_1}$. Equivalently,
$$
q={\widetilde{b}_1q_1\over \gamma_b},\ \ \widetilde{b}_1q_1\equiv
0\mod {\gamma_b} \tag{3.2.10}
$$
where $q_1$ is an integer.

By \thetag{3.2.9}, we have
$p-\mu\gamma_2\gamma_a\widetilde{b}_1q_1\equiv 0\mod
2a_1\widetilde{b}_1c$, and then
$\gamma_2\gamma_a\widetilde{b}_1\vert p$. Then
$$
p=\gamma_2\gamma_a\widetilde{b}_1p_1 \tag{3.2.11}
$$
where $p_1$ is an integer. The congruence \thetag{3.2.9} is then
equivalent to
$$
p_1-\mu q_1\equiv 0\mod {(2/\gamma_2)(a_1/\gamma_a)c}.
\tag{3.2.12}
$$

Denoting
$$
\alpha =\pm b_1/\widetilde{b}_1^2, \tag{3.2.13}
$$
we can rewrite $p^2-\gamma\delta q^2=4a_1b_1c/\alpha$ as
$$
\gamma p_1^2-\delta q_1^2=\pm 2
(2/\gamma_2)(a_1/\gamma_a)\gamma_bc. \tag{3.2.14}
$$

Now we can rewrite conditions $(G^\prime)$, $(A^\prime)$ and
$(B^\prime)$ using the introduced $(p_1,\,q_1)$.

As we have seen, the conditions \thetag{3.1.113}, \thetag{3.1.118}
and \thetag{3.1.125} are equivalent to \thetag{3.2.12} and
$$
\widetilde{b}_1q_1\equiv 0\mod \gamma_b. \tag{3.2.15}
$$

The condition \thetag{3.1.114} gives
$$
p_1q_1(\delta-\gamma \mu^2)\equiv 0\mod
(2/\gamma_2)(a_1/\gamma_a)c^2\gamma_b.
$$
Since $\text{g.c.d}(\gamma,c)=\text{g.c.d}(\gamma_b,a_1c^2)=1$ and
$\delta-\gamma\mu^2\equiv 0\mod 4a_1b_1c^2/\gamma$, it is
equivalent to
$$
p_1q_1(\delta-\gamma \mu^2)\equiv 0\mod (2/\gamma_2)\gamma_b.
\tag{3.2.16}
$$

The condition \thetag{3.1.115} gives
$$
q_1(p_1-\mu q_1)\equiv 0\mod (2/\gamma_2)(a_1/\gamma_a)c\gamma_b.
$$
Since $\text{g.c.d}(\gamma,c)=\text{g.c.d}(\gamma_b,a_1c^2)=1$ and
$p_1-\mu q_1\equiv 0\mod (2/\gamma_2)(a_1/\gamma_a)c$ (by
\thetag{3.2.12}), it is similarly equivalent to
$$
q_1(p_1-\mu q_1)\equiv 0\mod (2/\gamma_2)\gamma_b\,. \tag{3.2.17}
$$

The condition \thetag{3.1.120} gives
$$
(\delta-\gamma \mu^2)p_1q_1\equiv 0\mod
{(4/\gamma_2^2)b_1c^2}\,. \tag{3.2.18}
$$
Since $(\delta-\mu^2\gamma)\equiv 0\mod
{4b_1c^2/(\gamma_2\gamma_b)}$, the congruence \thetag{3.2.18} is
actually a congruence $\mod {\gamma_b}$ on $p_1q_1$. It also
implies \thetag{3.2.16}.

The condition \thetag{3.1.121} gives $b_1p_1(p_1-\mu q_1)\equiv 0
\mod {(2/\gamma_2)(a_1/\gamma_a)\gamma_b}$. It satisfies because
of \thetag{3.2.12} and $\gamma_b\vert b_1$.

The condition \thetag{3.1.122} gives
$$
q_1(\delta q_1+\mu\gamma p_1)\equiv 0\mod 2 \gamma_b^2.
\tag{3.2.19}
$$

The condition \thetag{3.1.123} gives $(b_1/\gamma_b)q_1(-p_1+\mu
q_1)\equiv 0\mod {(2/\gamma_2)(a_1/\gamma_a)}$. It satisfies by
\thetag{3.2.12}.

The condition \thetag{3.1.124} gives
$$
\gamma_a q_1(p_1+\mu q_1)\equiv 0\mod {(2/\gamma_2)\gamma_b}.
\tag{3.2.20}
$$

It is easy to see that \thetag{3.2.17}, \thetag{3.2.20} together
with $\delta\equiv \mu^2 \gamma \mod {4a_1b_1c^2/\gamma}$ imply
\thetag{3.2.18}.

Taking $\pm \gamma_a$ \thetag{3.2.17} plus \thetag{3.2.20}, we
obtain
$$
\gamma_2 p_1q_1\equiv 0\mod \gamma_b \tag{3.2.21}
$$
and $\gamma_2\mu q_1^2\equiv 0\mod \gamma_b$. Since $\mu$ can be
always taken coprime to $\gamma_b$, the last congruence is
equivalent to
$$
\gamma_2 q_1^2\equiv 0\mod \gamma_b. \tag{3.2.22}
$$
Any of them together with \thetag{3.2.17} can be taken to replace
\thetag{3.2.20}.

The condition \thetag{3.1.126} is equivalent to
$$
(\delta-\gamma \mu^2) q_1^2\equiv 0 \mod
{(4/\gamma_2)b_1c^2\gamma_b}. \tag{3.2.23}
$$
By \thetag{3.2.22} and  $(\delta-\mu^2\gamma)\equiv 0\mod
{4b_1c^2/(\gamma_2\gamma_b)}$, the congruence \thetag{3.2.23} is
actually a congruence $\mod {\gamma_b}$.

It is easy to see that \thetag{3.2.23} and \thetag{3.2.20} imply
\thetag{3.2.19}.

The condition \thetag{3.1.119} gives
$$
\mu\gamma
p_1^2-(\delta+\mu^2\gamma)p_1q_1+\mu\delta q_1^2 \equiv 0\mod
{(4/\gamma_2^2)(a_1^2/\gamma_a^2)\gamma_b c^2}.
$$
Since $a_1$ and $b_1$ are coprime, this is equivalent to
two congruences
$$
\mu\gamma
p_1^2-(\delta+\mu^2\gamma)p_1q_1+\mu\delta q_1^2 \equiv 0\mod
{(4/\gamma_2^2)c^2\gamma_b}
$$
and
$$
\mu\gamma p_1^2-(\delta+\mu^2\gamma)p_1q_1+\mu\delta q_1^2 \equiv
0\mod {(4/\gamma_2^2)c^2 (a_1^2/\gamma_a^2)}. \tag{3.2.24}
$$
By \thetag{3.2.23}, we have $\delta q_1^2 \equiv \gamma \mu^2
q_1^2 \mod {(4/\gamma_2^2)c^2\gamma_b}$, and the first congruence
gives
$$
\mu \gamma (p_1-\mu q_1)^2-(\delta -\mu^2\gamma)p_1q_1
\equiv 0\mod {(4/\gamma_2^2)c^2\gamma_b}
$$
It satisfies because of \thetag{3.2.12} and \thetag{3.2.18}. The
congruence \thetag{3.2.24} can be written as
$$
\mu \gamma (p_1-\mu q_1)^2+(\delta - \mu^2\gamma)(p_1-\mu q_1)q_1
\equiv 0\mod {(4/\gamma_2^2)c^2 (a_1^2/\gamma_a^2)}.
$$
It satisfies because of \thetag{3.2.12} and since
$\delta -\gamma \mu^2\equiv 0\mod {(4/\gamma_2)(a_1/\gamma_a)}c^2$.

The condition \thetag{3.1.116} is equivalent to
$$
\cases
(b_1/\gamma_b)(\gamma p_1^2+\delta q_1^2)\equiv 0\mod 2\gamma_b l\\
(b_1/\gamma_b)\gamma_2\gamma_a p_1q_1\equiv 0\mod l\\
\gamma^2 p_1^2+\gamma\delta q_1^2\equiv 2\gamma^2\mu p_1q_1 \mod
4a_1c^2\gamma_b^2 l
\endcases
\tag{3.2.25}
$$
is not satisfied for any prime $l$ such that $l^2\vert a_1b_1$.

By its meaning, the congruences \thetag{3.2.25} satisfies if we
formally put $l=1$. Assume that for a prime $l$ we have $l^2\vert
a_1b_1$ and $\text{g.c.d}(l,\gamma)=1$. Then \thetag{3.2.25} is
equivalent to
$$
\cases
b_1(\gamma p_1^2+\delta q_1^2)\equiv 0\mod 2l\\
b_1p_1q_1\equiv 0\mod l\\
\gamma p_1^2+\delta q_1^2\equiv 2\gamma \mu p_1q_1 \mod
(4/\gamma_2)(a_1/\gamma_a)c^2 l
\endcases
$$
is not satisfied. By Lemma 3.2.2, we can assume that \newline
$\delta\equiv \gamma \mu^2\mod (4/\gamma_2)(a_1/\gamma_a)l c^2$,
and the last condition is equivalent to
$$
\cases
b_1(\gamma p_1^2+\delta q_1^2)\equiv 0\mod 2l\\
b_1p_1q_1\equiv 0\mod l\\
\gamma (p_1-\mu q_1)^2\equiv 0 \mod (4/\gamma_2)(a_1/\gamma_a)l
c^2
\endcases
$$
is not satisfied. By \thetag{3.2.12}, this is equivalent to
$$
\cases
b_1(\gamma p_1^2+\delta q_1^2)\equiv 0\mod 2l\\
b_1p_1q_1\equiv 0\mod l\\
\gamma_a [(p_1-\mu q_1)/(2c/\gamma_2)]^2\equiv 0 \mod
(a_1/\gamma_a)l
\endcases
\tag{3.2.26}
$$
is not satisfied. Assume that $l\vert b_1$. By \thetag{3.2.14},
then the first and second congruences satisfy and \thetag{3.2.26}
is equivalent to
$$
p_1-\mu q_1 \equiv 0\mod {(2/\gamma_2)(a_1/\gamma_a)cl}
\tag{3.2.27}
$$
does not satisfy. Assume that $l\vert a_1$. By \thetag{3.2.12}
then third congruence in \thetag{3.2.26} satisfies, and by
\thetag{3.2.12} and \thetag{3.2.14} the condition \thetag{3.2.26}
is equivalent to
$$
l\not\vert p_1. \tag{3.2.28}
$$

The condition \thetag{3.1.117} is equivalent to
$$
\cases
(b_1/\gamma_b)(\gamma p_1^2+\delta q_1^2)\equiv 0\mod 2\gamma_b l\\
(b_1/\gamma_b) p_1q_1\equiv 0\mod \gamma_b l\\
\mu \gamma^2 p_1^2+\mu \gamma\delta q_1^2\equiv 2\delta\gamma
p_1q_1 \mod 4a_1c^2\gamma_b^2 l
\endcases
\tag{3.2.29}
$$
is not satisfied for any prime $l\,\vert\, (2a_1b_1/\gamma)$.

By its meaning, \thetag{3.2.29} satisfies if we formally put
$l=1$. Assume that \newline $\text{g.c.d}(l,\gamma)=1$ and $l\vert
2a_1b_1/\gamma$. Then \thetag{3.2.29} is equivalent to
$$
\cases
b_1(\gamma p_1^2+\delta q_1^2)\equiv 0\mod 2l\\
b_1p_1q_1\equiv 0\mod l\\
\mu \gamma p_1^2+\mu \delta q_1^2\equiv 2\delta p_1q_1 \mod
(4/\gamma_2)(a_1/\gamma_2)c^2 l
\endcases
$$
is not satisfied.

First assume that $l\vert a_1b_1$. By Lemma 3.2.2, we can assume
that
\newline $\delta\equiv \gamma \mu^2\mod
(4/\gamma_2)(a_1/\gamma_a)l c^2$, and the last condition is
equivalent to
$$
\cases
b_1(\gamma p_1^2+\delta q_1^2)\equiv 0\mod 2l\\
b_1p_1q_1\equiv 0\mod l\\
\mu\gamma (p_1-\mu q_1)^2\equiv 0 \mod (4/\gamma_2)(a_1/\gamma_a)l
c^2
\endcases
$$
is not satisfied. We obtain exactly the same conditions
\thetag{3.2.27} and \thetag{3.2.28} as above.

Now assume that $l\not\vert a_1b_1$. Then $l=2$ and $a_1$, $b_1$,
$\gamma$, $\delta$ are odd. We then obtain that
$$
\cases
\gamma p_1^2+\delta q_1^2\equiv 0\mod 2l\\
p_1q_1\equiv 0\mod l\\
\mu \gamma p_1^2+\mu \delta q_1^2\equiv 2\delta p_1q_1 \mod
(4/\gamma_2)(a_1/\gamma_2)c^2 l
\endcases
\tag{3.2.30}
$$
 is not satisfied. From first two congruences we get that
 $p_1\equiv q_1\equiv 0\mod 2$. Then we get
$q_1\delta\equiv q_1\gamma \mu^2\mod (4/\gamma_2)(a_1/\gamma_a)l
c^2$, and we can rewrite \thetag{3.2.30} again as
$$
\cases
\gamma p_1^2+\delta q_1^2\equiv 0\mod 2l\\
b_1p_1q_1\equiv 0\mod l\\
\mu\gamma (p_1-\mu q_1)^2\equiv 0 \mod (4/\gamma_2)(a_1/\gamma_a)l
c^2
\endcases
\tag{3.2.31}
$$
is not satisfied. From the last congruence in \thetag{3.2.31} we
get $\gamma (p_1-\mu q_1) \equiv 0\mod{4c}$ and $p_1+\mu q_1
\equiv 0 \mod{2}$. It follows that $\gamma p_1^2 -\gamma \mu^2
q_1^2 \equiv 0\mod{8c}$. Moreover, we have $q_1\delta\equiv
q_1\gamma \mu^2\mod{8c}$. Thus, we obtain $\gamma p_1^2-\delta
q_1^2\equiv 0\mod{8c}$. It leads to a contradiction since $\gamma
p_1^2-\delta q_1^2=\pm 2(2/\gamma_2)(a_1/\gamma_a)\gamma_b c\equiv
4c\mod{8c}$.

\smallpagebreak

In fact, \thetag{3.2.29} always follows from \thetag{3.2.29} if
$l\vert 2a_1b_1$ and $\text{g.c.d}(l,\gamma)=1$. Really, by our
construction, \thetag{3.2.29} means that the corresponding element
$\wth\in N(X)$ of Theorem 3.1.3 is not divisible by $l$. We have
$\mu\in (\bz/(2a_1b_1c^2/\gamma))^\ast$,  $\delta\equiv \gamma
\mu^2\mod {4a_1b_1c^2/\gamma}$, and $\det N(X)=-\gamma\delta$ is
not divisible by $l$. Then $N(X)\cdot \wth$ is not divisible by
$l$, and \thetag{3.2.29} follows. In particular, we can assume
that $l^2\vert a_1b_1$.

Thus, we can rewrite the conditions $(G^\prime)$
and $(A^\prime)$ respectively in the form

\noindent {\bf (A): The conditions of $a$-series:}

\noindent {\bf (AG): The general conditions of $a$-series:}
$$
p_1-\mu q_1\equiv 0\mod {(2/\gamma_2)(a_1/\gamma_a)c},
\tag{3.2.32}
$$
$$
\split &p_1-\mu q_1\not\equiv 0\mod
{(2/\gamma_2)(a_1/\gamma_a)c\,l}\\
&\text{for any prime $l$ such that}\ l^2\vert b_1\ \text{and}\
\text{g.c.d}(l,\gamma)=1, \endsplit \tag{3.2.33}
$$
$$
l \not\vert p_1\ \text{for any prime $l$ such that}\ l^2\vert a_1\
\text{and}\ \text{g.c.d}(l,\gamma)=1. \tag{3.2.34}
$$

\noindent {\bf (AS) The singular conditions of $a$-series:}
$$
\widetilde{b}_1q_1\equiv 0\mod \gamma_b, \tag{3.2.35}
$$
$$
q_1(p_1-\mu q_1)\equiv 0 \mod {(2/\gamma_2)\gamma_b}, \tag{3.2.36}
$$
$$
\gamma_2p_1q_1 \equiv \gamma_2q_1^2 \equiv 0\mod\gamma_b,
\tag{3.2.37}
$$
$$
(\delta-\gamma \mu^2) q_1^2\equiv 0 \mod
{(4/\gamma_2)b_1c^2\gamma_b}. \tag{3.2.38}
$$
$$
\cases
(b_1/\gamma_b)(\gamma p_1^2+\delta q_1^2)\equiv 0\mod 2\gamma_b l\\
(b_1/\gamma_b)\gamma_2\gamma_a p_1q_1\equiv 0\mod l\\
\gamma^2 p_1^2+\gamma\delta q_1^2\equiv 2\gamma^2\mu p_1q_1 \mod
4a_1c^2\gamma_b^2 l
\endcases
\tag{3.2.39}
$$
is not satisfied for any prime $l$ such that $l^2\vert a_1b_1$ and
$l\vert \gamma$,
$$
\cases
(b_1/\gamma_b)(\gamma p_1^2+\delta q_1^2)\equiv 0\mod 2\gamma_b l\\
(b_1/\gamma_b) p_1q_1\equiv 0\mod \gamma_b l\\
\mu \gamma^2 p_1^2+\mu \gamma\delta q_1^2\equiv 2\delta\gamma
p_1q_1 \mod 4a_1c^2\gamma_b^2 l
\endcases
\tag{3.2.40}
$$
is not satisfied for any prime $l\,\vert\, (2a_1b_1/\gamma)$ and
$l\vert \gamma$.

\smallpagebreak

Similarly we can rewrite the conditions $(G^\prime)$ and
$(B^\prime)$ respectively in the form

\noindent {\bf (B): The conditions of $b$-series:}

\noindent {\bf (BG): The general conditions of $b$-series:}
$$
p_1-\mu q_1\equiv 0\mod {(2/\gamma_2)(b_1/\gamma_b)c},
\tag{3.2.41}
$$
$$
\split &p_1-\mu q_1\not\equiv 0\mod
{(2/\gamma_2)(b_1/\gamma_b)c\,l}\\
&\text{for any prime $l$ such that}\ l^2\vert a_1\ \text{and}\
\text{g.c.d}(l,\gamma)=1, \endsplit \tag{3.2.42}
$$
$$
l \not\vert p_1\ \text{for any prime $l$ such that}\ l^2\vert b_1\
\text{and}\ \text{g.c.d}(l,\gamma)=1. \tag{3.2.43}
$$

\noindent {\bf (BS) The singular conditions of $b$-series:}
$$
\widetilde{a}_1q_1\equiv 0\mod \gamma_a, \tag{3.2.44}
$$
$$
q_1(p_1-\mu q_1)\equiv 0 \mod {(2/\gamma_2)\gamma_a}, \tag{3.2.45}
$$
$$
\gamma_2p_1q_1 \equiv \gamma_2q_1^2 \equiv 0\mod\gamma_a,
\tag{3.2.46}
$$
$$
(\delta-\gamma \mu^2) q_1^2\equiv 0 \mod
{(4/\gamma_2)a_1c^2\gamma_a}. \tag{3.2.47}
$$
$$
\cases
(a_1/\gamma_a)(\gamma p_1^2+\delta q_1^2)\equiv 0\mod 2\gamma_a l\\
(a_1/\gamma_a)\gamma_2\gamma_b p_1q_1\equiv 0\mod l\\
\gamma^2 p_1^2+\gamma\delta q_1^2\equiv 2\gamma^2\mu p_1q_1 \mod
4b_1c^2\gamma_a^2 l
\endcases
\tag{3.2.48}
$$
is not satisfied for any prime $l$ such that $l^2\vert a_1b_1$ and
$l\vert \gamma$,
$$
\cases
(a_1/\gamma_a)(\gamma p_1^2+\delta q_1^2)\equiv 0\mod 2\gamma_a l\\
(a_1/\gamma_a) p_1q_1\equiv 0\mod \gamma_a l\\
\mu \gamma^2 p_1^2+\mu \gamma\delta q_1^2\equiv 2\delta\gamma
p_1q_1 \mod 4b_1c^2\gamma_a^2 l
\endcases
\tag{3.2.49}
$$
is not satisfied for any prime $l\,\vert\, (2a_1b_1/\gamma)$ and
$l\vert \gamma$.

\smallpagebreak

We remind that here
$$
\mu\in (\bz/2a_1b_1c^2/\gamma)^\ast,\tag{3.2.50}
$$
$$
\delta\equiv \gamma\mu^2\mod{4a_1b_1c^2/\gamma},\tag{3.2.51}
$$
$$
\gamma p_1^2-\delta q_1^2=\pm 2(2/\gamma_2)(a_1/\gamma_a)\gamma_bc
\tag{3.2.52}
$$
for the $a$-series, and
$$
\gamma p_1^2-\delta q_1^2=\pm 2(2/\gamma_2)(b_1/\gamma_b)\gamma_ac
\tag{3.2.53}
$$
for the $b$-series.

\subhead
3.3. The resolution of the singular conditions (AS) and (BS)
\endsubhead
Here we completely resolve the singular conditions (AS) and (BS)
(assuming the corresponding general conditions \thetag{3.2.50},
\thetag{3.2.51}, \thetag{3.2.52}, \thetag{3.2.53}, \thetag{3.2.32},
\thetag{3.2.41} of
these series). It makes all our results very effective.

Below we consider the singular condition (AS) of $a$-series.
The singular condition (AS) consists of congruences and non-congruences
$\mod \gamma$. We denote by (AS)$^{(p)}$ the corresponding conditions
over the prime number $p\vert \gamma$. It is enough to satisfy all
conditions (AS)$^{(p)}$ for all $p\vert \gamma$. Below we consider
several cases which all together cover all possible ones. For a prime
$p$ and a natural number $n$, we denote as $n^{(p)}=p^{\nu_p}$ the
$p$-component of $n$.
That is $n^{(p)}=p^{\nu_p(n)} \vert n$, and $\text{g.c.d}(n^{(p)},
n/n^{(p)})=1$.

\smallpagebreak

There are several cases which we consider below.

\smallpagebreak

{\it Case $\gamma_2=2$, $\gamma_a\equiv \gamma_b\equiv 1\mod 2$.}
Then $c$, $a_1$, $b_1$ are odd and (AS)$^{(2)}$ obviously satisfies.

\smallpagebreak

{\it Case $\gamma_2=2$, $2\vert \gamma_a$.} Then
$\gamma^{(2)}=2a_1^{(2)}\ge 4$, $\gamma_b^{(2)}=1$ and $b_1\equiv 1\mod 2$.

By \thetag{3.2.51}, we obtain that
$$
\delta\equiv 0\mod 2. \tag{3.3.1}
$$
All conditions (AS) trivially satisfy over $2$ except
\thetag{3.2.39} for $l=2$ which gives $\delta q_1^2\not\equiv 0\mod
4$ if $4\vert a_1$. It is equivalent to $\delta\equiv 2\mod 4$ and
$q_1\equiv 1\mod 2$ if $4\vert a_1$. This follows from
\thetag{3.2.52} since $4\vert \gamma$. Thus (AS) over $2$ satisfies.

{\it Case $\gamma_2=2$, $2\vert \gamma_b$.} Then
$\gamma_a^{(2)}=1$ and $\gamma^{(2)}=2\gamma_b^{(2)}\ge 4$,
$b_1^{(2)}=\gamma_b^{(2)}$.  Denote $\gamma^{(2)}=2^{t}$, $t\ge
2$. For $l=2$, the condition \thetag{3.2.40} satisfies, and
\thetag{3.2.35}---\thetag{3.2.39} give over $2$ respectively
$$
q_1\equiv 0\mod 2^{[t/2]}, \tag{3.3.2}
$$
$$
q_1(p_1-\mu q_1)\equiv 0\mod 2^{t-1}, \tag{3.3.3}
$$
$$
p_1q_1\equiv q_1^2\equiv 0\mod 2^{t-2},\tag{3.3.4}
$$
$$
(\delta-\gamma\mu^2)q_1^2\equiv 0\mod 2^{2t-1}, \tag{3.3.5}
$$
$$
\gamma p_1^2+\delta q_1^2\not\equiv 0\mod 2^{t+1}\ if \ t\ge
3.\tag{3.3.6}
$$
By \thetag{3.2.51}, $\delta$ is even.

Assume $t=2$. Then \thetag{3.3.2}---\thetag{3.3.6} are equivalent to
$q_1\equiv 0\mod 2$. By \thetag{3.2.52}, then $p_1\equiv 1\mod 2$.

Assume that $t$ is even and $t\ge 4$. By \thetag{3.3.2}, we have
$q_1\equiv 0\mod 2^{t/2}$. Then \thetag{3.3.6} is equivalent to
$p_1\equiv 1\mod 2$. By \thetag{3.3.3}, then $q_1\equiv 0\mod
2^{t-1}$. It follows \thetag{3.3.5}. Thus,
\thetag{3.3.2}---\thetag{3.3.6} are equivalent to $p_1\equiv 1\mod 2$
and $q_1\equiv 0\mod 2^{t-1}$. We had the same for $t=2$.

Assume that $t$ is odd and $t\ge 3$. Let us suppose that
$p_1\equiv 0\mod 2$. Then \thetag{3.3.6} gives $\delta
q_1^2\not\equiv 0\mod 2^{t+1}$. By \thetag{3.3.2}, we have
$q_1\equiv 0\mod 2^{(t-1)/2}$. Moreover $\delta$ is even. Then
\thetag{3.3.6} is equivalent to $\delta\equiv 2\mod 4$ and and
$q_1\equiv 2^{(t-1)/2}\mod 2^{(t-1)/2+1}$. Then $(\delta-\gamma
\mu^2)q_1^2\equiv 2^t\mod 2^{t+1}$ and \thetag{3.3.5} is not valid.
This shows that $p_1\equiv 1\mod 2$. By \thetag{3.3.2}, $q_1$ is
even. By \thetag{3.3.3}, then $q_1\equiv 0\mod 2^{t-1}$. These imply
all conditions \thetag{3.3.2}---\thetag{3.3.6}. Thus, we obtain the
same conditions $p_1$ is odd and $q_1\equiv 0\mod 2^{t-1}$. Thus,
in this case, the condition (AS) over $2$ is
$$
If\ \gamma_2=2\ and\ 2\vert b_1\, ,\ then\ p_1\equiv 1\mod 2\ and\
q_1\equiv 0\mod \gamma^{(2)}/2. \tag{3.3.7}
$$

\smallpagebreak

{\it Case  $\gamma_2=1$, $2\vert \gamma$.}

{\it Case $\gamma_2=1$, $2\vert \gamma$ and $2\vert a_1$.} Then
$2\le \gamma^{(2)}=\gamma_a^{(2)}\vert a_1$, $\gamma_b^{(2)}=1$
and $b_1\equiv 1\mod 2$.

By \thetag{3.2.50}, \thetag{3.2.51} and \thetag{3.2.32}, we get respectively
$$
\mu\equiv 1\mod 2,   \tag{3.3.8}
$$
$$
\delta\equiv \gamma\mu^2\mod {4(a_1^{(2)}/\gamma_a^{(2)})},
\tag{3.3.9}
$$
and then $\delta\equiv \gamma \mod 4$,
$$
p_1 -\mu q_1\equiv 0\mod {2(a_1^{(2)}/\gamma_a^{(2)}}).\tag{3.3.10}
$$
It follows that over 2 all
conditions (AS) satisfy except \thetag{3.2.39} and \thetag{3.2.40}
which give respectively
$$
\gamma p_1^2+\delta q_1^2\equiv 2\gamma \mu p_1q_1\mod {8
(a_1^{(2)}/\gamma_a^{(2)})\tag{3.3.11}
$$
is not satisfied if $4\vert a_1$, and
$$
\cases
p_1q_1\equiv 0\mod 2\\
\mu \gamma p_1^2+\mu \delta q_1^2\equiv 2\delta p_1q_1\mod
{8(a_1^{(2)}/\gamma_a^{(2)})}}
\endcases
\tag{3.3.12}
$$
is not satisfied.

By \thetag{3.3.10}, we have $\gamma p_1^2+\gamma \mu q_1^2\equiv
2\gamma \mu p_1q_1\mod 4 {a_1 (a_1^{(2)}/\gamma_a^{(2)})}$. Since
$a_1$ is even, it follows
$$
\gamma p_1^2+\gamma \mu q_1^2\equiv 2\gamma \mu p_1q_1\mod 8
{(a_1^{(2)}/\gamma_a^{(2)})}.\tag{3.3.13}
$$
Then \thetag{3.3.11} is equivalent to
$(\delta-\gamma\mu^2)q_1^2\not\equiv 0\mod
{8(a_1^{(2)}/\gamma_a^{(2)})}$ if $4\vert a_1$. By \thetag{3.3.9},
this is equivalent to
$$
q_1\equiv 1\mod 2\ and \ \delta-\gamma\mu^2\not\equiv 0\mod
{8(a_1^{(2)}/\gamma_a^{(2)})} \tag{3.3.14}
$$
if $4\vert a_1$.

Similarly, one can see that \thetag{3.3.12} is equivalent to
$$
\cases p_1q_1\equiv 0\mod 2\\
(\delta-\gamma\mu^2)q_1^2\equiv 0\mod {8(a_1^{(2)}/\gamma_a^{(2)})}
\endcases
$$
is not valid. By above relations, it is equivalent to $p_1\equiv
1\mod 2$. Thus, in this case, (AS) over 2 is
equivalent to two conditions:
$$
if\  2\vert \gamma,\ \gamma_2=1\ and\ 2\vert a_1,\ then\ p_1\equiv
1\mod 2\tag{3.3.15}
$$
and
$$
if\ 2\vert \gamma,\ \gamma_2=1\ and\ 4\vert a_1,\ then\
\delta-\gamma\mu^2\not\equiv 0\mod {(8a_1b_1c^2/\gamma)}.
\tag{3.3.16}
$$

\smallpagebreak

{\it Case $\gamma_2=1$, $2\vert \gamma$ and $2\vert b_1$.} Then
$2\le \gamma^{(2)}=\gamma_b^{(2)}\vert b_1$, $\gamma_a^{(2)}=1$
and $a_1\equiv 1\mod 2$.

By \thetag{3.2.50}, \thetag{3.2.51} and \thetag{3.2.32}, we get
respectively
$$
\mu\equiv 1\mod 2 , \tag{3.3.17}
$$
$$
\delta\equiv \gamma\mu^2\mod {4(b_1^{(2)}/\gamma_b^{(2)})},
\tag{3.3.18}
$$
and then $\delta\equiv \gamma \mod 4$,
$$
p_1 -\mu q_1\equiv 0\mod 2. \tag{3.3.19}
$$

Over $2$, the conditions \thetag{3.2.35}---\thetag{3.2.40} give
respectively
$$
\widetilde{b}_1q_1\equiv 0\mod \gamma_b^{(2)}, \tag{3.3.20}
$$
$$
q_1(p_1-\mu q_1)\equiv 0 \mod {2\gamma_b^{(2)}, \tag{3.3.21}
$$
$$
p_1q_1 \equiv q_1^2 \equiv 0\mod\gamma_b^{(2)}, \tag{3.3.22}
$$
$$
(\delta-\gamma \mu^2) q_1^2\equiv 0 \mod
{4(\gamma_b^{(2)})^2(b_1^{(2)}/\gamma_b^{(2)})}, \tag{3.3.23}
$$
$$
\gamma p_1^2+ \delta q_1^2\equiv 2\gamma \mu p_1q_1 \mod
{8\gamma_b^{(2)}}  \tag{3.3.24}
$$
is not satisfied if $4\vert b_1$,
$$
\cases
(b_1/\gamma_b) p_1q_1\equiv 0\mod 2\gamma_b^{(2)}\\
\mu \gamma p_1^2+\mu \delta q_1^2\equiv 2\delta p_1q_1 \mod
{8\gamma_b^{(2)}}} \endcases \tag{3.3.25}
$$
is not satisfied.

By \thetag{3.3.23}, we have $\delta q_1^2\equiv \gamma\mu^2q_1^2\mod
{8\gamma_b^{(2)}}$. It follows that \thetag{3.3.24} is equivalent to
$$
p_1-\mu q_1\not\equiv 0\mod 4 \tag{3.3.26}
$$
if $4\vert b_1$. Similarly \thetag{3.3.25} is equivalent to
$$
\cases
(b_1/\gamma_b) p_1q_1\equiv 0\mod 2\gamma_b^{(2)}\\
p_1-\mu q_1\equiv 0\mod 4
\endcases \tag{3.3.27}
$$
is not satisfied.

Assume that
$$
(b_1/\gamma_b) p_1q_1\not\equiv 0\mod 2\gamma_b^{(2)}. \tag{3.3.28}
$$
By \thetag{3.3.22}, we have $p_1q_1\equiv 0\mod{\gamma_b^{(2)}}$.
Then \thetag{3.3.28} is equivalent to $(b_1/\gamma_b)$ is odd and
$p_1q_1\equiv \gamma_b^{(2)}\mod 2\gamma_b^{(2)}$. By
\thetag{3.3.21}, then $q_1^2\equiv \gamma_b^{(2)}\mod
2\gamma_b^{(2)}$. Then $\gamma_b^{(2)}$ is a square and $4\vert
b_1$. By \thetag{3.3.26}, then $p_1-\mu q_1\not\equiv 0\mod 4$.
Thus, \thetag{3.3.26} and \thetag{3.3.27} are equivalent to
$$
p_1-\mu q_1\equiv 2\mod 4 \tag{3.3.29}
$$
since $p_1\equiv q_1\equiv 0\mod 2$ by \thetag{3.3.17},
\thetag{3.3.19} and \thetag{3.3.22}. By \thetag{3.3.21}, then $q_1\equiv
0\mod \gamma_b^{(2)}$, and all other conditions of (AS) over $2$
follow from here.

Thus, in this case, (AS) over $2$ is equivalent to
$$
if\ 2\vert \gamma,\ \gamma_2=1,\ and\ 2\vert b_1,\ then\ p_1-\mu
q_1\not\equiv 0\mod 4\ and\ q_1\equiv 0\mod {\gamma_b^{(2)}}.
\tag{3.3.30}
$$

{\it Case a prime odd $l\vert \gamma$ and $l\vert a_1$.} Then
$l\le \gamma^{(l)}=\gamma_a^{(l)}\vert a_1$, $\gamma_b^{(l)}=1$
and $b_1\not\equiv 0\mod l$.

By \thetag{3.2.50}, \thetag{3.2.51}, \thetag{3.2.32} we have respectively
$$
\mu\in (\bz/(a_1^{(l)}/\gamma_a^{(l)}))^\ast,  \tag{3.3.31}
$$
$$
\delta\equiv \gamma\mu^2\mod {(a_1^{(l)}/\gamma_a^{(l)})},
\tag{3.3.32}
$$
$$
p_1 -\mu q_1\equiv 0\mod {(a_1^{(l)}/\gamma_a^{(l)})}.\tag{3.3.33}
$$
Conditions (AS) satisfy over $l$ except \thetag{3.2.39} and
\thetag{3.2.40} which give respectively
$$
\cases
\delta q_1^2\equiv 0\mod l\\
\gamma p_1^2+\delta q_1^2\equiv 2\gamma \mu p_1q_1 \mod
{(a_1^{(l)}/\gamma_a^{(l)})l}
\endcases
\tag{3.3.34}
$$
is not satisfied if $l^2\vert a_1$,
$$
\cases \delta q_1^2\equiv 0\mod l\\
p_1q_1\equiv 0\mod l\\
\mu \gamma p_1^2+\mu \delta q_1^2\equiv 2\delta p_1q_1 \mod
{(a_1^{(l)}/\gamma_a^{(l)}) l}
\endcases
\tag{3.3.35}
$$
is not satisfied if $l\,\vert\, (a_1^{(l)}/\gamma_a^{(l)})$.

Taking square of \thetag{3.3.33}, we obtain
$$
\gamma p_1^2+\gamma \mu^2 q_1^2\equiv 2\gamma\mu p_1q_1\mod
(a_1^{(l)}/\gamma_a^{(l)})l. \tag{3.3.36}
$$
This shows that \thetag{3.3.34} is equivalent to
$$
\cases
\delta q_1^2\equiv 0\mod l\\
(\delta-\gamma \mu^2) q_1^2\equiv 0 \mod
{(a_1^{(l)}/\gamma_a^{(l)})l}
\endcases
\tag{3.3.37}
$$
is not satisfied if $l^2\vert a_1$. By \thetag{3.3.33}, this is
equivalent to
$$
\split &If\ odd\ prime\ l\vert \gamma\ and\ l^2\vert a_1,\
then\
q_1\not\equiv 0\mod l\ and\\
&either\ \delta\not\equiv 0\mod l\
or\ (\delta-\gamma \mu^2)\not\equiv 0\mod
(a_1^{(l)}/\gamma_a^{(l)})l.
\endsplit
\tag{3.3.38}
$$

If $l\,\vert\, (a_1^{(l)}/\gamma_a^{(l)})$, then $l^2\vert a_1$ and
\thetag{3.3.38} is valid. Then $q_1\not\equiv 0\mod l$. By
\thetag{3.3.31} and \thetag{3.3.33}, then $p_1\not\equiv 0\mod l$.
Thus \thetag{3.3.35} follows from \thetag{3.3.38}.

Finally we get that (AS) over $l$ is equivalent to \thetag{3.3.38}
in this case.

{\it Case a prime odd $l\vert \gamma$ and $l\vert b_1$.} Then
$l\le \gamma^{(l)}=\gamma_b^{(l)}\vert b_1$, $\gamma_a^{(l)}=1$
and $a_1\not\equiv 0\mod l$.

By \thetag{3.2.50}, \thetag{3.2.51}, \thetag{3.2.32} we have respectively
$$
\mu\in (\bz/(b_1^{(l)}/\gamma_b^{(l)}))^\ast, \tag{3.3.39}
$$
$$
\delta\equiv \gamma\mu^2\mod {(b_1^{(l)}/\gamma_b^{(l)})},
\tag{3.3.40}
$$
$$
\gamma p_1^2-\delta q_1^2=\pm 2(2/\gamma_2)(a_1/\gamma_a)\gamma_bc
\tag{3.3.41}
$$
Over $l$, the conditions \thetag{3.2.35}---\thetag{3.2.40} give
respectively
$$
\widetilde{b}_1q_1\equiv 0\mod \gamma_b^{(l)}, \tag{3.3.42}
$$
$$
q_1(p_1-\mu q_1)\equiv 0 \mod {\gamma_b^{(l)}}, \tag{3.3.43}
$$
$$
p_1q_1 \equiv q_1^2 \equiv 0\mod\gamma_b^{(l)}, \tag{3.3.44}
$$
$$
(\delta-\gamma \mu^2) q_1^2\equiv 0 \mod
{(b_1^{(l)}/\gamma_b^{(l)})(\gamma_b^{(l)})^2}, \tag{3.3.45}
$$
$$
\cases
(b_1^{(l)}/\gamma_b^{(l)})(\gamma p_1^2+\delta q_1^2)\equiv 0\mod 
{\gamma_b^{(l)} l}\\
\gamma p_1^2+\delta q_1^2\equiv 2\gamma\mu p_1q_1 \mod
{\gamma_b^{(l)} l}
\endcases
\tag{3.3.46}
$$
is not satisfied if $l^2\vert b_1$,
$$
\cases
(b_1^{(l)}/\gamma_b^{(l)})(\gamma p_1^2+\delta q_1^2)\equiv 0\mod 
{\gamma_b^{(l)} l}\\
\mu \gamma p_1^2+\mu \delta q_1^2\equiv 2\delta p_1q_1 \mod
{\gamma_b^{(l)} l}
\endcases
\tag{3.3.47}
$$
is not satisfied if $l\,\vert\, (b_1^{(l)}/\gamma_b^{(l)})$.

By \thetag{3.3.45}, we have $\delta q_1^2\equiv \gamma\mu^2q_1^2\mod
\gamma_b l$. Then \thetag{3.3.46} is equivalent to
$$
\cases
(b_1^{(l)}/\gamma_b^{(l)})(p_1^2+\mu^2 q_1^2)\equiv 0\mod {l}\\
p_1-\mu q_1\equiv 0 \mod {l}
\endcases
\tag{3.3.48}
$$
is not satisfied if $l^2\vert b_1$. Similarly (using also
\thetag{3.3.39}) one can see that \thetag{3.3.47} is equivalent to
$$
\cases
(b_1^{(l)}/\gamma_b^{(l)})(p_1^2+\mu^2 q_1^2)\equiv 0\mod {l}\\
p_1-\mu q_1 \equiv 0 \mod {l}
\endcases
\tag{3.3.49}
$$
is not satisfied if $l\,\vert\, (b_1^{(l)}/\gamma_b^{(l)})$. Thus,
\thetag{3.3.48} implies \thetag{3.3.49}, and it is enough to satisfy
\thetag{3.3.48}.

Assume that $p_1-\mu q_1\not\equiv 0\mod l$. By \thetag{3.3.43} we
get $q_1\equiv 0\mod \gamma_b^{(l)}$, and all other conditions
follow.

Assume that $p_1-\mu q_1\equiv 0\mod l$. Then
$p_1^2+\mu^2q_1^2\equiv 2p_1^2\mod l$, and \thetag{3.3.48} implies
that $p_1\not\equiv 0\mod l$. By \thetag{3.3.44}, $q_1\equiv 0\mod
l$, and we get a contradiction. Since $q_1\equiv 0\mod l$, the
condition $p_1\not\equiv \mu q_1\mod l$ can be replaced by
$p_1\not\equiv 0\mod l$.

Thus, in this case, the condition (AS) over $l$ is equivalent to
two conditions
$$
If\ odd\ prime\ l\vert \gamma\ and\ l\vert b_1,\ then\ q_1\equiv 0\mod
{\gamma_b^{(l)}}.
\tag{3.3.50}
$$
$$
If\ odd\ prime\ l\vert \gamma\ and\ l^2\vert b_1,\ then\
p_1\not\equiv 0\mod l. \tag{3.3.51}
$$

Thus, we obtain

\proclaim{Theorem 3.3.1} The singular condition (AS)
is equivalent to
$$
\split &if\ odd\ prime\ l\vert \gamma\ and\ l^2\vert a_1,\
then\
q_1\not\equiv 0\mod l\ and\\
&either\ \delta\not\equiv 0\mod l\
or\ (\delta-\gamma \mu^2)\not\equiv 0\mod
(a_1^{(l)}/\gamma_a^{(l)})l;\\
&if\ odd\ prime\ l\vert \gamma\ and\ l\vert b_1,\ then\ q_1\equiv 0\mod
{\gamma_b^{(l)}};\\
&if\ odd\ prime\ l\vert \gamma\ and\ l^2\vert b_1,\  
then\ p_1 \not\equiv 0\mod l;\\
&if\  2\vert \gamma,\ \gamma_2=1\ and\ 2\vert a_1,\ then\ p_1\equiv
1\mod 2;\\
&if\ 2\vert \gamma,\ \gamma_2=1\ and\ 4\vert a_1,\ then\
\delta-\gamma\mu^2\not\equiv 0\mod {(8a_1b_1c^2/\gamma)};\\
&if\ 2\vert \gamma,\ \gamma_2=1,\ and\ 2\vert b_1,\ then\ p_1-\mu
q_1\not\equiv 0\mod 4\ and\ q_1\equiv 0\mod {\gamma_b^{(2)}};\\
&if\ 2\vert \gamma,\ \gamma_2=2\ and\ 2\vert b_1\, ,\
then\ p_1\equiv 1\mod 2\ and\
q_1\equiv 0\mod \gamma^{(2)}/2.
\endsplit
\tag{3.3.52}
$$

The singular condition (BS) is equivalent to
$$
\split
&if\ odd\ prime\ l\vert \gamma\ and\ l\vert a_1,\ then\ q_1\equiv 0\mod
{\gamma_a^{(l)}};\\
&if\ odd\ prime\ l\vert \gamma\ and\ l^2\vert a_1,\  
then\ p_1 \not\equiv 0\mod l;\\
&if\ odd\ prime\ l\vert \gamma\ and\ l^2\vert b_1,\
then\  q_1\not\equiv 0\mod l\ and\\
&either\ \delta\not\equiv 0\mod l\
or\ (\delta-\gamma \mu^2)\not\equiv 0\mod
(b_1^{(l)}/\gamma_b^{(l)})l;\\
&if\ 2\vert \gamma,\ \gamma_2=1,\ and\ 2\vert a_1,\ then\ p_1-\mu
q_1\not\equiv 0\mod 4\ and\ q_1\equiv 0\mod {\gamma_a^{(2)}};\\
&if\  2\vert \gamma,\ \gamma_2=1\ and\ 2\vert b_1,\ then\ p_1\equiv
1\mod 2;\\
&if\ 2\vert \gamma,\ \gamma_2=1\ and\ 4\vert b_1,\ then\
\delta-\gamma\mu^2\not\equiv 0\mod {(8a_1b_1c^2/\gamma)};\\
&if\ 2\vert\gamma,\ \gamma_2=2\ and\ 2\vert a_1\, ,\
then\ p_1\equiv 1\mod 2\ and\
q_1\equiv 0\mod \gamma^{(2)}/2.
\endsplit
\tag{3.3.53}
$$
\endproclaim

Thus, we can finally rewrite the conditions (A) of $a$-series, and
the conditions (B) of $b$-series in the very efficient form which makes
all our results very effective.

\smallpagebreak

\noindent {\bf (A): The conditions of $a$-series:}

\noindent {\bf (AG): The general conditions of $a$-series:}
$$
p_1-\mu q_1\equiv 0\mod {(2/\gamma_2)(a_1/\gamma_a)c},
\tag{3.3.54}
$$
$$
\split &p_1-\mu q_1\not\equiv 0\mod
{(2/\gamma_2)(a_1/\gamma_a)c\,l}\\
&\text{for any prime $l$ such that}\ l^2\vert b_1\ \text{and}\
\text{g.c.d}(l,\gamma)=1, \endsplit \tag{3.3.55}
$$
$$
l \not\vert p_1\ \text{for any prime $l$ such that}\ l^2\vert a_1\
\text{and}\ \text{g.c.d}(l,\gamma)=1. \tag{3.3.56}
$$

\noindent {\bf (AS) The singular conditions of $a$-series:}
$$
\split &if\ odd\ prime\ l\vert \gamma\ and\ l^2\vert a_1,\
then\
q_1\not\equiv 0\mod l\ and\\
&either\ \delta\not\equiv 0\mod l\
or\ (\delta-\gamma \mu^2)\not\equiv 0\mod
(a_1^{(l)}/\gamma_a^{(l)})l;\\
&if\ odd\ prime\ l\vert \gamma\ and\ l\vert b_1,\ then\ q_1\equiv 0\mod
{\gamma_b^{(l)}};\\
&if\ odd\ prime\ l\vert \gamma\ and\ l^2\vert b_1,\  
then\ p_1 \not\equiv 0\mod l;\\ 
&if\  2\vert \gamma,\ \gamma_2=1\ and\ 2\vert a_1,\ then\ p_1\equiv
1\mod 2;\\
&if\ 2\vert \gamma,\ \gamma_2=1\ and\ 4\vert a_1,\ then\
\delta-\gamma\mu^2\not\equiv 0\mod {(8a_1b_1c^2/\gamma)};\\
&if\ 2\vert \gamma,\ \gamma_2=1,\ and\ 2\vert b_1,\ then\ p_1-\mu
q_1\not\equiv 0\mod 4\ and\ q_1\equiv 0\mod {\gamma_b^{(2)}};\\
&if\ 2\vert \gamma,\ \gamma_2=2\ and\ 2\vert b_1\, ,\
then\ p_1\equiv 1\mod 2\ and\
q_1\equiv 0\mod \gamma^{(2)}/2.
\endsplit
\tag{3.3.57}
$$

\smallpagebreak

\noindent {\bf (B): The conditions of $b$-series:}

\noindent {\bf (BG): The general conditions of $b$-series:}
$$
p_1-\mu q_1\equiv 0\mod {(2/\gamma_2)(b_1/\gamma_b)c},
\tag{3.3.58}
$$
$$
\split &p_1-\mu q_1\not\equiv 0\mod
{(2/\gamma_2)(b_1/\gamma_b)c\,l}\\
&\text{for any prime $l$ such that}\ l^2\vert a_1\ \text{and}\
\text{g.c.d}(l,\gamma)=1, \endsplit \tag{3.3.59}
$$
$$
l \not\vert p_1\ \text{for any prime $l$ such that}\ l^2\vert b_1\
\text{and}\ \text{g.c.d}(l,\gamma)=1. \tag{3.3.60}
$$

\noindent {\bf (BS) The singular conditions of $b$-series:}

$$
\split
&if\ odd\ prime\ l\vert \gamma\ and\ l\vert a_1,\ then\ q_1\equiv 0\mod
{\gamma_a^{(l)}};\\
&if\ odd\ prime\ l\vert \gamma\ and\ l^2\vert a_1,\  
then\ p_1\not\equiv 0\mod l;\\
&if\ odd\ prime\ l\vert \gamma\ and\ l^2\vert b_1,\
then\  q_1\not\equiv 0\mod l\ and\\
&either\ \delta\not\equiv 0\mod l\
or\ (\delta-\gamma \mu^2)\not\equiv 0\mod
(b_1^{(l)}/\gamma_b^{(l)})l;\\
&if\ 2\vert \gamma,\ \gamma_2=1,\ and\ 2\vert a_1,\ then\ p_1-\mu
q_1\not\equiv 0\mod 4\ and\ q_1\equiv 0\mod {\gamma_a^{(2)}};\\
&if\  2\vert \gamma,\ \gamma_2=1\ and\ 2\vert b_1,\ then\ p_1\equiv
1\mod 2;\\
&if\ 2\vert \gamma,\ \gamma_2=1\ and\ 4\vert b_1,\ then\
\delta-\gamma\mu^2\not\equiv 0\mod {(8a_1b_1c^2/\gamma)};\\
&if\ 2\vert\gamma,\ \gamma_2=2\ and\ 2\vert a_1\, ,\
then\ p_1\equiv 1\mod 2\ and\
q_1\equiv 0\mod \gamma^{(2)}/2.
\endsplit
\tag{3.3.61}
$$

\smallpagebreak

We remind that here $\gamma\vert 2a_1b_1$ and
$$
\mu\in (\bz/2a_1b_1c^2/\gamma)^\ast,\tag{3.3.62}
$$
$$
\delta\equiv \gamma\mu^2\mod{4a_1b_1c^2/\gamma},\tag{3.3.63}
$$
$$
\gamma p_1^2-\delta q_1^2=\pm 2(2/\gamma_2)(a_1/\gamma_a)\gamma_bc
\tag{3.3.64}
$$
for the $a$-series, and
$$
\gamma p_1^2-\delta q_1^2=\pm 2(2/\gamma_2)(b_1/\gamma_b)\gamma_ac
\tag{3.3.65}
$$
for the $b$-series.

 \head 4.  
Final results for $\rho=2$
\endhead
Now we can formulate the main results which follow from Theorem
3.1.3 and the calculations above.

\proclaim{Theorem 4.1} Let $X$ be a K3 surface with $\rho
(X)=2$, and $H$ a polarization of $X$ of degree $H^2=2rs$ where
$r,\,s\in \bn$. Assume that the Mukai vector $(r,H,s)$ is primitive.
Let $Y$ be the moduli space of sheaves on $X$ with the isotropic
Mukai vector $v=(r,H,s)$. Let $\wH=H/d$ be the corresponding
primitive polarization, and $\wH\cdot N(X)=\gamma \bz$. We denote
by $\mu$ the invariant of the pair $(N(X),\wH)$ and use notations
of Proposition 3.1.1.

We have $Y\cong X$ if
$$
\text{g.c.d}(c,d\gamma )=1
$$
and $X$ belongs either to $a$-series or to $b$-series.

Here $X$ belongs to $a$-series if for one of $\epsilon=\pm 1$ the
equation
$$
\gamma p_1^2-\delta q_1^2=\epsilon 2
(2/\gamma_2)(a_1/\gamma_a)\gamma_b c. \tag{4.1}
$$
has an integral solution $(p_1,q_1)$ satisfying the conditions (A)
of $a$-series \thetag{3.3.54}---\thetag{3.3.57}.

These solutions $(p_1,q_1)$ of \thetag{4.1} give all solution
$(x,y)$ of Theorem 3.1.3 from $a$-series as associated solutions
$$
(x,y)=\pm \left({-2a_1b_1c\over \gamma} + {\epsilon b_1\gamma_2\gamma_a
p_1^2\over \gamma_b}\, , \, {\epsilon b_1\gamma_2\gamma_a
p_1q_1\over \gamma_b}\right). \tag{4.2}
$$

Here $X$ belongs to $b$-series if for one of signs $\epsilon=\pm
1$ the equation
$$
\gamma p_1^2-\delta q_1^2=\epsilon 2
(2/\gamma_2)(b_1/\gamma_b)\gamma_a c. \tag{4.3}
$$
has an integral solution $(p_1,q_1)$ satisfying the conditions (B)
of $b$-series \thetag{3.3.58}---\thetag{3.3.61}.

These solutions $(p_1,q_1)$ of \thetag{4.3} give all solutions
$(x,y)$ of Theorem 3.1.3 of $b$-series as associated solutions
$$
(x,y)=\pm \left({-2a_1b_1c\over \gamma} + {\epsilon a_1\gamma_2\gamma_b
p_1^2\over \gamma_a}\, , \, {\epsilon a_1\gamma_2\gamma_b
p_1q_1\over \gamma_a}\right). \tag{4.4}
$$

These conditions are necessary to have $Y\cong X$ if $X$ is a
general K3 surface with $\rho (X)=2$, i. e. the automorphism group
of the transcendental periods $(T(X),H^{2,0}(X))$ is $\pm 1$.
\endproclaim

Now we want to interpret solutions $(p_1,q_1)$ of Theorem 4.1 as
elements of the Picard lattice $N(X)$.

Let $(p_1,q_1)$ be a solution of Theorem 4.1 from $a$-series.
Then
$$
\gamma p_1^2-\delta q_1^2= \epsilon 2
(2/\gamma_2)(a_1/\gamma_a)\gamma_b c. \tag{4.5}
$$
By \thetag{3.3.54}, we have
$$
p_1-\mu q_1\equiv 0 \mod{(2/\gamma_2)(a_1/\gamma_a)c}. \tag{4.6}
$$
Let us put
$$
t=(b_1/\gamma_b)c. \tag{4.7}
$$
Then
$$
tp_1-\mu tq_1\equiv 0\mod 2a_1b_1c^2/\gamma. \tag{4.8}
$$
and
$$
\wth_1={t(p_1\wH+q_1f(\wH))\over 2a_1b_1c^2/\gamma}\in N(X).
\tag{4.9}
$$
We have
$$
\wth_1^2={t^2(\gamma p_1^2-\delta q_1^2)\over
2a_1b_1c^2/\gamma}=\epsilon 2b_1c \tag{4.10}
$$
and
$$
\wH\cdot \wth_1=\gamma(b_1/\gamma_b)c p_1 \equiv
0\mod{\gamma(b_1/\gamma_b)c}\ \text{and}\ p_1= {\wH\cdot
\wth_1\over \gamma (b_1/\gamma_b)c}. \tag{4.11}
$$
Also
$$
-f(\wH)\cdot \wth_1={b_1c \delta q_1\over \gamma_b}\ \text{and}\
q_1=-{f(\wH)\cdot \wth_1 \over \delta (b_1/\gamma_b)c}. \tag{4.12}
$$
Here
$$
-\gamma \delta =\det N(X). \tag{4.13}
$$

Another calculation of $p_1$ and $q_1$ is as follows. We have
$$
\wth_1={u\wH+vf(\wH)\over
2a_1b_1c^2/\gamma}={p_1\wH+q_1f(\wH)\over
(2/\gamma_2)(a_1/\gamma_a)c} \tag{4.14}
$$
where
$$
u\equiv 0\mod (b_1/\gamma_b)c,\ \text{and}\  p_1={u\over
(b_1/\gamma_b)c}, \tag{4.15}
$$
$$
v\equiv 0\mod (b_1/\gamma_b)c,\ \text{and}\ q_1={v\over
(b_1/\gamma_b)c}. \tag{4.16}
$$
We remind that here $\bz f(\wH)$ is the orthogonal complement to
$\wH$ in $N(X)$. Both these calculations of $p_1$ and $q_1$ are
equivalent.

By construction, \thetag{3.3.54} is equivalent to $\wth_1\in
N(X)$, \thetag{3.3.55} is equivalent to $\wth_1/l\not\in N(X)$,
\thetag{3.3.56} is equivalent to $\wH\cdot \wth_1 \not\equiv 0\mod
\gamma (b_1/\gamma_b)cl$.

Changing the letters $a$ and $b$ places we get the same results
for the $b$-series.

\smallpagebreak

Thus, we get

\proclaim{Theorem 4.2} Let $X$ be a K3 surface with $\rho
(X)=2$, and $H$ a polarization of $X$ of degree $H^2=2rs$ where
$r,\,s\in \bn$. Assume that the Mukai vector $(r,H,s)$ is primitive.
Let $Y$ be the moduli space of sheaves on $X$ with the isotropic
Mukai vector $v=(r,H,s)$. Let $\wH=H/d$ be the corresponding
primitive polarization, and $\wH\cdot N(X)=\gamma \bz$. We denote
by $\mu$ the invariant of the pair $(N(X),\wH)$ and use notations
of Proposition 3.1.1.

We have $Y\cong X$ if
$$
\text{g.c.d}(c,d\gamma )=1,
$$
and at least for one $\epsilon=\pm 1$ there exists $\wth_1\in
N(X)$ which belongs to the $a$-series or to the $b$-series
described below:

$\wth_1$ belongs to the $a$-series if
$$
\wth_1^2=\epsilon2b_1c\ and\ \wH\cdot \wth_1\equiv
0\mod{\gamma(b_1/\gamma_b)c}, \tag{4.17}
$$
$$
\wH\cdot \wth_1\not\equiv 0\mod{\gamma(b_1/\gamma_b)cl_1},\
\wth_1/l_2\not\in N(X) \tag{4.18}
$$
for any prime $l_1$ such that $l_1^2\vert a_1$ and
$\text{g.c.d}(l_1,\gamma)=1$, and for any prime $l_2$ such that
$l_2^2\vert b_1$ and $\text{g.c.d}(l_2,\gamma)=1$,
and
$$p_1= {\wH\cdot
\wth_1\over \gamma (b_1/\gamma_b)c},\ \ q_1=-{f(\wH)\cdot \wth_1
\over \delta (b_1/\gamma_b)c} \tag{4.19}
$$
satisfy the singular condition (AG) (conditions \thetag{3.3.57}
$\mod \gamma$) of $a$-series.

$\wth_1$ belongs to the $b$-series if
$$
\wth_1^2=\epsilon2a_1c\ and\ \wH\cdot \wth_1\equiv
0\mod{\gamma(a_1/\gamma_a)c}, \tag{4.20}
$$
$$
\wH\cdot \wth_1\not\equiv 0\mod{\gamma(a_1/\gamma_a)cl_1},\
\wth_1/l_2\not\in N(X) \tag{4.21}
$$
for any prime $l_1$ such that $l_1^2\vert b_1$ and
$\text{g.c.d}(l_1,\gamma)=1$, and
for any prime $l_2$ such that $l_2^2\vert
a_1$ and $\text{g.c.d}(l_2,\gamma)=1$,
and
$$p_1= {\wH\cdot
\wth_1\over \gamma (a_1/\gamma_a)c},\ \ q_1=-{f(\wH)\cdot \wth_1
\over \delta (a_1/\gamma_a)c} \tag{4.22}
$$
satisfy the singular condition (BG) (conditions \thetag{3.3.61}
$\mod \gamma$) of $b$-series.

These conditions are necessary to have $Y\cong X$ if $X$ is a
general K3 surface with $\rho (X)=2$, i. e. the automorphism group
of the transcendental periods $(T(X),H^{2,0}(X))$ is $\pm 1$.
\endproclaim

\remark{Important Remark 4.3} Applying Theorem 3.1.3 and the
formulae \thetag{4.2}, \thetag{4.4} for the associated
solution, we get the following formulae in terms of $X$ for the canonical
primitive $nef$ element $\wth$ of $Y$ defined by $(-a,0,b)\mod \bz v$\ :
$$
\wth^\prime = \cases {-\wH\over c}+{\epsilon(\wH\cdot
\wth_1)\wth_1\over b_1c^2}
&\text{if $\wth_1$ is from $a$-series,}\\
{-\wH\over c}+{\epsilon(\wH\cdot \wth_1)\wth_1\over a_1c^2}
&\text{if $\wth_1$ is from $b$-series}\\
\endcases
\tag{4.23}
$$
belongs to $N(X)$ and
$$
(Y,\wth)\cong (X,\pm w(\wth^\prime ))\ \text{for some}\ w\in
W^{(-2)}(N(X)). \tag{4.24}
$$
\endremark

Specializing  (by Lemma 2.2.1) the theorem 4.2, we get the following
sufficient condition to have $Y\cong X$ which is valid for $X$ with any $\rho
(X)$. This is one of the main results of the paper.

In Theorem 4.4 below, for $\wH\in N$ we apply the same
definitions and notations: $f(\wH)$, $\delta$, $\mu$, as for
$\wH\in N=N(X)$ of Proposition 3.1.1.

\proclaim{Theorem 4.4} Let $X$ be a K3 surface and $H$ a
polarization of $X$ of degree $H^2=2rs$ where $r,\,s\in \bn$. Assume that
the Mukai vector $(r,H,s)$ is primitive. Let $Y$ be the moduli
space of sheaves on $X$ with the isotropic Mukai vector
$v=(r,H,s)$. Let $\wH=H/d$ be the corresponding primitive polarization.

We have $Y\cong X$ if there exists $\wth_1\in N(X)$ such that
$\wH$, $\wth_1$ belong to a 2-dimensional primitive sublattice
$N\subset N(X)$ such that $\wH\cdot N=\gamma \bz$, $\gamma >0$,
and
$$
\text{g.c.d}(c,d\gamma)=1,
\tag{4.25}
$$
moreover, for one of $\epsilon=\pm 1$ the element $\wth_1$ belongs
to the $a$-series or to the $b$-series described below:

$\wth_1$ belongs to the $a$-series if
$$
\wth_1^2=\epsilon2b_1c\ and\ \wH\cdot \wth_1\equiv
0\mod{\gamma(b_1/\gamma_b)c},
\tag{4.26}
$$
$$
\wH\cdot \wth_1\not\equiv 0\mod{\gamma(b_1/\gamma_b)cl_1},\
\wth_1/l_2\not\in N(X)
\tag{4.27}
$$
for any prime $l_1$ such that $l_1^2\vert a_1$ and
$\text{g.c.d}(l_1,\gamma)=1$,
and any prime $l_2$ such that $l_2^2\vert b_1$ and
$\text{g.c.d}(l_2,\gamma)=1$,
and
$$p_1= {\wH\cdot
\wth_1\over \gamma (b_1/\gamma_b)c},\ \ q_1=-{f(\wH)\cdot \wth_1
\over \delta(b_1/\gamma_b)c} \tag{4.28}
$$
satisfy the singular condition (AS) of $a$-series:
$$
\split &if\ odd\ prime\ l\vert \gamma\ and\ l^2\vert a_1,\
then\
q_1\not\equiv 0\mod l\ and\\
&either\ \delta\not\equiv 0\mod l\
or\ (\delta-\gamma \mu^2)\not\equiv 0\mod
(a_1^{(l)}/\gamma_a^{(l)})l;\\
&if\ odd\ prime\ l\vert \gamma\ and\ l\vert b_1,\ then\ q_1\equiv 0\mod
{\gamma_b^{(l)}};\\
&if\ odd\ prime\ l\vert \gamma\ and\ l^2\vert b_1,\  
then\ p_1\not\equiv 0\mod l;\\
&if\  2\vert \gamma,\ \gamma_2=1\ and\ 2\vert a_1,\ then\ p_1\equiv
1\mod 2;\\
&if\ 2\vert \gamma,\ \gamma_2=1\ and\ 4\vert a_1,\ then\
\delta-\gamma\mu^2\not\equiv 0\mod {(8a_1b_1c^2/\gamma)};\\
&if\ 2\vert \gamma,\ \gamma_2=1,\ and\ 2\vert b_1,\ then\ p_1-\mu
q_1\not\equiv 0\mod 4\ and\ q_1\equiv 0\mod {\gamma_b^{(2)}};\\
&if\ 2\vert \gamma,\ \gamma_2=2\ and\ 2\vert b_1\, ,\
then\ p_1\equiv 1\mod 2\ and\
q_1\equiv 0\mod \gamma^{(2)}/2.
\endsplit
$$

$\wth_1$ belongs to the $b$-series if
$$
\wth_1^2=\epsilon2a_1c\ and\ \wH\cdot \wth_1\equiv
0\mod{\gamma(a_1/\gamma_a)c},
\tag{4.29}
$$
$$
\wH\cdot \wth_1\not\equiv 0\mod{\gamma(a_1/\gamma_a)cl_1},\
\wth_1/l_2\not\in N(X) \tag{4.30}
$$
for any prime $l_1$ such that
$l_1^2\vert b_1$ and $\text{g.c.d}(l_1,\gamma)=1$ and
any prime $l_2$ such that $l_2^2\vert a_1$ and
$\text{g.c.d}(l_2,\gamma)=1$,
and
$$p_1= {\wH\cdot
\wth_1\over \gamma (a_1/\gamma_a)c},\ \ q_1=-{f(\wH)\cdot \wth_1
\over \delta(a_1/\gamma_a)c} \tag{4.31}
$$
satisfy the singular condition (BS) of $b$-series:

$$
\split
&if\ odd\ prime\ l\vert \gamma\ and\ l\vert a_1,\ then\ q_1\equiv 0\mod
{\gamma_a^{(l)}};\\
&if\ odd\ prime\ l\vert \gamma\ and\ l^2\vert a_1,\  
then\ p_1\not\equiv 0\mod l;\\
&if\ odd\ prime\ l\vert \gamma\ and\ l^2\vert b_1,\
then\  q_1\not\equiv 0\mod l\ and\\
&either\ \delta\not\equiv 0\mod l\
or\ (\delta-\gamma \mu^2)\not\equiv 0\mod
(b_1^{(l)}/\gamma_b^{(l)})l;\\
&if\ 2\vert \gamma,\ \gamma_2=1,\ and\ 2\vert a_1,\ then\ p_1-\mu
q_1\not\equiv 0\mod 4\ and\ q_1\equiv 0\mod {\gamma_a^{(2)}};\\
&if\  2\vert \gamma,\ \gamma_2=1\ and\ 2\vert b_1,\ then\ p_1\equiv
1\mod 2;\\
&if\ 2\vert \gamma,\ \gamma_2=1\ and\ 4\vert b_1,\ then\
\delta-\gamma\mu^2\not\equiv 0\mod {(8a_1b_1c^2/\gamma)};\\
&if\ 2\vert\gamma,\ \gamma_2=2\ and\ 2\vert a_1\, ,\
then\ p_1\equiv 1\mod 2\ and\
q_1\equiv 0\mod \gamma^{(2)}/2.
\endsplit
$$

Moreover, one has formulae \thetag{4.23} and \thetag{4.24} in terms of $X$
for the canonical primitive $nef$ element $\wth$ of $Y$ defined by
$(-a,0,b)\mod \bz v$.

These conditions are necessary to have $Y\cong X$ if $\rho (X)\le
2$ and $X$ is a general K3 surface with its Picard lattice,
i. e. the automorphism group of the transcendental periods
$(T(X),H^{2,0}(X))$ is $\pm 1$.
\endproclaim

See Sect. 6 about the cases $\gamma=1$ and $\gamma=2$.

\head
5. Divisorial conditions on moduli of $(X,H)$ which imply $Y\cong X$
and $\gamma (\wH)=\gamma$
\endhead

Using notations above, assuming $\text{g.c.d}(c,d\gamma)=1$ for
$$
\overline{\mu}=\{\mu,-\mu\}\subset
(\bz/(2a_1b_1c^2/\gamma))^\ast,\ \epsilon=\pm 1
$$
we denote by
$$
\Da(r,s,d,\gamma;A)^{\overline{\mu}}_\epsilon \tag{5.1}
$$
the set of all $\delta \in \bn$ such that $\delta\equiv \gamma
\mu^2\mod 4a_1b_1c^2/\gamma$ and the equation $\gamma p_1^2-\delta
q_1^2=\epsilon 2(2/\gamma_2)(a_1/\gamma_a)\gamma_b c$ has an
integral solution $(p_1,q_1)$ satisfying the condition (A)
\thetag{3.3.54}---\thetag{3.3.57} of the $a$-series. Similarly, we
define
$$
\Da(r,s,d,\gamma;B)^{\overline{\mu}}_\epsilon \tag{5.2}
$$
the $b$-series changing $a$ and $b$ places (see the equation
\thetag{3.3.65} and conditions (B)
\thetag{3.3.58}---\thetag{3.3.61}).

We denote
$$
\Da(r,s,d,\gamma)^{\overline{\mu}}= \left(\bigcup_{\epsilon\in
\{-1,1\}} {\Da(r,s,d,\gamma;A)^{\overline{\mu}}_\epsilon}\right)
$$
$$
\bigcup \left(\bigcup_{\epsilon \in \{-1,1\}}
{\Da(r,s,d,\gamma;B)^{\overline{\mu}}_\epsilon}\right). \tag{5.3}
$$
By Theorem 4.1, the set $\Da(r,s,d,\gamma)^{\overline{\mu}}$
describes all possible pairs $H\in N(X)$ of general polarized K3
surfaces $(X,H)$ with $\rk N(X)=2$, the primitive polarization
$\wH=H/d\in N(X)$, the invariant $\gamma (\wH)=\gamma$ (i. e.
$\wH\cdot N(X)=\gamma \bz$) and the invariant $\pm \mu$  for
$\wH\in N(X))$, such that $Y\cong X$. By general results of
\cite{N1} and \cite{N2} the pair $\wH\in N(X)$ defines the
irreducible 18-dimensional moduli of such pairs $(X,\wH)$, i. e. a
(irreducible) divisorial condition on 19-dimensional moduli of
polarized K3 surfaces $(X,H)$ which implies that $\gamma
(\wH)=\gamma$ and $Y\cong X$. Thus, we can interpret our results
as follows.

\proclaim{Theorem 5.1}
The set
$$
\Da(r,s,d,\gamma)= \{(\overline{\mu},\,\delta)\ |\ \{\mu,
-\mu\}\subset \left(\bz/(2a_1b_1c^2/\gamma)\right)^\ast,\ \delta
\in \Da(r,s,d,\gamma)^{\overline{\mu}}\} \tag{5.4}
$$
describes all irreducible divisorial conditions on moduli of
polarized K3 surfaces $(X,H)$ with $H^2=2rs$ and the primitive
polarization $\wH=H/d$ (here $d^2\vert rs$), which imply $Y\cong
X$ for any $X$, and $\wH\cdot N(X)=\gamma \bz$ for a general $X$.

We have (see \thetag{5.3})
$$
\Da(r,s,d,\gamma)^{\overline{\mu}}= \left(\bigcup_{\epsilon\in
\{-1,1\}} {\Da(r,s,d,\gamma;A)^{\overline{\mu}}_\epsilon}\right)
$$
$$
\bigcup \left(\bigcup_{\epsilon\in \{-1,1\}}
{\Da(r,s,d,\gamma;B)^{\overline{\mu}}_\epsilon}\right). \tag{5.5}
$$
where each set $\Da(r,s,d,\gamma;A)^{\overline{\mu}}_\epsilon$ and
$\Da(r,s,d,\gamma;B)^{\overline{\mu}}_\epsilon$ is infinite if it
is not empty.
\endproclaim

\demo{Proof} We need to prove the last statement only.
Assume that $\Da(r,s,d,\gamma;A)^{\overline{\mu}}_\epsilon$ is not
empty. Thus, there exist integral $(p_0,q_0)$ such that
$$
\cases {\gamma p_0^2-\epsilon 2(2/\gamma_2)(a_1/\gamma_a)\gamma_b
c\over q_0^2}\equiv \gamma\mu^2
\mod {{4a_1b_1c^2\over \gamma}}\\
\gamma p_0^2-\epsilon 2(2/\gamma_2)(a_1/\gamma_a)\gamma_b c >0\\
(p_0,q_0)\ \text{satisfies}\ (A)
\endcases .
\tag{5.6}
$$
Then
$$
\delta_0={\gamma p_0^2-\epsilon
2(2/\gamma_2)(a_1/\gamma_a)\gamma_b c \over q_0^2}\in
\Da(r,s,d,\gamma;A)^ {\overline{\mu}}_\epsilon. \tag{5.7}
$$
The \thetag{5.6} is equivalent to
$$
\cases \gamma p_0^2-\epsilon 2(2/\gamma_2)(a_1/\gamma_a)\gamma_b c
\equiv \gamma\mu^2 q_0^2
\mod 4a_1b_1c^2q_0^2/\gamma \\
\gamma p_0^2-\epsilon 2(2/\gamma_2)(a_1/\gamma_a)\gamma_b c>0\\
(p_0,q_0)\ \text{satisfies}\ (A)
\endcases .
\tag{5.8}
$$
Clearly, $(p,q_0)$ where
$$
p\equiv p_0\mod 8a_1b_1c^2q_0^2,\ \text{and}\ \gamma
p_0^2-\epsilon 2(2/\gamma_2)(a_1/\gamma_a)\gamma_b c>0 \tag{5.9}
$$
also satisfies \thetag{5.8} and defines
$$
\delta={\gamma p_0^2-\epsilon 2(2/\gamma_2)(a_1/\gamma_a)\gamma_b
c \over q_0^2}\in \Da(r,s,d,\gamma;A)^ {\overline{\mu}}_\epsilon.
\tag{5.10}
$$
Obviously, their number is infinite. This proves the statement.
\enddemo

The key question is:

\proclaim{Problem 5.2} When $\Da(r,s,d,\gamma;A)^
{\overline{\mu}}_\epsilon$ and $\Da(r,s,d,\gamma;B)^
{\overline{\mu}}_\epsilon$ are non-empty?
\endproclaim

We hope to consider this question in further publications on the
subject. It was shown in  \cite{N4} (see also \cite{MN1},
\cite{MN2}) that at least one of these sets is not empty if $d=1$
and $\gamma=1$. Theorem 4.1 and exactly the same considerations
as in \cite{N4} show that it is also valid for $\gamma=1$ and any
$d$ because singular conditions (AS) and (BS) satisfy if
$\gamma=1$. Thus we have

\proclaim{Theorem 5.3} At least for one of $\overline{\mu}$,
$\epsilon$ one of sets $\Da(r,s,d,\gamma=1;A)^
{\overline{\mu}}_\epsilon$ or \newline $\Da(r,s,d,\gamma=1;B)^
{\overline{\mu}}_\epsilon$ is not empty.

In particular, for any primitive isotropic Mukai vector $(r,H,s)$
the set of divisorial conditions on moduli of $X$ which imply that
$Y\cong X$ and $\gamma=1$ is not empty and is then infinite.
\endproclaim

We hope to consider Problem 5.2 for other $\gamma$ in further
publications on the subject.

\head
6. Examples of $\gamma=1$ and $\gamma=2$
\endhead
As concrete examples of results of Sect. 4, we consider cases of
$\gamma=1$ and $\gamma=2$.

When $\gamma=1$, then singular conditions (AS) and (BS) are obviously
valid, and we obtain especially simple results.
We formulate only the analogy of Theorem 4.4.

\proclaim{Theorem 6.1} Let $X$ be a K3 surface and $H$ a
polarization of $X$ of degree $H^2=2rs$ where $r,\,s\in \bn$. Assume that
the Mukai vector $(r,H,s)$ is primitive. Let $Y$ be the moduli
space of sheaves on $X$ with the isotropic Mukai vector
$v=(r,H,s)$. Let $\wH=H/d$ be the corresponding primitive polarization.

We have $Y\cong X$ if there exists $\wth_1\in N(X)$ such that
$\wH$, $\wth_1$ belong to a 2-dimensional primitive sublattice
$N\subset N(X)$ such that
$$
\wH\cdot N=\bz
\tag{6.1}
$$
(i. e. $\gamma=1$),
moreover, for one of $\epsilon=\pm 1$ the element $\wth_1$ belongs
to the $a$-series or to the $b$-series described below:

$\wth_1$ belongs to the $a$-series if
$$
\wth_1^2=\epsilon2b_1c\ and\ \wH\cdot \wth_1\equiv
0\mod{b_1c},
\tag{6.2}
$$
$$
\wH\cdot \wth_1\not\equiv 0\mod{b_1cl_1},\
\wth_1/l_2\not\in N(X)
\tag{6.3}
$$
for any prime $l_1$ such that $l_1^2\vert a_1$,
and any prime $l_2$ such that $l_2^2\vert b_1$.

$\wth_1$ belongs to the $b$-series if
$$
\wth_1^2=\epsilon2a_1c\ and\ \wH\cdot \wth_1\equiv
0\mod{a_1 c},
\tag{6.4}
$$
$$
\wH\cdot \wth_1\not\equiv 0\mod{a_1cl_1},\
\wth_1/l_2\not\in N(X) \tag{6.5}
$$
for any prime $l_1$ such that $l_1^2\vert b_1$ and
any prime $l_2$ such that $l_2^2\vert a_1$.

Moreover, one has formulae \thetag{4.23} and \thetag{4.24} in terms of $X$
for the canonical primitive $nef$ element $\wth$ of $Y$ defined by
$(-a,0,b)\mod \bz v$.

These conditions are necessary to have $Y\cong X$ and $\wH\cdot N(X)=\bz$
if $\rho (X)\le 2$ and $X$ is a general K3 surface with its Picard lattice,
i. e. the automorphism group of the transcendental periods
$(T(X),H^{2,0}(X))$ is $\pm 1$.
\endproclaim

This generalizes results of \cite{MN1}, \cite{MN2} and \cite{N4} where
the condition $H\cdot N=\bz$ had been imposed (i. e. $d=\gamma=1$).

\smallpagebreak

Now let us assume that $\gamma=2$. By Theorem 4.4, we obtain three cases
which all together cover all possibilities for $\gamma=2$.

When $\gamma=2$ and $a_1\equiv b_1\equiv 1\mod 2$, then
$\gamma_2=2$, and singular conditions (AS) and (BS) satisfy.
Theorem 4.4 gives then

\proclaim{Theorem 6.2} Let $X$ be a K3 surface and $H$ a
polarization of $X$ of degree $H^2=2rs$ where $r,\,s\in \bn$. Assume that
the Mukai vector $(r,H,s)$ is primitive. Let $Y$ be the moduli
space of sheaves on $X$ with the isotropic Mukai vector
$v=(r,H,s)$. Let $\wH=H/d$ be the corresponding primitive polarization.
Assume that
$$
\text{g.c.d}(2,c)=1 \tag{6.6}
$$
and
$$
a_1\equiv b_1\equiv 1\mod 2. \tag{6.7}
$$

We have $Y\cong X$ if there exists $\wth_1\in N(X)$ such that
$\wH$, $\wth_1$ belong to a 2-dimensional primitive sublattice
$N\subset N(X)$ such that
$$
\wH\cdot N=2 \bz
\tag{6.8}
$$
(then $\gamma=2$, $\gamma_a=\gamma_b=1$ and $\gamma_2=2$),
moreover, for one of $\epsilon=\pm 1$ the element $\wth_1$ belongs
to the $a$-series or to the $b$-series described below:

$\wth_1$ belongs to the $a$-series if
$$
\wth_1^2=\epsilon2b_1c\ and\ \wH\cdot \wth_1\equiv
0\mod{2b_1c},
\tag{6.9}
$$
$$
\wH\cdot \wth_1\not\equiv 0\mod{2b_1cl_1},\
\wth_1/l_2\not\in N(X)
\tag{6.10}
$$
for any prime $l_1$ such that $l_1^2\vert a_1$,
and any prime $l_2$ such that $l_2^2\vert b_1$.

$\wth_1$ belongs to the $b$-series if
$$
\wth_1^2=\epsilon2a_1c\ and\ \wH\cdot \wth_1\equiv
0\mod{2a_1 c},
\tag{6.11}
$$
$$
\wH\cdot \wth_1\not\equiv 0\mod{2a_1cl_1},\
\wth_1/l_2\not\in N(X) \tag{6.12}
$$
for any prime $l_1$ such that $l_1^2\vert b_1$ and
any prime $l_2$ such that $l_2^2\vert a_1$.

Moreover, one has formulae \thetag{4.23} and \thetag{4.24} in terms of $X$
for the canonical primitive $nef$ element $\wth$ of $Y$ defined by
$(-a,0,b)\mod \bz v$.

These conditions are necessary (for odd $c$, $a_1$ and $b_1$) to have
$Y\cong X$ and $\wH\cdot N(X)=2\bz$
if $\rho (X)\le 2$ and $X$ is a general K3 surface with its Picard lattice,
i. e. the automorphism group of the transcendental periods
$(T(X),H^{2,0}(X))$ is $\pm 1$.
\endproclaim

In \cite{MN2} the primitive isotropic Mukai vector $(c,H,c)$ where $H^2=2c^2$
had been considered. Then $a=b=1$, $d=1$, $a_1=b_1=1$ and $\gamma\vert 2$.
The case $\gamma=1$ had been described in \cite{MN2}.
Theorem 6.2 describes the remaining case $\gamma=2$
and then $c$ is odd which was not considered in \cite{MN2}.

Now assume that $\gamma=2$ and $2\vert a_1$. Then $\gamma_2=1$,
$\gamma_a=2$ and $\gamma_b=1$. The singular condition (AS) gives
then \thetag{6.18} and \thetag{6.19} below. The singular condition
(BS) gives \thetag{6.22} and \thetag{6.23} below. Thus, Theorem
4.4 implies the following.

\proclaim{Theorem 6.3} Let $X$ be a K3 surface and $H$ a
polarization of $X$ of degree $H^2=2rs$ where $r,\,s\in \bn$. Assume that
the Mukai vector $(r,H,s)$ is primitive. Let $Y$ be the moduli
space of sheaves on $X$ with the isotropic Mukai vector
$v=(r,H,s)$. Let $\wH=H/d$ be the corresponding primitive polarization.
Assume that
$$
\text{g.c.d}(2,c)=1 \tag{6.13}
$$
and
$$
a_1\equiv 0\mod 2,\ b_1\equiv 1\mod 2. \tag{6.14}
$$

We have $Y\cong X$ if there exists $\wth_1\in N(X)$ such that
$\wH$, $\wth_1$ belong to a 2-dimensional primitive sublattice
$N\subset N(X)$ such that
$$
\wH\cdot N=2 \bz
\tag{6.15}
$$
(then $\gamma=2$, $\gamma_a=2, \gamma_b=1$ and $\gamma_2=1$),
moreover, for one of $\epsilon=\pm 1$ the element $\wth_1$ belongs
to the $a$-series or to the $b$-series described below:

$\wth_1$ belongs to the $a$-series if
$$
\wth_1^2=\epsilon2b_1c\ and\ \wH\cdot \wth_1\equiv
0\mod{2b_1c},
\tag{6.16}
$$
$$
\wH\cdot \wth_1\not\equiv 0\mod{2b_1cl_1},\
\wth_1/l_2\not\in N(X)
\tag{6.17}
$$
for any prime $l_1$ such that $l_1^2\vert a_1$ and $\text{g.c.d}(l_1,2)=1$,
and any prime $l_2$ such that $l_2^2\vert b_1$; moreover (singular conditions),
$$
\wH\cdot \wth_1\not\equiv 0\mod 4b_1c
\tag{6.18}
$$
and
$$
\delta\not\equiv 2\mu^2\mod4a_1\ if\ 4\vert a_1.
\tag{6.19}
$$

$\wth_1$ belongs to the $b$-series if
$$
\wth_1^2=\epsilon2a_1c\ and\ \wH\cdot \wth_1\equiv
0\mod{a_1 c},
\tag{6.20}
$$
$$
\wH\cdot \wth_1\not\equiv 0\mod{a_1cl_1},\
\wth_1/l_2\not\in N(X) \tag{6.21}
$$
for any prime $l_1$ such that $l_1^2\vert b_1$ and
any prime $l_2$ such that $l_2^2\vert a_1$ and $\text{g.c.d}(l_2,2)=1$;
moreover (singular conditions),
$$
\wH\cdot \wth_1\equiv 0\mod {2a_1c}
\tag{6.22}
$$
and
$$
\wth_1/2\not \in N(X).
\tag{6.23}
$$

Moreover, one has formulae \thetag{4.23} and \thetag{4.24} in terms of $X$
for the canonical primitive $nef$ element $\wth$ of $Y$ defined by
$(-a,0,b)\mod \bz v$.

These conditions are necessary (for odd $c$, even $a_1$ and odd $b_1$) to have
$Y\cong X$ and $\wH\cdot N(X)=2\bz$
if $\rho (X)\le 2$ and $X$ is a general K3 surface with its Picard lattice,
i. e. the automorphism group of the transcendental periods
$(T(X),H^{2,0}(X))$ is $\pm 1$.
\endproclaim

Changing $a$ and $b$ places, we get from Theorem 4.4 the remaining
case.

\proclaim{Theorem 6.4} Let $X$ be a K3 surface and $H$ a
polarization of $X$ of degree $H^2=2rs$ where $r,\,s\in \bn$. Assume that
the Mukai vector $(r,H,s)$ is primitive. Let $Y$ be the moduli
space of sheaves on $X$ with the isotropic Mukai vector
$v=(r,H,s)$. Let $\wH=H/d$ be the corresponding primitive polarization.
Assume that
$$
\text{g.c.d}(2,c)=1 \tag{6.24}
$$
and
$$
a_1\equiv 1\mod 2,\ b_1\equiv 0\mod 2. \tag{6.25}
$$

We have $Y\cong X$ if there exists $\wth_1\in N(X)$ such that
$\wH$, $\wth_1$ belong to a 2-dimensional primitive sublattice
$N\subset N(X)$ such that
$$
\wH\cdot N=2 \bz
\tag{6.26}
$$
(then $\gamma=2$, $\gamma_a=1, \gamma_b=2$ and $\gamma_2=1$),
moreover, for one of $\epsilon=\pm 1$ the element $\wth_1$ belongs
to the $a$-series or to the $b$-series described below:

$\wth_1$ belongs to the $a$-series if
$$
\wth_1^2=\epsilon2b_1c\ and\ \wH\cdot \wth_1\equiv
0\mod{b_1 c},
\tag{6.27}
$$
$$
\wH\cdot \wth_1\not\equiv 0\mod{b_1cl_1},\
\wth_1/l_2\not\in N(X) \tag{6.28}
$$
for any prime $l_1$ such that $l_1^2\vert a_1$ and
any prime $l_2$ such that $l_2^2\vert b_1$ and $\text{g.c.d}(l_2,2)=1$;
moreover (singular conditions),
$$
\wH\cdot \wth_1\equiv 0\mod {2b_1c}
\tag{6.29}
$$
and
$$
\wth_1/2\not \in N(X).
\tag{6.30}
$$

$\wth_1$ belongs to the $b$-series if
$$
\wth_1^2=\epsilon2a_1c\ and\ \wH\cdot \wth_1\equiv
0\mod{2a_1c},
\tag{6.31}
$$
$$
\wH\cdot \wth_1\not\equiv 0\mod{2a_1cl_1},\
\wth_1/l_2\not\in N(X)
\tag{6.32}
$$
for any prime $l_1$ such that $l_1^2\vert b_1$ and $\text{g.c.d}(l_1,2)=1$,
and any prime $l_2$ such that $l_2^2\vert a_1$; moreover (singular conditions),
$$
\wH\cdot \wth_1\not\equiv 0\mod 4a_1c
\tag{6.33}
$$
and
$$
\delta\not\equiv 2\mu^2\mod4b_1\ if\ 4\vert b_1.
\tag{6.34}
$$

Moreover, one has formulae \thetag{4.23} and \thetag{4.24} in terms of $X$
for the canonical primitive $nef$ element $\wth$ of $Y$ defined by
$(-a,0,b)\mod \bz v$.

These conditions are necessary (for odd $c$, odd $a_1$ and even $b_1$)
to have $Y\cong X$ and $\wH\cdot N(X)=2\bz$
if $\rho (X)\le 2$ and $X$ is a general K3 surface with its Picard lattice,
i. e. the automorphism group of the transcendental periods
$(T(X),H^{2,0}(X))$ is $\pm 1$.
\endproclaim

Theorems 6.2 --- 6.4 cover all types of a primitive isotropic
Mukai vector when it is in principle possible to have $Y\cong X$
and $\gamma=2$.

\smallpagebreak

Using results of Sect. 4, one can write down similar very concrete
and effective results for any $\gamma$.

\enddocument

\end